\documentclass{article}
\usepackage{amsmath}                     
\usepackage{amssymb}                     
\usepackage{amsthm}                     
\newtheorem{theorem}{Theorem}[section]
\newtheorem{proposition}[theorem]{Proposition}
\newtheorem{lemma}[theorem]{Lemma}
\newtheorem{corollary}[theorem]{Corollary}

\theoremstyle{remark}
\newtheorem{example}[theorem]{Example}
\newtheorem{remark}[theorem]{Remark}
\input xy
\xyoption{all}
\numberwithin{equation}{section}
\newcommand{\RR}{\mathbb{R}}
\newcommand{\ZZ}{\mathbb{Z}}
\newcommand{\CC}{\mathbb{C}}
\newcommand{\FF}{\mathbb{F}}
\newcommand{\HH}{\mathbb{H}}
\newcommand{\CL}[1]{\ensuremath{\operatorname{Cl}(#1)}}
\newcommand{\CCL}[1]{\ensuremath{\mathbb{C}\operatorname{l}(#1)}}
\newcommand{\CLS}[1]{\ensuremath{\operatorname{Cl}^*(#1)}}
\newcommand{\SP}[1]{\operatorname{Spin}(#1)}
\newcommand{\SO}[1]{\operatorname{SO}(#1)}
\newcommand{\so}[1]{\mathfrak{so}(#1)}
\newcommand{\spin}[1]{\mathfrak{spin}(#1)}
\newcommand{\Ad}{\ensuremath{\operatorname{Ad}}}
\newcommand{\ad}{\ensuremath{\operatorname{ad}}}
\newcommand{\tg}{\mathfrak{t}}
\newcommand{\End}{\operatorname{End}}
\newcommand{\GL}{\operatorname{GL}}
\newcommand{\M}[1]{\operatorname{M}_{#1}}
\newcommand{\hg}{\mathfrak{h}}
\newcommand{\g}{\mathfrak{g}}
\newcommand{\pg}{\mathfrak{p}}
\newcommand{\sign}{\operatorname{sgn}}
\newcommand{\trace}{\operatorname{Tr}}
\newcommand{\rs}{\tilde{r}_*}
\newcommand{\lr}{\tilde{l}}
\newcommand{\Diff}{\operatorname{Diff}}
\newcommand{\Hom}{\operatorname{Hom}}
\newcommand{\ch}{\operatorname{ch}}
\title{\vspace{5cm}\textbf{The Dirac operator on compact symmetric spaces}}
\author{\textbf{Emiko Dupont}\vspace{4cm}\\
Master thesis for\\ 
the Cand.Scient degree in mathematics\\
at the University of Copenhagen\\
Advisor: Henrik Schlichtkrull}
\date{\textbf{June 2003}}
\begin{document}
\maketitle
\thispagestyle{empty}
\newpage
\thispagestyle{empty}
\phantom{1}

\newpage
\thispagestyle{empty}
\section*{Acknowledgements}
Firstly, I wish to thank my advisor, Henrik Schlichtkrull, for his guidance
throughout my work with this thesis and for introducing me to the field of
representation theory in which I have benefitted greatly from his
expertise. In addition, I am grateful for his advice and support
during my application for a Ph.D. position.

I also wish to thank the mathematics department at the University of
Utah where part of this work was carried out. In particular, I am
grateful to Henryk Hecht for his generous efforts and guidance, and to
Dragan Mili\v{c}i\'{c} for his advice. I have also benefitted from my stay at
the University of Oxford during which I was introduced to the field of
differential geometry. In particular, I would like to thank Ulrike
Tillmann and Andrew Dancer.

Finally, I am grateful to my family and friends for their support and
encouragement. Special thanks to my friends at the University of Copenhagen for creating a pleasant and
friendly working environment throughout my studies, and to my
parents, Johan Dupont and Mariko Hayashi, my sister, Yoko Dupont, and my fianc\'{e}, Alastair Craw,
for putting up with me, especially in recent months.
\newpage
\thispagestyle{empty}
\phantom{1}

\newpage
\thispagestyle{empty}
\tableofcontents
\newpage
\thispagestyle{empty}
\phantom{1}

\newpage
\begin{abstract}
Let $G$ be a compact connected semisimple Lie group and let $H\subset G$ be a
closed connected subgroup such that
$\operatorname{rank}G=\operatorname{rank}H$ and $G/H$ is a symmetric
space.  Given an irreducible representation of $H$, we define a Dirac
operator $D$ and determine the representations of $G$ in the kernel of $D$. Moreover, we show
that any
irreducible representation of $G$ can be constructed in this way. Our
approach is 
similar to that of Parthasarathy~\cite{Pa}. 
\end{abstract}

\section{Introduction}
Let $G$ be a compact connected semisimple Lie group and let $H\subset G$ be a
closed connected subgroup such that
$\operatorname{rank}G=\operatorname{rank}H$ and $G/H$ is a symmetric
space. For each irreducible representation $(V,\tau)$ of $H$ we define a
Dirac operator $D$ on the sections of a vector bundle
$\underline{S\otimes V}$ on
$G/H$. The operator $D$, which depends on the representation $(V,\tau)$, splits into two parts $D^+$ and
$D^-$, and the kernels of $D^{\pm}$ are finite-dimensional
representations $\tilde{\pi}^{\pm}$ of $G$ under the left regular
action. The objective of this paper is to show that one of the
representations 
$\tilde{\pi}^{\pm}$ is an irreducible
representation of $G$ while the other is zero (except in some cases
when both are zero). Moreover, every
irreducible representation of $G$ can be constructed as either
$\tilde{\pi}^+$ or $\tilde{\pi}^-$ for some irreducible representation
$(V,\tau)$ of $H$. This gives a geometric realization of the representations of
$G$ in terms of representations of $H$.

We now describe the contents of the paper in more detail. 
In sections \ref{chap:spinor}-\ref{chap:dirac} we review the theory necessary to  define
the Dirac operator and study its properties. Given an $n$-dimensional
real vector space, one defines a representation $(S,\sigma)$
of the spin group $\SP{n}$ called the spinor representation
 which we investigate in section \ref{chap:spinor}. Then in
section \ref{chap:spin} we look at (not necessarily
compact or symmetric) homogeneous spaces $G/H$ . In certain cases, the spinor
representation induces a representation $(S,\chi)$ of $H$ which we
also call the spinor representation, and we obtain an induced vector
bundle $\underline{S}$ on $G/H$ called the spinor bundle. Thus if
$(V,\tau)$ is any representation on $H$ we get a vector bundle
$\underline{S\otimes V}$ on $G/H$. Under additional assumptions on $G$
and $H$ we also have the half spinor representations $(S^{\pm},\chi^{\pm})$ and
induced bundles $\underline{S^{\pm}\otimes V}$ on $G/H$. In section
\ref{chap:dirac} we define the Dirac operator $D$ and study some of
its properties. By definition, $D$ is an operator on the space
of sections $\Gamma(\underline{S\otimes V})$ of $\underline{S\otimes
  V}$. Since $D$ is an elliptic first order differential operator which
commutes with the left regular action of $G$ on
$\Gamma(\underline{S\otimes V})$, it follows that the kernel of the Dirac
operator is a finite-dimensional representation $\tilde{\pi}$ of $G$
under the left regular action. If the half spinor representations
exist, we also have the
operators $D^{\pm}$ on $\Gamma(\underline{S^{\pm}\otimes V})$ and we
obtain finite-dimensional representations $\tilde{\pi}^{\pm}$ of $G$ on
the kernel of $D^{\pm}$.

In sections \ref{chap:symm}-\ref{chap:ex} we specialize to the case of compact symmetric spaces
$G/H$ where $\operatorname{rank}G=\operatorname{rank}H$. By
finding the irreducible parts of the spinor representation $\chi$
of $H$ in section \ref{chap:symm}, we determine
the action of the Casimir element of $H$ which in turn is an
important step in determining the square of the Dirac operator in
section \ref{chap:dirac}. The square of the Dirac operator consists of a constant term plus the action
of the Casimir element of $G$ under the left regular action. Thus elements of the kernel of $D$ are eigenvectors
of the Casimir operator and this gives a criterion which must be satisfied by the subrepresentations of
$\tilde{\pi}^{\pm}$. As a consequence, the
irreducible parts of $\tilde{\pi}^{\pm}$ have multiplicity at most one and, furthermore,
$\tilde{\pi}^+$ and
$\tilde{\pi}^-$ have no common irreducible part. On the other hand, by
calculating the difference between the trace of $\tilde{\pi}^+$ and
that of  $\tilde{\pi}^-$, we see that
$\tilde{\pi}^+$ and
$\tilde{\pi}^-$ differ in at most one irreducible part. Hence one of
$\tilde{\pi}^+$ or $\tilde{\pi}^-$ is irreducible and the other is zero
(unless both are zero). The main result, theorem
\ref{thm:repr}, then follows. Finally, in section \ref{chap:ex}, we study the specific case
$\SO{2m+1}/\SO{2m}$.

For $G/H$ a symmetric space where $G$ is 
non-compact and $H$ a maximally compact subgroup, 
Parthasarathy~\cite{Pa} uses
representations on the kernels of Dirac
operators $D^{\pm}$ to construct the so-called discrete series of $G$, i.e., the
representations of $G$ which are equivalent to the left or right
regular representation of $G$ on a closed invariant subspace of
$L_2(G)$. The
proof of our main result, theorem \ref{thm:repr}, follows an approach similar to that of \cite{Pa} but is simplified by the compactness
assumption. More recent work by Slebarski~\cite{Sl} considers the case where $G$ is compact and $H$ is a
maximal torus of $G$, i.e., $G/H$ need not be symmetric. A result similar to theorem \ref{thm:repr}
holds but Slebarski's proof is complicated by the presence of an extra
term in the square
of the Dirac operator. This extra term depends on
the torsion of the reductive connection; in the symmetric case the
reductive connection is torsion-free. 

\newpage
\section{The spinor representation}\label{chap:spinor}
\subsection{Clifford algebras}
In this section we recall some properties of Clifford
algebras. There is a classification of Clifford
algebras which shows that they are in fact familiar matrix algebras or
a direct sum of two matrix algebras. This simplifies their
representation theory and it shows that Clifford algebras can be
thought of as Lie algebras. Studying Clifford algebras will enable us
to define the spinor representation of the spin group in
section \ref{sec:spingroup}. The exposition is mainly based on Lawson
and Michelsohn~\cite{LM}
and Gilbert and Murray~\cite{GM}.
 
Let $V$ be an $n$-dimensional vector space over $\FF$ (with $\FF=\RR$ or
$\CC$) with a quadratic
form $q$. Let $\CL{V,q}$ denote the (universal) Clifford algebra associated to
$(V,q)$. This is an associative algebra with unit which is generated
by the identity $1$ and $V$ subject to the relation
\[
v\cdot v+q(v)1=0\quad\textrm{for all }v\in V. 
\]
The Clifford algebra has the following universal property.
\begin{proposition}\label{prop:extend}
Let $f\colon V\rightarrow A$ be a linear map into an associative algebra
over $\FF$ with
unit such that
\[
f(v)\cdot f(v)+q(v)1=0\quad\textrm{for all }v\in V.
\]
Then $f$ extends uniquely to an algebra homomorphism
$f\colon \CL{V,q}\rightarrow A$. $\CL{V,q}$ is the unique associative algebra over
$\FF$ with
this property.
\end{proposition}
\proof See proposition I.1.1 \cite{LM}. \qed

Let $B$ be the bilinear form associated to $q$, i.e., 
\[
B(v,w)=\tfrac{1}{2}\left( q(v)+q(w)-q(v-w)\right)\quad\textrm{for all }v,w\in V.
\]
We say that $(V,q)$ is non-degenerate if
\[
V^{\perp}=\{w\in V\mid B(v,w)=0 \textrm{ for all }v\in V\}=\{0\}.
\]
In the following we assume that $V$ is non-degenerate.
Let $\{e_1,\ldots,e_n\}$ be an orthonormal basis of $V$ in the sense
that\[
B(e_i,e_j)=0\quad\textrm{for }i\ne j\]
\[
q(e_i)\in\{\pm 1\}.
\]
If $\FF=\CC$ we may assume that $q(e_i)=1$.
The elements of the orthonormal basis satisfy the
relations
\[
e_je_k+e_ke_j=-2q(e_j)\delta_{jk}
\]
and the products
\[
e_1^{m_1}\cdots e_n^{m_n},\quad m_j\in\{0,1\}
\]
form a basis of $\CL{V,q}$ (where $e_1^0\cdots e_n^0=1$). In particular, $\CL{V,q}$
has dimension $2^n$.

If $(V,q)$ is an $n$-dimensional quadratic space over $\RR$, we can
choose a basis of $V\cong\RR^n$ such that $q$ is of the form 
\[
q_{r,s}(x)=x_1^2+\cdots+x_r^2-x_{r+1}^2-\cdots -x_n^2\quad\textrm{for }x=(x_1,\ldots,x_n)\in\RR^n
\]
where $r+s=n$.   
We denote the Clifford algebra of $(\RR^{r+s},q_{r,s})$ by $\CL{r,s}$. In the special
case of $s=0$, the associated bilinear form $B$ is just the usual
inner product on $\RR^n$ and we denote the Clifford algebra $\CL{n,0}$ by $\CL{n}$.

Similarly if $(V,q)$ is an $n$-dimensional quadratic space over $\CC$, we can
choose a basis of $V\cong\CC^n$ such that $q$ is of the form 
\[
q_{\CC}(z)=z_1^2+\cdots+z_n^2\quad\textrm{for }z=(z_1,\ldots,z_n)\in\CC^n. 
\]
We denote the Clifford algebra of $(\CC^n,q_{\CC})$ by $\CCL{n}$. Now consider
the complexification $\CL{r,s}\otimes_{\RR}\CC$ of $\CL{r,s}$. By
proposition \ref{prop:extend} this is just the Clifford algebra of the
complex quadratic space $(\CC^n,q_{r,s}\otimes 1)$ for $r+s=n$
(where we think of $\CC^n$ as $\RR^{r+s}\otimes_{\RR}\CC$) and so
\[
\CCL{n}\cong\CL{r,s}\otimes_{\RR}\CC\quad\textrm{for }r+s=n.
\]

We claim that $\CL{V,q}$ is a $\ZZ_2$-graded algebra. To see this, let
$\overline{\alpha}\colon V\rightarrow V$ be the linear map given by
$\overline{\alpha}(v)=-v$. This extends by proposition \ref{prop:extend} to a map
$\overline{\alpha}\colon \CL{V,q}\rightarrow \CL{V,q}$. Since $\overline{\alpha}^2=\operatorname{Id}$, we get a splitting
of $\CL{V,q}$ into the $+1$ and $-1$ eigenspaces $\CL{V,q}^0$ and $\CL{V,q}^1$
which we call the even part and odd part of $\CL{V,q}$
respectively. Note that $\CL{V,q}^0$ is spanned by the even products
$\{e_1\cdots e_{2k}\mid 2k\le n\}$ and $\CL{V,q}^1$ is spanned by the odd
products $\{e_1\cdots e_{2k-1}\mid 2k-1\le n\}$. This splitting is a $\ZZ_2$-grading of $\CL{V,q}$ in the sense that
\[
\CL{V,q}^0\CL{V,q}^0=\CL{V,q}^1\CL{V,q}^1=\CL{V,q}^0
\]
\[
\CL{V,q}^1\CL{V,q}^0=\CL{V,q}^0\CL{V,q}^1=\CL{V,q}^1.
\]
Note that $\CL{V,q}^0$ is a subalgebra of $\CL{V,q}$. In the case of $\CL{n}$,
we get the following result:
\begin{proposition}\label{prop:lower}
There is an algebra isomorphism
\[
\CL{n}\cong\CL{n+1}^0.
\]
It follows that also
\[
\CCL{n}\cong\CCL{n+1}^0.
\]
\end{proposition}
\proof Let $f\colon \RR^n\rightarrow\CL{n+1}^0$ be given by
\[
f(e_i)=e_i e_{n+1}\quad\textrm{for }i\in\{1,\ldots n\}.
\]
This extends by proposition \ref{prop:extend} to the desired algebra
isomorphism $\CL{n}\rightarrow\CL{n+1}^0$. (See theorem I.3.7 of
\cite{LM}). \qed

It turns out that the Clifford algebras $\CL{r,s}$ and $\CCL{n}$ are
familiar matrix algebras over $\RR,\CC$ or $\HH$ (considered as
algebras over $\RR$). We now focus
our attention on the algebras $\CL{n}$
and $\CCL{n}$ which are classified in the following theorem. (A similar
classification of $\CL{r,s}$ can be made).
\begin{theorem}\label{thm:class}
For $n\ge 0$ we have that
\begin{eqnarray*}
\CL{n+8}&\cong&\CL{n}\otimes_{\RR}\CL{8}\\
\CCL{n+2}&\cong&\CCL{n}\otimes_{\CC}\CCL{2}.
\end{eqnarray*}
Hence we get a classification of the algebras $\CL{n}$ and $\CCL{n}$ due to
the following table:
\begin{center}
\begin{tabular}{|c|c|c|}
\hline
&$\CL{n}$&$\CCL{n}$\\
\hline
$1$&$\CC$&$\CC\oplus\CC$\\
\hline
$2$&$\HH$&$\M{2}(\CC)$\\
\hline
$3$&$\HH\oplus\HH$&$\M{2}(\CC)\oplus
\M{2}(\CC)$\\
\hline
$4$&$\M{2}(\HH)$&$\M{4}(\CC)$\\
\hline
$5$&$\M{4}(\CC)$&$\M{4}(\CC)\oplus
\M{4}(\CC)$\\
\hline
$6$&$\M{8}(\RR)$&$\M{8}(\CC)$\\
\hline
$7$&$\M{8}(\RR)\oplus
\M{8}(\RR)$&$\M{8}(\CC)\oplus \M{8}(\CC)$\\
\hline
$8$&$\M{16}(\RR)$&$\M{16}(\CC)$\\
\hline
\end{tabular}
\end{center}

\end{theorem}
\proof See theorem I.4.3 of \cite{LM}. \qed

The algebras $\CL{n}$ and $\CCL{n}$ for $n\ge 3$ are
shown to be isomorphic to $\M{k}(\FF)$ or $\M{k}(\FF)\oplus\M{k}(\FF)$
(with $\FF=\RR,\CC$ or $\HH$) by using the following isomorphisms.
\begin{eqnarray*}
\M{l}(\RR)\otimes_{\RR}\M{m}(\RR)&\cong&\M{lm}(\RR)\\
\M{l}(\RR)\otimes_{\RR}\FF&\cong&\M{l}(\FF)\\
\M{l}(\CC)\otimes_{\CC}\M{m}(\CC)&\cong&\M{lm}(\CC)
\end{eqnarray*}

The classification shows that $\CL{n}$ splits into a sum of two
matrix algebras exactly when $n\equiv 3\mod 4$ and similarly
$\CCL{n}$ splits into a sum of two matrix algebras exactly when $n$ is odd. This splitting can be understood by
considering the volume element $\omega$ given by
\[
\omega=e_1\cdots e_n
\]
where $\{e_1,\ldots,e_n\}$ is an orthonormal basis of $\RR^n$. Note that $\omega$ is independent of the choice of basis once we have
fixed an orientation of $\RR^n$ and $e_1,\ldots,e_n$ are chosen to be positively oriented.
For $n$ odd, define the complex volume element $\omega_{\CC}$ by
\[
\omega_{\CC}=i^{\tfrac{n+1}{2}}e_1\cdots e_n.
\]
Suppose that $n\equiv 3\mod 4$. Then
\[
\omega^2=(-1)^{(n-1)+(n-2)+\cdots +1}e_1^2\cdots e_n^2=(-1)^{\tfrac{n(n-1)}{2}}(-1)^n=(-1)^{\tfrac{n(n+1)}{2}}=1
\]
and
\[
e_i\omega=(-1)^{i-1}e_1\cdot e_ie_i\cdots
e_n=(-1)^{i-1}(-1)^{n-i}e_1\cdots e_n e_i=\omega e_i\quad\textrm{for }i\in\{1,\ldots,n\}.
\]
Hence, if $\CL{n}^+$ and $\CL{n}^-$ denote the $+1$ and $-1$ eigenspaces
of $\CL{n}$ under multiplication by $\omega$, then
\[
\CL{n}=\CL{n}^+\oplus\CL{n}^-
\]
and $\CL{n}^{\pm}$ are ideals since if $x^{\pm}\in\CL{n}^{\pm}$ then
\[
\omega e_i x^{\pm}=e_i\omega x^{\pm}=\pm e_i x^{\pm}\quad\textrm{for }i\in\{1,\ldots,n\}.
\]
Since $\omega\in \CL{n}^1$, we get that
\[
\overline{\alpha}(\CL{n}^{\pm})=\CL{n}^{\mp}
\]
and so $\CL{n}^+$ and $\CL{n}^-$ are isomorphic.
Similarly, if $n$ is odd, $\omega_{\CC}^2=1$ and $\omega_{\CC}$ lies
in the center of $\CCL{n}$ and therefore we get a splitting into ideals
\[
\CCL{n}=\CCL{n}^+\oplus\CCL{n}^-
\]
where $\CCL{n}^{\pm}$ denote the $\pm 1$ eigenspaces under multiplication
by $\omega_{\CC}$. These splittings correspond to the splittings given
in theorem \ref{thm:class}.

The classification in theorem \ref{thm:class} shows that the representation theory of $\CL{n}$ and
$\CCL{n}$ is particularly simple, since up to equivalence, $\M{k}(\FF)$ (with $\FF=\RR,\CC$
or $\HH$) only has one irreducible representation over $\RR$, namely
the standard representation
$\rho\colon \M{k}(\FF)\rightarrow\End_{\FF}(\FF^k)$, and similarly
$\M{k}(\FF)\oplus \M{k}(\FF)$
only has two irreducible representations, namely $\rho_1$ and $\rho_2$
where for $(\varphi_1,\varphi_2)\in \M{k}(\FF)\oplus \M{k}(\FF)$
\[
\rho_i(\varphi_1,\varphi_2)=\rho(\varphi_i)\quad\textrm{for }i=1,2.
\]
Note that a representation over $\CC$ of a complex matrix algebra is any representation over $\RR$
which commutes with multiplication by $i$.
We therefore get the following result:
\begin{proposition}\label{prop:rep}
When $n\equiv 0,1,2\mod 4$, then up to equivalence there is exactly one
irreducible real representation $\sigma$ of $\CL{n}$.

When $n\equiv 3\mod 4$, then up to equivalence there are exactly two
irreducible real representations $\sigma_1$ and
$\sigma_2$ of $\CL{n}$ where if
$(\varphi_1,\varphi_2)\in\CL{n}^+\oplus\CL{n}^-$ then
\[
\sigma_i(\varphi_1,\varphi_2)=\sigma(\varphi_i)\quad\textrm{for }i=1,2
\]
where $\sigma$ is an irreducible representation of $\CL{n}^+\cong\CL{n}^{-}$.

When $n$ is even, then up to equivalence there is exactly one
irreducible complex representation $\sigma$ of $\CCL{n}$.

When $n$ is odd, then up to equivalence there are exactly two
irreducible complex representations $\sigma_1$ and
$\sigma_2$ of $\CCL{n}$ where if
$(\varphi_1,\varphi_2)\in\CCL{n}^+\oplus\CCL{n}^-$ then
\[
\sigma_i(\varphi_1,\varphi_2)=\sigma(\varphi_i)\quad\textrm{for }i=1,2
\]
where $\sigma$ is an irreducible representation of $\CCL{n}^+\cong\CCL{n}^{-}$.  
\end{proposition}

Since $\CL{n}$ is either of the form $\M{k}(\FF)$ or a direct sum
$\M{k}(\FF)\oplus \M{k}(\FF)$, we can think of $\CL{n}$ as a Lie algebra with
Lie bracket given by
\[
[x,y]=xy-yx\quad\textrm{for all }x,y\in\CL{n}.
\]
Let
\[
\CLS{n}=\{x\in\CL{n}\mid x\textrm{ is invertible}\}
\]
be the multiplicative group of units in $\CL{n}$. Then $\CLS{n}$ will be a
Lie group with Lie algebra $\CL{n}$. The usual exponential mapping $\exp\colon \CL{n}\rightarrow\CLS{n}$ given by
\[
\exp(y)=\sum_{m=0}^{\infty}\frac{1}{m!}y^m
\]
is well-defined and we get the adjoint representation
$\Ad_{\operatorname{Cl}}\colon \CLS{n}\rightarrow\textrm{Aut}(\CL{n})$ with differential
$\ad_{\operatorname{Cl}}\colon \CL{n}\rightarrow\End(\CL{n})$ where
\[
\Ad_{\operatorname{Cl}}(g)x=gxg^{-1}\quad\textrm{for }g\in\CLS{n},x\in\CL{n}
\]
\[
\ad_{\operatorname{Cl}}(x)(y)=[x,y]\quad\textrm{for }x,y\in\CL{n}.
\]
\subsection{The spin group}\label{sec:spingroup}
Now we study the spin group which is a Lie group
that lies in $\CL{n}$. We see that there is an important representation
of the spin group called the spinor representation. We consider some
properties  of the spinor representation and find its weights.

Let
\[
\SP{n}=\{v_1\cdots v_{2k}\in \CL{n}\mid v_j\in V,q(v_j)\in\{\pm
1\}\}.
\]
This is
a compact subgroup of $\CLS{n}$ and is called the spin
group. 
Let
\[
\SO{n}=\{x\in \GL(\RR^n)\mid q(xv)=q(v)\textrm{ for all }v\in\RR^n,\det(x)=1\}.
\]
We then have the following important characterization of
$\SP{n}$.
\begin{theorem}
The map
\[
\psi\colon \SP{n}\rightarrow \SO{n}, \quad \psi(x)v=xvx^{-1}
\]
is a two-fold covering homomorphism.
\end{theorem}
\proof See theorem I.6.3 \cite{GM} or theorem I.2.9 \cite{LM}. \qed

$\SO{n}$ is a connected Lie group for $n\ge 1$ and its fundamental
group is $\pi_1(\SO{n})=\ZZ_2$ for $n\ge 3$ (see Knapp~\cite{Kna}
proposition 1.136). Hence $\SP{n}$ is the universal
covering group of $\SO{n}$ for $n\ge 3$. Since $\SO{2}=S^1$, we have that
$\SP{2}=S^1$ and since $\SO{1}=\{1\}$, $\SP{1}=\{\pm 1\}$. We
have the
following result (see 3.24 Warner~\cite{War}):
\begin{theorem}
There is a unique differentiable structure on $\SP{n}$ such that 
$\psi\colon \SP{n}\rightarrow\SO{n}$ is
differentiable and non-singular. With this differentiable structure,
$\SP{n}$ becomes a Lie group and $\psi$ a Lie group homomorphism.
\end{theorem}
So when we think of $\SP{n}$ as a Lie group in this way, proposition
3.26 of \cite{War} gives us that its Lie
algebra is isomorphic to
\[
\so{n}=\{X\in \M{n}(\RR)\mid X+X^t=0\}.
\] 
Now since $\SP{n}$ is also a closed subgroup of $\CLS{n}$, it inherits a Lie
group structure from $\CLS{n}$. The following theorem shows that this is
in fact the same.
\begin{proposition}\label{prop:spin}
As a Lie subgroup of $\CLS{n}$, $\SP{n}$ has Lie algebra $\spin{n}$
given by the bivectors in $\CL{n}$:
\[
\spin{n}=\operatorname{span}_{\RR}\{e_ie_j\mid i<j\}\subset \CL{n}.
\]
With this differentiable structure, $\psi\colon \SP{n}\rightarrow\SO{n}$ is
differentiable and non-singular. Furthermore,
\[
\ad_{\operatorname{Cl}}(\spin{n})\subset\End(\RR^n)
\]
and can in this way be thought of as $\so{n}$ (acting on $\RR^n$ as the
differential of the standard representation of $\SO{n}$ on $\RR^n$). 
Thus the following diagram commutes:
\[
\xymatrix{
&\CL{n}\ar[ddl]_{\exp}\ar[r]^{\ad_{\operatorname{Cl}}}&\End(\CL{n})\\
&\spin{n}\ar[d]^{\exp}\ar@{^{(}->}[u]\ar[r]^{\psi_*}_{\cong}&\so{n}\ar@{^{(}->}[u]\ar[d]^{\exp}\\
\CLS{n}&\SP{n}\ar@{_{(}->}[l]\ar[r]^{\psi}&\SO{n}
}
\]
\end{proposition}
\proof See theorem I.8.10 of \cite{GM}, lemma 1.1. of \cite{Pa}. \qed

Now we consider the representation theory of $\SP{n}$. Some
obvious representations of $\SP{n}$ are the ones induced by the
map $\psi\colon \SP{n}\rightarrow \SO{n}$. All such representations must be trivial
on the element $-1$. However, using the fact that
$\SP{n}\subset\CL{n}\subset\CCL{n}$, we obtain a representation of
$\SP{n}$ which is non-trivial on $-1$. Proposition \ref{prop:rep}
showed that there are at most two inequivalent irreducible
representations of $\CCL{n}$ over $\CC$. Take any such representation and restrict it
to $\SP{n}\subset\CL{n}$. We see (proposition \ref{prop:sigma} below) that up to equivalence, this gives
exactly one representation of $\SP{n}$ and this is called the (complex) spinor
representation of $\SP{n}$ which we denote by
$\sigma\colon \SP{n}\rightarrow\GL_{\CC}(S)$. The elements of $S$ are called
spinors. We
note that $\sigma(-1)=-\operatorname{Id}$.
Because of the classification in theorem
\ref{thm:class}, it is clear that $S$ has complex dimension 
\[
\dim_{\CC}S=2^m\quad\textrm{for }n=2m, n=2m+1.
\]
\begin{proposition}\label{prop:sigma}
Let $\sigma\colon \SP{n}\rightarrow\GL_{\CC}(S)$ be the complex
spinor representation of $\SP{n}$. 

When $n$ is odd, $\sigma$ is independent of the choice of irreducible
representation of $\CCL{n}$ used to define $\sigma$ and furthermore
$\sigma$ is an irreducible representation of $\SP{n}$. 

When $n=2m$ is even, $\sigma$ splits into a direct
sum of two inequivalent irreducible complex representations
$\sigma^{\pm}\colon \SP{n}\rightarrow \GL_{\CC}(S^{\pm})$ where
$\dim_{\CC}S^{\pm}=2^{m-1}$.
\end{proposition}
\proof
(Similar to that of proposition I.5.12 of \cite{LM}).
We note that due to proposition \ref{prop:lower},
\[
\SP{n}\subset\CCL{n}^0\cong\CCL{n-1},
\]
and since $\SP{n}$ contains an additive basis of $\CCL{n}^0$, any
irreducible representation of $\CCL{n}^0$ restricts to an irreducible
representation of $\SP{n}$.

Suppose that $n=2m+1$ is odd. We saw that in this case
$\CCL{n}=\CCL{n}^+\oplus\CCL{n}^-$ where
$\overline{\alpha}(\CCL{n}^{\pm})=\CCL{n}^{\mp}$. We must have that
\[
\CCL{n}^0=\{(\varphi,\overline{\alpha}(\varphi))\in\CCL{n}^+\oplus\CCL{n}^-\mid \varphi\in\CCL{n}^+\}
\]
and therefore the two inequivalent irreducible representations of $\CCL{n}$ restricted to
$\CCL{n}^0$ become $\sigma_1$ and
$\sigma_2=\sigma_1\circ \overline{\alpha}$. Hence, as representations of
$\CCL{n}^0$, they are equivalent. This gives a representation of $\CCL{n}^0\cong\CCL{2m}\cong\M{2^m}(\CC)$ of complex
dimension $2^m$ and therefore it must be irreducible.

Suppose that $n$ is even. Let $\omega_{\CC}=i^{\tfrac{n}{2}}e_1\cdots e_{n-1}$ be the volume element of
$\CCL{n-1}$. Then under the isomorphism $\CCL{n}^0\cong\CCL{n-1}$,
$\omega_{\CC}$ is taken to the element
\begin{eqnarray}
\omega_{\CC}'&=&i^{\tfrac{n}{2}}(e_1e_n)\cdots(e_{n-1}e_n)=i^{\tfrac{n}{2}}(-1)^{1+\cdots +(n-2)}e_1\cdots
e_{n-1}e_n^{n-1}\nonumber\\
&=&i^{\tfrac{n}{2}}(-1)^{\tfrac{(n-1)(n-2)}{2}+\tfrac{n-2}{2}}e_1\cdots
e_n = i^{\tfrac{n}{2}}e_1\cdots e_n\nonumber
\end{eqnarray}
and so $(\omega'_{\CC})^2=1$ and $\omega'_{\CC}$ commutes with the
elements of $\CCL{n}^0$. Hence if $S^{\pm}$ denote the $\pm 1$
eigenspaces of $\sigma(\omega'_{\CC})$, $\sigma$ restricts to
representations $\sigma^{\pm}$ of $\CCL{n}^0$ on $S^{\pm}$ that correspond to the two inequivalent
representations of $\CCL{n-1}$. \qed

The representations $\sigma^{\pm}$ are called the half spinor
representations of $\SP{n}$.

\begin{proposition}\label{prop:sigmas}
Let $\sigma_*\colon \spin{n}\rightarrow \End(S)$ denote the
differential of the spinor representation. Then the following diagram commutes.
\[
\xymatrix{
\spin{n}\ar@{^(->}[d]\ar[r]^{\sigma_*}&\End(S)\\
\CCL{n}\ar[r]^{\sigma}&\End(S)\ar[u]^{\operatorname{Id}}\\
\SP{n}\ar[r]^{\sigma}\ar@{_(->}[u]&\GL(S)\ar@{_(->}[u]
}
\]
where $\operatorname{Id}$ denotes the identification
\[
v=\tfrac{d}{dt}\arrowvert_{t=0}\exp(tv)\quad\textrm{for }v\in\End(S).
\]
\end{proposition}
\proof See corollary 1.1 \cite{Pa}. \qed

Assume that $n=2m$ or $n=2m+1$. We now find the weights of the complex spinor
representation $\sigma\colon \SP{n}\rightarrow
\GL_{\CC}(S)$. In order to do so we consider
$\sigma_*\colon \spin{n}\rightarrow\End_{\CC}(S)$ which we can also
denote by $\sigma$ due to proposition \ref{prop:sigmas}. For each
$k\in\{1,\ldots,m\}$ define
\[
\omega_k=-ie_{2k-1}e_{2k}
\]
where $\{e_1,\ldots,e_n\}$ is an orthonormal basis of $\RR^n$.
It is easy to see that
\begin{eqnarray}
\omega_i\omega_j&=&\omega_j\omega_i\quad\textrm{for }i,j\in \{1,\ldots,m\}\label{commutes}\\
\omega_k^2&=&1\quad\textrm{for }k\in \{1,\ldots,m\}\label{squared}\\
\omega_k e_j&=&\left\{\begin{array}{ll}-e_j\omega_k&\textrm{for }j=2k,
    j=2k-1\\
e_j\omega_k&\textrm{otherwise}\end{array}\right. .\label{omegaej}
\end{eqnarray}
Hence if we define
\[
\tg=\textrm{span}_{\RR}\{i\omega_k\mid k=1,\ldots,m\},
\]
then clearly $\tg$ is an abelian Lie subalgebra of $\spin{n}$ of
dimension $m$. We know that the Cartan subalgebra of $\so{n}\cong\spin{n}$ has
dimension $m$ and therefore $\tg$ is the Cartan subalgebra of
$\spin{n}$.

We now consider the action of the $\omega_k$'s on $S$ under
$\sigma\colon \spin{n}\rightarrow\End(S)$ (where we have extended $\sigma$
complex linearly to $\spin{n}\oplus i\spin{n}$).
Since $\sigma$ is a homomorphism, (\ref{squared}) shows that $S$ splits into a direct sum $S=S_{+}\oplus S_{-}$ where $S_{\pm}$ are
the $\pm 1$ eigenspaces of $\sigma(\omega_1)$. Since
\[
\sigma(e_j)^2=\sigma(e_j^2)=-\operatorname{Id}\quad\textrm{for }j\in\{1,\ldots,n\},
\]
$\sigma(e_j)$ is an isomorphism for all $j$ and (\ref{omegaej}) gives us that
\[
\sigma(e_1)(S_{\pm})=S_{\mp}.
\]
Hence
\[
\dim S_+=\dim S_-=\tfrac{1}{2}\dim S=2^{m-1}.
\]
Since $\sigma(\omega_1)$
commutes with $\sigma(\omega_k)$ for all $k$, $S_{\pm}$ are
invariant under $\sigma(\omega_k)$ for all $k$. We can therefore
repeat this process with $\omega_2$ in place of $\omega_1$, $S_{\pm}$
in place of $S$ and $e_3$ in place of $e_1$. This gives us a
splitting
\[
S=S_{++}\oplus S_{+-}\oplus S_{-+}\oplus S_{--}
\]
where $S_{\pm\pm}$ are simultaneous eigenspaces of $\omega_1$ and
$\omega_2$, each of dimension $2^{m-2}$. Continueing in this fashion
gives us the splitting
\[
S=\bigoplus_{\varepsilon\in\{\pm 1\}^m}S_{\varepsilon}
\]
where for each
$\varepsilon=(\varepsilon_1,\ldots,\varepsilon_m)\in\{\pm 1\}^m$, $S_{\varepsilon}$ is the simultaneous $\varepsilon_k$ eigenspace
of $\omega_k$. Note that $\dim S_{\varepsilon}=1$ for all
$\varepsilon$.

For $n=2m$, we saw in the proof of proposition \ref{prop:sigma} that
$S^{\pm}$ are the ${\pm} 1$ eigenspaces of $\sigma(\omega_{\CC}')$ where
\[
\omega_{\CC}'=i^m e_1\cdots e_n=\omega_1\cdots \omega_m.
\]
For each $s\in S_{\varepsilon}$
\[
\sigma(\omega_{\CC}')s=\varepsilon_1\cdots\varepsilon_ms.
\] 
So if
\begin{eqnarray}\label{e+-}
E^+&=&\{\varepsilon\in\{\pm 1\}^m\mid \varepsilon_k=-1\textrm{ for an even number
of }\varepsilon_k\}\nonumber\\
E^-&=&\{\pm 1\}^m\backslash E^+,
\end{eqnarray}
it is clear that
\[
S^+=\bigoplus_{\varepsilon\in E^+}S_{\varepsilon},\quad S^-=\bigoplus_{\varepsilon\in E^-}S_{\varepsilon}.
\]
This leads us to the following result.
\begin{proposition}\label{prop:sigmaweights}
The weights of the spinor representation $\sigma\colon \SP{n}\rightarrow
\GL_{\CC}(S)$ for $n=2m$, $n=2m+1$ are given by
\[
\lambda_{\varepsilon}=\tfrac{1}{2}\sum_{k=1}^m\varepsilon_k
\eta_k\circ\psi_*,\quad\textrm{for }\varepsilon=(\varepsilon_1,\ldots,\varepsilon_m)\in\{\pm 1\}^m
\]
where $\{\eta_k\mid 1\le k\le m\}$ are the weights of the
standard representation of $\SO{n}$ on $\RR^n$.

For $n=2m$ the weights of the representation $\sigma^+$ are
$\{\lambda_{\varepsilon}\mid \varepsilon\in E^+\}$ and the weights of $\sigma^-$ are
$\{\lambda_{\varepsilon}\mid \varepsilon\in E^-\}$.

The multiplicity of each $\lambda_{\varepsilon}$ is the number of ways
in which $\lambda_{\varepsilon}$ can be written in the above form.
\end{proposition}
\proof Let $\tg_{\CC}$ be the complexification of $\tg$ and consider the
action of $\tg_{\CC}$ on $S$ under $\sigma$. For 
$x=\sum_{k=1}^m x_k \omega_k\in\tg_{\CC}$ we have that
\[
\sigma(x)(s)=\left(\sum_{k=1}^m x_k \varepsilon_k\right)s\quad\textrm{for }s\in
S_{\varepsilon}.
\]
Hence if $\lambda_{\varepsilon}\in\tg_{\CC}^*$ is defined by
\[
\lambda_{\varepsilon}(\sum_{k=1}^m x_k \omega_k)=\sum_{k=1}^m \varepsilon_k x_k
\]
then 
\[
\sigma(x)(s)=\lambda_{\varepsilon}(x)(s)\quad\textrm{for }s\in
S_{\varepsilon}.
\]
So for each $\varepsilon\in\{\pm 1\}^m$, $S_{\varepsilon}$ is a weight space of
$\sigma$ with
weight $\lambda_{\varepsilon}$. It is clear that if $n=2m$, then for $\varepsilon\in
E^{\pm}$, $S_{\varepsilon}$ is a weight space of $\sigma^{\pm}$ with
weight $\lambda_{\varepsilon}$. Since $\dim S_{\varepsilon}=1$, the
multiplicity of each $\lambda_{\varepsilon}$ is the number of ways in
which $\lambda_{\varepsilon}$ can be expressed in the above form.

Now proposition \ref{prop:spin} showed that the differential of the
standard representation of $\SO{n}$ is the representation
$\ad_{\operatorname{Cl}}\circ\psi_*^{-1}\colon \so{n}\rightarrow\End(\RR^n)$ which we extend complex
linearly to get $\ad_{\operatorname{Cl}}\circ\psi_*^{-1}\colon \so{n}\oplus i\so{n}\rightarrow\End(\CC^n)$. We
have that for $j,k\in\{1,\ldots,m\}$
\[
\ad_{\operatorname{Cl}}(\omega_j)(e_{2k-1}\pm ie_{2k})=\left\{\begin{array}{ll}
\mp 2(e_{2k-1}\pm ie_{2k})&\textrm{if }j=k\\
0&\textrm{otherwise}
\end{array}\right. .
\]
So if for $k\in\{1,\ldots,m\}$, $\eta_k\in\tg_{\CC}^*$ is defined as
\[
\eta_k\circ\psi_*(\sum_{j=1}^m x_j \omega_j)=2x_k,
\] 
then for all $x\in\tg_{\CC}$ and $k\in\{1,\ldots,m\}$ we have that
\[
\ad_{\operatorname{Cl}}(x)v=\left\{\begin{array}{ll}
\eta_k(\psi_*x)v&\textrm{if }v\in\operatorname{span}\{e_{2k-1}-ie_{2k}\}\\
-\eta_k(\psi_*x)v&\textrm{if
  }v\in\operatorname{span}\{e_{2k-1}+ie_{2k}\}\end{array}\right. .
\]
In the case $n=2m+1$ we also have that
\[
\ad_{\operatorname{Cl}}(\omega_j)(e_{2m+1})=0\quad\textrm{for } j\in\{1,\ldots,m\}
\]
and therefore
\[
\ad_{\operatorname{Cl}}(x)(v)=0\quad\textrm{for } x\in\tg_{\CC}, v\in\operatorname{span}\{e_{2m+1}\}.
\]
Since $\{e_{2k-1}\pm ie_{2k}\mid 1\le k\le m\}$ is a
basis of $\CC^{2m}$ and $\{e_{2k-1}\pm ie_{2k}\mid
  1\le k\le m\}\cup\{e_{2m+1}\}$ is a basis of $\CC^{2m
+1}$, we get that the weights of the standard
representation of $\SO{n}$ are $\{\pm\eta_k\mid
  1\le k\le m\}$. Since
\[
\lambda_{\varepsilon}=\tfrac{1}{2}\sum_{k=1}^m\varepsilon_k \eta_k\circ\psi_*,
\]
this completes the proof.
\qed

\section{Spin structure on $G/H$}\label{chap:spin}
\subsection{The tangent bundle on $G/H$}
We now turn our attention to homogeneous spaces $G/H$. In
section \ref{sec:spinorb} we see that in some cases it is possible to
define a so-called spin structure of the tangent bundle $T(G/H)$ of
$G/H$. This gives a vector bundle called the spinor bundle on
$G/H$ in which each fibre is isomorphic to the space of
spinors $S$. The construction is described in \cite{Pa} and \cite{LM}.

Let $G$ be a smooth connected real Lie group with Lie algebra
$\g$. We have the adjoint representation of $G$ denoted by
$\Ad\colon G\rightarrow\End(\g)$. Suppose that we have a closed
connected subgroup $H$ of $G$ with Lie algebra $\hg$ such
that $\g=\hg\oplus\pg$ (as vector spaces) where $\pg$ is invariant
under $\Ad(H)$ and there is an inner product $\langle\cdot,\cdot\rangle$ on $\pg$
such that $(\pg,\Ad)$ is an orthogonal representation of $H$. Then $G/H$
becomes a manifold such that the canonical projection $\pi\colon
G\rightarrow G/H$ is differentiable (see theorem 3.58 in \cite{War}) and the differential
of $\pi$ at the identity $e\in G$ is a vector space isomorphism
$\pi_*\colon \pg\rightarrow T_H (G/H)$ and $\pi_*(\hg)=0$. 

We now recall the notion of an induced vector bundle. Let
$(V,\tau)$ be a representation of $H$. Let $\underline{V}=G\times_H V$ denote the
set of equivalence classes of $G\times V$ under the equivalence
relation
\[
(g,v)\sim (gh,\tau(h^{-1})v)\quad\textrm{for }g\in G,v\in V,h\in H.
\]
Then $\underline{V}$ is a smooth vector bundle over $G/H$ and we call it
the induced vector bundle of $(V,\tau)$ (see Kobayashi and
Nomizu~\cite{KN} p.54). $G$ acts on the
left of $\underline{V}$ by
\[
g'[g,v]=[g'g,v]\quad\textrm{for }g'\in G,[g,v]\in\underline{V}
\]
where $[g,v]$ denotes the equivalence class of $(g,v)\in G\times V$.

The tangent bundle $T(G/H)$ of $G/H$ can be identified with the induced vector
bundle $\underline{\pg}=G\times_H \pg$ of the representation
$(\Ad,\pg)$. Let $L_g\colon G/H\rightarrow G/H$ denote left
translation by $g\in G$, i.e.,  $L_g(g'H)=gg'H$ for $g'\in G$. For each $g\in G$, let
$(L_g)_*\colon T_H(G/H)\rightarrow T_{gH}(G/H)$ denote the differential of
$L_g$ at the identity $H\in G/H$. Then $L_*\colon h\mapsto (L_h)_*$ is a representation
of $H$ on $T_H(G/H)$ and under the
identification $\pg\cong T_H(G/H)$, it is
identified with $(\Ad,\pg)$, i.e.,  the following diagram commutes for
all $h\in H$.
\begin{equation}\label{AdH}
\xymatrix{
\pg\ar[r]^{\pi_*\hspace{7mm}}\ar[d]_{\Ad(h)}& T_H(G/H)\ar[d]^{(L_h)_*}\\
\pg\ar[r]^{\pi_*\hspace{7mm}}&T_H(G/H)}
\end{equation}
If $X\in
T_{gH}(G/H)$ then $X=(L_g)_*(\pi_*(\xi))$ for some $\xi\in\pg$
and we therefore get a map $G\times_H \pg\rightarrow T_{gH}(G/H)$ where
$[g,\xi]\mapsto (L_g)_*(\pi_*(\xi))$. This map is well-defined since
(\ref{AdH}) shows that for $h\in H$
\[
[gh,\Ad(h^{-1})\xi]\mapsto (L_{gh})_*(\pi_*(\Ad(h^{-1})\xi))=(L_{gh})_*(L_{h^{-1}})_*(\pi_*(\xi))=(L_g)_*(\pi_*(\xi)).
\]
Similarly, the complexified tangent bundle can be identified with
the induced bundle $\underline{\pg_{\CC}}$ (where the induced bundle
is obtained by extending $\Ad$ complex linearly to $\pg_{\CC}$).

We extend the inner product $\langle\cdot,\cdot\rangle$ on $\pg$ $G$-invariantly to $T(G/H)$, i.e., 
\[
\langle (L_g)_*X,(L_g)_*Y\rangle=\langle X,Y\rangle.
\]
for $X,Y\in T_H(G/H)\cong\pg$. This gives a $G$-invariant Riemannian
metric on $T(G/H)$. Similarly, the complexification of
$\langle\cdot,\cdot\rangle$ gives a hermitian inner product on
$\pg_{\CC}$ and this gives the complexified Riemannian
metric on the complexified tangent bundle $T(G/H)_{\CC}$.

\subsection{The spinor bundle on $G/H$}\label{sec:spinorb}
Now let $\SO{\pg}$ denote the space 
\[
\SO{\pg}=\{\varphi\in\End(\pg)\mid\langle\varphi\xi,\varphi\xi\rangle=\langle\xi,\xi\rangle\textrm{
  for all }\xi\in\pg,\det\varphi=1\}.
\]
Note
that since $H$ is connected and $\Ad$ continuous, we must have that
$\Ad(H)\subset\SO{\pg}$. Let
$\underline{\SO{\pg}}=G\times_H\SO{\pg}$ be the fibre bundle induced
by the representation in the group of diffeomorphisms given by 
\[
\xymatrix{
H\ar[r]^{\Ad\hspace{3mm}}&\SO{\pg}\ar[r]^{l\hspace{4mm}}&\Diff(\SO{\pg})
}
\]
where $l$ denotes left multiplication in $\SO{\pg}$. Then we can think
of $\underline{\SO{\pg}}$ as the principal $\SO{n}$ bundle of oriented
orthogonal frames of $T(G/H)$ where $n=\dim\pg$.

A spin structure of $T(G/H)$ is a principal $\SP{n}$ bundle
$P\rightarrow G/H$ such
that there exists a two fold covering bundle map
\[
\Psi\colon P\rightarrow \underline{\SO{\pg}}
\]
such that if $\psi\colon\SP{n}\rightarrow\SO{n}$ is the two fold
covering homomorphism of $\SO{n}$ then
\begin{equation}\label{Psi}
\Psi(px)=\Psi(p)\psi(x)\quad\textrm{for }p\in P,x\in \SP{n}.
\end{equation}

We now construct a specific spin structure of $T(G/H)$. Let $\SP{\pg}$
denote the two-fold cover $\SO{\pg}$ with covering map $\psi$. Suppose that the adjoint representation $\Ad\colon H\rightarrow\SO{\pg}$ lifts
to $\SP{\pg}$, i.e.,  there exists a homomorphism $\overline{\rho}\colon
H\rightarrow\SP{\pg}$ such that the following diagram commutes
\[
\xymatrix{
&\SP{\pg}\ar[d]^{\psi}\\
H\ar[r]^{\Ad}\ar[ur]^{\overline{\rho}}&\SO{\pg}.
}
\]
Note that if $\overline{\rho}$
exists, it is unique since if there exists another homomorphism
$\overline{\rho}'$ with this property, then for all $h\in H$
\[
\psi\left(\overline{\rho}(h)\overline{\rho}'(h^{-1})\right)=\psi(\overline{\rho}(h))\psi(\overline{\rho}'(h^{-1}))=\Ad(h)\Ad(h^{-1})=\operatorname{Id}
\]
and therefore
\[
\overline{\rho}(h)\overline{\rho}'(h^{-1})\in\{\pm 1\}\quad\textrm{for
  all }h\in H.
\]
Since $h\mapsto\overline{\rho}(h)\overline{\rho}'(h^{-1})$ is
continuous on $H$ and $H$ is connected, we must have that
$h\mapsto\overline{\rho}(h)\overline{\rho}'(h^{-1})$ is
constant. Since $\overline{\rho}(e)\overline{\rho}'(e^{-1})=1$, we
conclude that $\overline{\rho}'=\overline{\rho}$.

The fibre bundle $\underline{\SP{\pg}}=G\times_H\SP{\pg}$ induced by the representation
$l\circ\overline{\rho}$, where $l$ denotes left multiplication in
$\SP{\pg}$, is a principal $\SP{n}$ bundle. The covering homomorphism
$\psi\colon \SP{n}\rightarrow\SO{n}$ induces a bundle map $\Psi\colon \underline{\SP{\pg}}\rightarrow \underline{\SO{\pg}}$ given by
\[
\Psi[g,a]=[g,\psi(a)]\quad\textrm{for }[g,a]\in \underline{\SP{\pg}}.
\]
This is well-defined since
\[
\Psi[gh,\overline{\rho}(h^{-1})a]=[gh,\psi(\overline{\rho}(h^{-1})a)]=[gh,\Ad(h^{-1})\psi(a)]=[g,\psi(a)]\quad\textrm{for
  }h\in H
\]
and clearly $\Psi$ satisfies the condition (\ref{Psi}). Hence,
$\underline{\SP{\pg}}$ is a spin structure of $T(G/H)$. The spin
structure is $G$-invariant in the sense that $\Psi$ commutes with the
action of $G$.

Recall that we have the spinor representation
$\sigma\colon\SP{n}\rightarrow\GL(S)$. When the adjoint representation $\Ad\colon H\rightarrow\SO{\pg}$ lifts
to $\SP{\pg}$, we get the spinor representation
\begin{equation}\label{chi}
\chi=\sigma\circ\overline{\rho}
\end{equation}
of $H$ and when $n$ is even, we also have
the half spinor representations
\[
\chi^{\pm}=\sigma^{\pm}\circ\overline{\rho}.
\]

Now we define the complex spinor bundle $\underline{S}$ to
be the vector bundle induced on $G/H$ by $(S,\chi)$,
i.e.,  $\underline{S}=G \times_H S$. Note that this is equivalent to the
induced bundle $\underline{\SP{\pg}}\times_H S$ on the principal bundle
$\underline{\SP{\pg}}$ under the representation $(S,\sigma)$. Similarly, if $n$
is even, we get the half spinor bundles $\underline{S^{\pm}}$ induced by
$(S^{\pm},\chi^{\pm})$ and in this case we have that
\[
\underline{S}=\underline{S^+}\oplus\underline{S^-}.
\]

\section{The Dirac operator}\label{chap:dirac}
\subsection{Sections of induced bundles}
If $G/H$ has a spin structure, it is possible to define the Dirac
operator on $G/H$. This is an operator on the sections $\Gamma
(\underline{S\otimes V})$ of an induced vector bundle $\underline{S\otimes V}$
where $(S,\chi)$ is the spinor representation of $H$ and $(V,\tau)$ is
any complex representation of $H$. The Dirac operator is the
composition of reductive covariant differentiation which we describe
in section \ref{sec:connection} and Clifford multiplication which we
describe in section \ref{sec:dirac}. First we study the
sections of an induced bundle in more detail.

Let $(V,\tau)$ be a representation of $H$ over $\FF$ (where $\FF=\RR$
or $\CC$) and let $\underline{V}$ be the induced vector bundle on $G/H$. If $\varphi\in\Gamma
(\underline{V})$ is a section of $\underline{V}$, then  
it is of the form
\[
\varphi(gH)=[g,\hat{\varphi}(g)]\quad\textrm{for }g\in G
\]
where $\hat{\varphi}\colon G\rightarrow V$ is a smooth map satisfying
\[
\hat{\varphi}(gh)=\tau(h^{-1})\hat{\varphi}(g)\quad\textrm{for }g\in G,h\in H.
\] 
Hence we can make the identification
\[
\Gamma(\underline{V})\cong\{\hat{\varphi}\colon
G\rightarrow V\mid\hat{\varphi}\textrm{ smooth, }\hat{\varphi}(gh)=\tau(h^{-1})\hat{\varphi}(g)\textrm{ for }g\in
G,h\in H\}.
\]
Identifying the smooth maps $G\rightarrow V$ with $
C^{\infty}(G)\otimes V$ where $C^{\infty}(G)$ denotes the smooth maps
$G\rightarrow \FF$, we therefore get the identification
\begin{equation}\label{sections}
\Gamma(\underline{V})\cong\{\hat{\varphi}\in 
C^{\infty}(G)\otimes V\mid (\tilde{r}\otimes \tau)(h)\hat{\varphi}=\hat{\varphi}\textrm{ for }h\in H\}
\end{equation}
where $\tilde{r}$ denotes the right regular action of $G$ on
$C^{\infty}(G)$, i.e., 
\[
\tilde{r}(g)f(g')=f(g'g)\quad\textrm{for }g,g'\in
G,f\in C^{\infty}(G).
\] 
In
this identification, an element $\hat{\varphi}$
belonging to the right hand side of (\ref{sections}) is of the form
\[
\hat{\varphi}=\sum_i \tilde{\varphi}_i\otimes v_i
\]
where $\{v_i\}$ is a basis of $V$. By abuse of notation we write
$\tilde{\varphi}\otimes v$ to denote such an element and a section
$\varphi$ is identified with the element
$\hat{\varphi}=\tilde{\varphi}\otimes v$.

We now see that $\g$ acts on $\Gamma(\underline{V})$. 
$G$ acts on $\Gamma(\underline{V})$ by the left regular action $\lr$
given by
\[
\lr(g)\varphi(g'H)=g\cdot\varphi(g^{-1}g'H)\quad\textrm{for }g,g'\in G,\varphi\in\Gamma(\underline{V}).
\]
We have that
\[
\widehat{\lr(g)\varphi}(g')=\hat{\varphi}(g^{-1}g')
\]
so in terms of the isomorphism (\ref{sections}), the action is given by
\[
\lr(g)(\tilde{\varphi}\otimes
v)(g')=\tilde{\varphi}(g^{-1}g')\otimes v.
\]
The differential of $\lr$ gives an action $\lr_*$ of $\g$ on
$\Gamma(\underline{V})$ where
\[
\lr_{*}(\xi)(\tilde{\varphi}\otimes
v)(g')=\frac{d}{dt}\arrowvert_{t=0}\tilde{\varphi}(\exp(-t\xi)g')\otimes
v=(\tilde{\xi}\tilde{\varphi}\otimes v)(g')
\]
when $[g',\xi]=(L_{g'})_*\xi\in T_{g'H}(G/H)$ for $\xi\in\pg$ and
$\tilde{\xi}$ denotes the right invariant vector field on $G$ induced
by $\xi$. 

\subsection{Connections on induced bundles}\label{sec:connection}
In this section we define a way of differentiating sections of induced
vector bundles in terms of a vector field. This is called covariant
differentiation. In order to do so, we need the concept of a
connection. We start by considering connections and covariant
differentiation in general terms. We then study a specific
connection called the reductive connection.

Let $l_{g'}$, respectively $r_{g'}$ denote left and right translation
by $g'$ in $G$. A vector $X\in T_g(G)$ is called vertical if it is tangent to the
fibre $\pi^{-1}(gH)$, i.e.,  $\pi_*(l_{g^{-1}})_*X=0$. We denote the
space of vertical vectors in $T_g(G)$ by $V_g$. A connection on the
tangent bundle $T(G)$ of $G$ is
a choice of subspace $Q_g\subset T_g(G)$ for each $g\in G$ such that
\begin{description}
\item[(i)] $T_g(G)=V_g\oplus Q_g$
\item[(ii)] $Q_{gh}=(r_h)_*(Q_g)$
\item[(iii)]$Q_g$ depends differentiably on $g$.
\end{description}
Suppose we have a connection on $T(G)$. For each $g\in G$ the vectors of  $Q_g$ are
called horizontal. If $\gamma(t),t\in [0,t_0]$ is a (smooth) curve in $G/H$, a curve
$\tilde{\gamma}$ is called a horizontal lift of $\gamma$ if
$\pi(\tilde{\gamma}(t))=\gamma(t)$ for all $t$ and the
tangent vectors of $\tilde{\gamma}$ are all horizontal. Given a curve
$\gamma$ in $G/H$ and an element $g\in \pi^{-1}(\gamma(0))$, there is a
unique horizontal lift $\tilde{\gamma}$ of $\gamma$ such that
$\tilde{\gamma}(0)=g$. Hence for each curve $\gamma$ in $G/H$, we have
a well-defined map
$p(\gamma)_0^{t_0}\colon\pi^{-1}(\gamma(0))\rightarrow\pi^{-1}(\gamma(t_0))$
taking a point $g\in \pi^{-1}(\gamma(0))$ to $\tilde{\gamma}(t_0)$
where $\tilde{\gamma}$ is the horizontal lift of $\gamma$ with
starting point $g$. It can be shown that $p({\gamma})_0^{t_0}$ is a linear
isomorphism which is independent of the parametrization of
$\gamma$. It is called the parallel displacement along $\gamma$ from
$\gamma(0)$ to $\gamma(t_0)$. The inverse of $p(\gamma)_0^{t_0}$ is
parallel displacement along $\gamma^{-1}$. (See \cite{KN} section II.3).

Suppose that $\underline{V}$ is an induced real vector bundle on $G/H$ and
that $\gamma$ is a curve in $G/H$. Given an element
$[g,v]\in\underline{V}$, the horizontal lift of $\gamma$ to
$\underline{V}$
with starting point $[g,v]$ is $[\tilde{\gamma},v]$ where $\tilde{\gamma}$ is
the horizontal lift of $\gamma$ with starting point $g$. As before, we
get the parallel displacement of fibres along $\gamma$ from
$\gamma(0)$ to $\gamma(t_0)$ given by
\[
p(\gamma)_0^{t_0}[g,v]=[p(\gamma)_0^{t_0}g,v].
\]

We now define covariant differentiation. Let
$\varphi\in\Gamma(\underline{V})$ and $X\in T_{gH}(G/H)$. The
covariant derivative of $\varphi$ in the direction of $X$ is denoted
$\nabla_X\varphi$ and is defined as follows. Let $\gamma(t)$, $t\in
[-t_0,t_0]$ be a curve in $G/H$ with tangent vector $X$ at
$\gamma(0)=gH$. Then
\[
\nabla_X\varphi=\lim_{s\rightarrow 0}\frac{(p(\gamma)_0^s)^{-1}(\varphi(\gamma(s)))-\varphi(gH)}{s}.
\]
This is independent of the choice of $\gamma$. Now if $X$ is a vector
field on $G/H$ we define the covariant derivative of $\varphi$ with
respect to $X$ to be the element
$\nabla_X\varphi\in\Gamma(\underline{V})$ given by
\[
\nabla_X\varphi(gH)=\nabla_{X_{gH}}\varphi\quad\textrm{for }g\in G.
\]
This gives a map $\nabla\colon
\Gamma(\underline{V})\rightarrow\Gamma(T^*(G/H)\otimes\underline{V})$ which
we call the covariant derivative. We have the following result (see
\cite{KN} proposition III.1.2).
\begin{proposition}
The covariant derivative $\nabla\colon
\Gamma(\underline{V})\rightarrow\Gamma(T^*(G/H)\otimes\underline{V})$ is
linear and satisfies the Leibnitz rule, i.e.,  if 
$f\in C^{\infty}(G/H)$ and $\varphi\in\Gamma(\underline{V})$ then
\[
\nabla(f\varphi)=df\otimes\nabla\varphi+f\nabla\varphi.
\]
\end{proposition}

In our case the splitting $\g=\hg \oplus\pg$ gives a connection on $T(G)$
by defining $Q_g=(l_g)_*\pg$ for all $g\in G$ . We call this the reductive connection. Now we consider the covariant derivative given by this
connection. If $X=(L_g)_*\xi\in T_{gH}(G/H)$ where $\xi\in\pg$ (where we have identified $\pg$ with $T_H(G/H)$), then
\[
\gamma(t)=g\exp(t\xi)H,\quad\textrm{for }t\in [-t_0,t_0]
\]
defines a curve in $G/H$ with tangent $X$ at $\gamma(0)=gH$. In order to find the covariant derivative of a section
$\varphi\in \Gamma(\underline{V})$ in the direction of
$X$ we will need to compute
\[
(p(\gamma)_0^s)^{-1}(\varphi(\gamma(s)))
\]
for $s\in [0,t_0]$.
For $h\in
H$, the
unique horizontal lift $\tilde{\gamma}_h$ of $\gamma$ with
$\tilde{\gamma}_h(0)=gh$ is given by
\[
\tilde{\gamma}_h(t)=gh\exp(t\Ad(h^{-1})\xi)=g\exp(t\xi)h.
\]
It is clear that $\pi\circ\tilde{\gamma}_h(t)=\gamma(t)$ and
$\tilde{\gamma}$ is horizontal because of the $\Ad H$-invariance of $\pg$.
Writing $\varphi(gH)$ as $[g,\hat{\varphi}(g)]$ we have
that
\[
\varphi(\gamma(s))=[g\exp(s\xi),\hat{\varphi}(g\exp(s\xi))]=[\tilde{\gamma}_e(s),\hat{\varphi}(g\exp(s\xi))].
\]
Hence
\[
(p(\gamma)_0^s)^{-1}(\varphi(\gamma(s)))=[g,\hat{\varphi}(g\exp(s\xi))].
\]
Note that $(p(\gamma)_0^s)^{-1}(\varphi(\gamma(s)))$ is well-defined
since if we chose a different representative of $\varphi(\gamma(s))$, then
\[
\varphi(\gamma(s))=[g\exp(s\xi)h,\tau^{-1}(h)\hat{\varphi}(g\exp(s\xi)h)]=[\tilde{\gamma}_h(s),\tau^{-1}(h)\hat{\varphi}(g\exp(s\xi)h)]
\]
and
\begin{eqnarray*}
(p(\gamma)_0^s)^{-1}(\varphi(\gamma(s)))&=&(p(\gamma)_0^s)^{-1}[g\exp(s\xi)h,\tau^{-1}(h)\hat{\varphi}(g\exp(s\xi)h)]\\
&=&[gh,\tau^{-1}(h)\hat{\varphi}(g\exp(s\xi)h)]=[g,\hat{\varphi}(g\exp(s\xi))].
\end{eqnarray*}
We conclude that
\[
\nabla_X\varphi=\lim_{s\rightarrow 0}\frac{[g,\hat{\varphi}(g\exp(s\xi))]-[g,\hat{\varphi}(g)]}{s}=\left[g,\tfrac{d}{ds}\arrowvert_{s=0}\hat{\varphi}(g\exp(s\xi))\right].
\]
Now we describe $\nabla$ in terms of the isomorphism
(\ref{sections}). If $\tilde{\varphi}\otimes v\in
C^{\infty}(G)\otimes V$, i.e., 
$\tilde{\varphi}(g)\otimes v=\hat{\varphi}(g)$ for all $g\in G$ we have
that
\[
\tfrac{d}{dt}\arrowvert_{t=0}\hat{\varphi}(g\exp(t\xi))=(\rs(\xi)\otimes
1)(\tilde{\varphi}\otimes v)(g)
\]
where $\rs(\xi)$ denotes the differential of the right regular
action of $G$ on $C^{\infty}(G)$.
Hence
\[
\nabla_{[g,\xi]}(\tilde{\varphi}\otimes v)(g)=(\rs(\xi)\otimes
1)(\tilde{\varphi}\otimes v)(g).
\]

Now let $\{\xi_i\}$
be an orthonormal basis of $\pg$ and let $\{\xi_i^*\}$ be its dual
basis, i.e., 
\[
\xi_i^*(\xi)=\langle\xi_i,\xi\rangle\quad\textrm{for }\xi\in\pg.
\]
The map $\xi_i\mapsto\xi_i^*$ gives an isomorphism of $\pg$ to $\pg^*$.
Using the fact that $T(G/H)\cong\underline{\pg}$, a vector field $X$ on
$G/H$ can be thought of as an element $\tilde{X}\otimes \xi\in 
C^{\infty}(G)\otimes\pg$ such that $\tilde{X}(g)\otimes \xi=\hat{X}(g)$ where
$[g,\hat{X}(g)]=X(gH)$. Let $\hat{X}(g)=\xi\in\pg$.

We get that
\begin{eqnarray*}
\nabla_X(
\tilde{\varphi}\otimes v)(g)&=&\rs(\xi)\tilde{\varphi}(g)\otimes v=\sum_i
\rs(\xi_i^*(\xi)) \tilde{\varphi}(g)\otimes v\\
&=&\sum_i(\rs(\xi_i)\otimes\xi_i^*\otimes 1)(\tilde{\varphi}\otimes\xi\otimes v)(g).
\end{eqnarray*}
So we may think of
$\nabla\colon\Gamma(\underline{V})\rightarrow\Gamma(\underline{V}\otimes
T^*(G/H))\cong\Gamma(\underline{
  V\otimes\pg})$
as
\begin{equation}\label{nabla}
\nabla(\tilde{\varphi}\otimes v)=\sum_i \left((
  \rs(\xi_i)\otimes 1)(\tilde{\varphi}\otimes v)\right)\otimes \xi_i^*\cong\sum_i
 \left((\rs(\xi_i)\otimes 1)(\tilde{\varphi}\otimes v)\right)\otimes\xi_i
\end{equation}
when we think of sections in terms of (\ref{sections}).

Similarly, when $\underline{V}$ is an induced complex vector bundle
over $G/H$, we get the reductive covariant derivative
$\nabla\colon \Gamma(\underline{V})\rightarrow\Gamma(\underline{V}\otimes
(T^*(G/H))_{\CC})\cong\Gamma(\underline{
  V\otimes\pg_{\CC}})$ where $(T^*(G/H))_{\CC}$ denotes the
complexification of the cotangent bundle of $G/H$ and in terms of the
isomorphism (\ref{sections}), $\nabla$ is given by (\ref{nabla}).

\subsection{The Dirac operator on $G/H$}\label{sec:dirac}
Let $(V,\tau)$ be a complex representation of $H$. We then have the representation
$(S\otimes V,\chi\otimes\tau)$ of $H$ (see (\ref{chi})), and we can form the induced
bundle $\underline{S\otimes V}$ which has the reductive covariant
derivative $\nabla$.  In this section we define Clifford
multiplication and this
enables us to define the Dirac operator. In proposition
\ref{prop:dirac} we see that the
Dirac operator has a simple expression in terms of the isomorphism
(\ref{sections}). Finally, we look at some properties of the
Dirac operator.

Recall that the complex spinor representation
$\sigma\colon\SP{\pg}\rightarrow\GL(S)$ is the restriction of a
representation $\sigma\colon \CCL{\pg}\rightarrow\End(S)$. Since
$\pg_{\CC}\subset\CCL{\pg}$ we therefore
get a map $c\colon S\otimes \pg_{\CC}\rightarrow S$ defined by
\[
c(s\otimes \xi)=\sigma(\xi)s\quad\textrm{for }\xi\in\pg_{\CC},s\in S.
\]
We call this Clifford multiplication. Note that for $n=\dim\pg$ even, Clifford
multiplication maps $S^{\pm}\otimes \pg_{\CC}$ into $S^{\mp}$ respectively:
Recall that $S^{\pm}$ are the $\pm 1$ eigenspaces of
$\sigma(\omega_{\CC}')$ where $\omega_{\CC}'=i^{\tfrac{n}{2}}e_1\cdots
e_n$ for an oriented orthormal basis $\{e_j\}_{1\le j\le n}$ of $\pg$. Since
$n$ is even,
\[
\omega_{\CC}'e_j=i^{\tfrac{n}{2}}e_1\cdots
e_ne_j=(-1)^{n-1}i^{\tfrac{n}{2}}e_je_1\cdots
e_n=-e_j\omega_{\CC}'\quad\textrm{for }j\in\{1,\ldots,n\}
\]
and therefore for $s^{\pm}\in S^{\pm}$ and $\xi\in\pg_{\CC}$ 
\begin{eqnarray*}
\sigma(\omega_{\CC}')(c(
s^{\pm}\otimes\xi))&=&\sigma(\omega_{\CC}')\sigma(\xi)s^{\pm}=-\sigma(\xi)\sigma(\omega_{\CC}')s^{\pm}\\
&=&\mp\sigma(\xi)s^{\pm}=\mp
c(s^{\pm}\otimes\xi).
\end{eqnarray*}
Hence we get the maps $c^{\pm}\colon S^{\pm}\otimes\pg_{\CC}\rightarrow
S^{\mp}$ by restricting $c$.

Now construct the induced bundle
$\underline{S\otimes \pg_{\CC}}$ by the representation $(S\otimes\pg_{\CC},\chi\otimes
\Ad)$. Clifford multiplication then induces a bundle map
$c\colon\underline{S\otimes\pg_{\CC}}\rightarrow \underline{S}$ given by
\[
c[g,s\otimes \xi]=[g,\sigma(\xi)s]\quad\textrm{for }s\in S,\xi\in\pg.
\]
This is well-defined. To see this, observe that since
$\chi=\sigma\circ\overline{\rho}$, we have that for $h\in H$
\[
\sigma(\Ad(h^{-1})\xi)\chi(h^{-1})s=\sigma(\Ad(h^{-1})(\xi)\overline{\rho}(h^{-1}))s.
\]
Now since
\[
\Ad(h^{-1})\xi=\psi\circ\overline{\rho}(h^{-1})(\xi)=\overline{\rho}(h^{-1})\xi\overline{\rho}(h)
\]
we get that
\[
\sigma(\Ad(h^{-1})\xi)\chi(h^{-1})s=\sigma(\overline{\rho}(h^{-1})\xi)s=\chi(h^{-1})\sigma(\xi)s.
\]
Hence,
\[
c[gh,\Ad(h^{-1})\xi\otimes\chi(h^{-1})s]=[gh,\chi(h^{-1})\sigma(\xi)s]=c[g,\xi\otimes s].
\]
Similarly when $n$ is even, the maps $c^{\pm}\colon S^{\pm}\otimes \pg_{\CC}\rightarrow
S^{\mp}$ give bundle maps $c^{\pm}\colon\underline{S^{\pm}\otimes \pg_{\CC}}\rightarrow
\underline{S^{\mp}}$ on the bundles induced by the representations 
$(S^{\pm}\otimes\pg_{\CC},\chi^{\pm}\otimes\Ad)$ and
$(S^{\mp},\chi^{\mp})$.

Clifford multiplication $c\colon\underline{S\otimes
  \pg_{\CC}}\rightarrow\underline{S}$ induces a map
$c\colon\Gamma(\underline{S\otimes V\otimes
  \pg_{\CC}})\rightarrow\Gamma(\underline{S\otimes V})$, which we also call Clifford
multiplication, given by
\[
c(\tilde{\varphi}\otimes s\otimes v\otimes
\xi)(g)=\tilde{\varphi}(g)\otimes\sigma(\xi)s\otimes v
\]
where $\tilde{\varphi}(g)\otimes s\otimes v\otimes \xi=\hat{\varphi}(g)$ and $\varphi(gH)=[g,\hat{\varphi}(g)]$.

The Dirac operator
is 
the operator $D\colon
\Gamma(\underline{S\otimes V})\rightarrow\Gamma(\underline{S\otimes V})$ given by
\[
D=c\circ\nabla.
\]
Similarly, if $n=\dim\pg$ is even we have the reductive covariant derivatives
$\nabla^{\pm}$ on $\underline{S^{\pm}\otimes V}$ and Clifford
multiplications $c^{\pm}\colon\Gamma(\underline{S^{\pm}\otimes V\otimes
  \pg_{\CC}})\rightarrow\Gamma(\underline{S^{\mp}\otimes V})$ and we define the operators $D^{\pm}\colon
\Gamma(\underline{ S^{\pm}\otimes
  V})\rightarrow\Gamma(\underline{S^{\mp}\otimes V})$ by
\[
D^{\pm}=c^{\pm}\circ\nabla^{\pm}.
\]

Considering $\Gamma(\underline{S\otimes V})$ as in (\ref{sections}),
we see that $D$ has a simple expression.
\begin{proposition}\label{prop:dirac}
The Dirac operator $D\colon\Gamma(\underline{S\otimes V})\rightarrow
\Gamma(\underline{S\otimes V})$ is in terms of $(\ref{sections})$
given by
\[
D=\sum_i\rs(\xi_i)\otimes \sigma(\xi_i)\otimes 1
\]
where $\{\xi_i\}$ is an orthonormal basis of $\pg$ and $\rs$ denotes
the right regular action of $\g$ on $C^{\infty}(G)$.
\end{proposition}
\proof We have that
\begin{eqnarray*}
D(\tilde{\varphi}\otimes s\otimes v)&=& c\left(\sum_i
  (\rs(\xi_i)\otimes 1)
  (\tilde{\varphi}\otimes(s\otimes v))\otimes \xi_i\right)\\
&=&c\left(\sum_i
  \rs(\xi_i)\tilde{\varphi}\otimes s\otimes
  v\otimes\xi_i\right)\\
&=&\sum_i\rs(\xi_i)\tilde{\varphi}\otimes \sigma(\xi_i)s\otimes v\\
&=&\left(\sum_i\rs(\xi_i)\otimes \sigma(\xi_i)\otimes 1\right)(\tilde{\varphi}\otimes
s\otimes v).
\end{eqnarray*}
This completes the proof.
\qed

We see that $D$ is a first order differential operator which is
homogeneous (i.e.,  it commutes with the action $\tilde{l}$ of $G$ on
$\Gamma(\underline{S\otimes V})$). The symbol
$\sigma_D$ of $D$ is as follows. For $\xi\in\pg$, $\sigma_D(\xi)\colon
S\otimes V\rightarrow S\otimes V$ is the linear map given by
\[
\sigma_D(\xi)(s\otimes v)=(\sigma(\xi)\otimes 1)(s\otimes v).
\]
Since $D$ is homogeneous, this determines the symbol at all points
$g\in G$ (see \cite{Wal} lemma 5.5.1). 
$\sigma_D(\xi)$ is clearly an isomorphism for $\xi\ne 0$, since in this case we have that
\[
\sigma(\xi)\sigma(-q(\xi)^{-1}\xi)=\sigma(-q(\xi)^{-1}\xi^2)=\sigma(1)=\operatorname{Id}.
\]
Hence $D$ is elliptic. Therefore $\operatorname{ker}D$ is
finite-dimensional (see \cite{LM} theorem III.5.2). Since $D$ is homogeneous, $\operatorname{ker}D$
is invariant under the left regular action $\tilde{l}$ of $G$
on $\Gamma(\underline{S\otimes V})$ and so the restriction of $\tilde{l}$
to $\operatorname{ker}D$ gives a finite-dimensional representation of
$G$ which we denote by $\tilde{\pi}$. Similarly, if $n=\dim\pg$ is
even, the restriction of the left regular action of $G$ on 
$\Gamma(\underline{S^{\pm}\otimes V})$ to $\operatorname{ker}D^{\pm}$ give finite-dimensional
representations $\tilde{\pi}^{\pm}$ of $G$.

It is shown in \cite{Pa} lemma 4.1 that we can fix a hermitian inner product
$\langle\cdot,\cdot\rangle_S$ on $S$
with respect to which Clifford multiplication is skew-hermitian and $\chi$ is
unitary. Now by fixing a hermitian inner product $\langle\cdot,\cdot\rangle_V$ on $V$ such that $\tau$ is
unitary, we get an inner product on $S\otimes V$ given by
\[
\langle\sum_i s_i\otimes v_i,\sum_j s_j'\otimes v_j'\rangle=\sum_{i,j}\langle s_i,s_j'\rangle_S\langle v_i,v_j'\rangle_V
\]
when $s_i,s_j'\in S,v_i,v_j'\in V$ and
with respect to this inner product, 
$\chi\otimes\tau$ is unitary. We therefore have a hermitian inner product on each
fibre of  $\underline{S\otimes V}$ and we get an inner product on
$\Gamma(\underline{S\otimes V})$ given by
\[
\langle\varphi,\varphi'\rangle=\int_{G/H}\langle\varphi(x),\varphi'(x)\rangle
dx
\]
where we integrate with respect to a $G$-invariant measure on $G/H$
(such a measure exists according to corollary A4 of Dupont~\cite{Pju}). With respect to this inner product,
$\tilde{l}$ is unitary and lemma $4.2$ of \cite{Pa} shows that $D$ is
formally 
self-adjoint. If $\dim\pg$ is even, $D^+$ and $D^-$ are formal
adjoints of each other.

\section{Symmetric spaces}\label{chap:symm}
\subsection{Basic notions}
In the rest of this paper we consider the special case of a compact
symmetric space. In section \ref{sec:square} we see that in this
case, the square of the Dirac operator has a particularly simple
expression. This enables us to determine the representations in the
kernel of the Dirac operator in section \ref{sec:ker}. We start by establishing some basic
concepts concerning compact symmetric spaces. We then study the spinor
representation in this case.

Let $G$ be a compact connected semisimple Lie group with Lie algebra
$\g$. A Lie algebra automorphism
$\theta$ of $\g$ is called an involution of
$\g$ if $\theta^2=\operatorname{Id}$. Let $H$ be the analytic subgroup
of $G$ corresponding to the $+1$ eigenspace $\hg$ of $\g$. Then $G/H$
is called a symmetric space. We assume that
$\theta\ne\operatorname{Id}$, i.e., $\hg$ is a proper Lie
subalgebra of $\g$. Let $\pg$ denote the $-1$ eigenspace of
$\theta$ and let $B$ be the Killing form on $\g$. 
$-B$ is a positive definite bilinear
form on $\g$ and therefore it is an inner product
$\langle\cdot,\cdot\rangle$ on $\g$. (See Knapp~\cite{Kna} theorem 1.45 and
corollary 4.26). 
$\pg$ is the
orthogonal complement of $\hg$ with respect to
$\langle\cdot,\cdot\rangle$ and we have that
\[
\g=\hg\oplus\pg,
\]
\begin{equation}\label{symmetric}
[\hg,\hg]\subset\hg,\quad [\pg,\pg]\subset\hg,\quad [\hg,\pg]\subset\pg.
\end{equation}
For $h\in H$ and $\xi,\xi'\in\g$ we have that
\[
\langle\Ad(h)\xi,\Ad(h)\xi'\rangle=-B(\Ad(h)\xi,\Ad(h)\xi')=-B(\xi,\xi')=\langle\xi,\xi'\rangle
\]
(see \cite{Kna} proposition 1.119). Hence $\Ad(h)$ is orthogonal for all $h$. Since $\hg$
is $\Ad H$-invariant, we have that $\pg$ is $\Ad H$-invariant.

$H$ is not necessarily semisimple but since it is compact we have that
$\hg=[\hg,\hg]\oplus Z_{\hg}$ where $Z_{\hg}$ is the center of
$\hg$ (see corollary 4.25 of \cite{Kna}). We note that
\begin{equation}\label{orthogonal}
Z_{\hg}=[\hg,\hg]^{\perp}\cap\hg.
\end{equation}
To see this, let $\xi\in[\hg,\hg]^{\perp}\cap\hg$ and $\xi'\in\hg$. Then
\[
\langle[\xi,\xi'],[\xi,\xi']\rangle=\langle-\ad(\xi')(\xi),[\xi,\xi']\rangle=\langle\xi,\ad(\xi')[\xi,\xi']\rangle=\langle\xi,[\xi',[\xi,\xi']]\rangle=0,
\]
i.e.,  $[\xi,\xi']=0$ and therefore $\xi\in Z_{\hg}$. Hence
$[\hg,\hg]^{\perp}\cap\hg\subset Z_{\hg}$ and since $\dim
Z_{\hg}=\dim([\hg,\hg]^{\perp}\cap\hg)$, we conclude (\ref{orthogonal}).

In the rest of this paper we 
make the assumption that the rank of $G$ equals the rank of $H$. We
now see
that this implies that $\dim\pg$ is even. $\g_{\CC}$ has Cartan subalgebra $\tg_{\CC}$ where $\tg$ is the
Lie algebra of a maximal torus $T$ in $G$. Let $\tg'$ be
the Lie algebra of a maximal torus $T'$ of $H$ and suppose 
that $\dim\tg=\dim\tg'$. Then we may assume that $T=T'$.  Since $\g_{\CC}$ is semisimple, we get a root space
decomposition of $\g_{\CC}$
\[
\g_{\CC}=\tg_{\CC}\oplus\bigoplus_{\alpha\in\Delta^+}\g_{\pm\alpha}
\]
where $\Delta^+$ denotes the positive roots of $\g_{\CC}$ and for each
$\alpha\in\tg_{\CC}^*$,
\[
\g_{\alpha}=\{\xi\in\g_{\CC}\mid \ad(X)\xi=\alpha(X)\xi\textrm{ for all
  }X\in\tg_{\CC}\}.
\]
Since for each root $\alpha$, $\g_{\alpha}$ is the simultaneous eigenspace of
$\ad(\tg_{\CC})$ and $\hg$ and $\pg$ are $\ad(\tg_{\CC})$-invariant, we
have that either $\g_{\alpha}\subset\hg_{\CC}$ or
$\g_{\alpha}\subset\pg_{\CC}$. When $\g_{\alpha}\subset\hg_{\CC}$, it
is the root space $\hg_{\alpha}$ of $\alpha$ in $\hg_{\CC}$ and we get that
\[
\hg_{\CC}=\tg_{\CC}\oplus\bigoplus_{\alpha\in\Delta_{\hg}^+}\g_{\pm\alpha},\quad
\pg_{\CC}=\bigoplus_{\alpha\in\Delta_{\pg}^+}\g_{\pm\alpha}
\]
where the disjoint union $\Delta_{\hg}^+\cup\Delta_{\pg}^+=\Delta^+$.
Since each root space is one-dimensional, we get that
\begin{eqnarray*}
\dim(\g_{\CC}/\tg_{\CC})&=&\dim\bigoplus_{\alpha\in\Delta^+}\g_{\pm\alpha}\textrm{
  is even}\\
\dim(\hg_{\CC}/\tg_{\CC})&=&\dim\bigoplus_{\alpha\in\Delta_{\hg}^+}\g_{\pm\alpha}\textrm{
  is even}.
\end{eqnarray*}
Hence
\[
\dim\pg_{\CC}=\dim(\g_{\CC}/\hg_{\CC})=\dim(\g_{\CC}/\tg_{\CC})-\dim(\hg_{\CC}/\tg_{\CC})\textrm{
  is even}.
\]
The following are examples of symmetric spaces with the above
properties. We return to these examples in section \ref{chap:ex}.

\begin{example}\label{ex:3}
Let $G=\SO{3}$ and $H=\SO{2}\cong S^1$. $G$ and $H$ are compact and 
connected and $G$ is semisimple. (See \cite{Kna} section II.1 and \cite{Kna} proposition 1.136). The Lie algebras of $G$ and $H$ are given by
\begin{eqnarray*}
\g&=&\so{3}=\{X\in\M{3}(\RR)\mid X+X^t=0\}\\
\hg&=&\so{2}=\{X\in\M{2}(\RR)\mid X+X^t=0\}
\end{eqnarray*}
and their complexifications are
\begin{eqnarray*}
\g_{\CC}&=&\so{3,\CC}=\{X\in\M{3}(\CC)\mid X+X^t=0\}\\
\hg_{\CC}&=&\so{2,\CC}=\{X\in\M{2}(\CC)\mid X+X^t=0\}.
\end{eqnarray*}
Elements of $H$ are of the form
\[
h=\left(\begin{array}{cc}\cos\theta_1&\sin\theta_1\\
-\sin\theta_1&\cos\theta_1\end{array}\right)
\]
where $\theta_1\in\RR$ and we may think of $H$ as a subgroup of $G$ by
identifying an element $h\in H$ with the element
\[
h=\left(\begin{array}{cc}h&0\\0&1\end{array}\right)\in\SO{3}.
\]
Under this identification, we identify an element 
\[
X=\left(\begin{array}{cc}0&x\\-x&0\end{array}\right)\in\hg
\]
with the element
\[
X=\left(\begin{array}{ccc}0&x&0\\-x&0&0\\0&0&0\end{array}\right)\in\so{3}.
\]
Let $\theta\colon\g\rightarrow\g$ be defined by
\[
\theta\left(\begin{array}{ccc}0&y_{12}&y_{13}\\-y_{12}&0&y_{23}\\-y_{13}&-y_{23}&0\end{array}\right)=\left(\begin{array}{ccc}0&y_{12}&-y_{13}\\-y_{12}&0&-y_{23}\\y_{13}&y_{23}&0\end{array}\right).
\]
$\theta$ clearly is an involution of $\g$ with $+1$ eigenspace
$\hg$. Thus $G/H=S^2$ is a symmetric space. The maximal torus of $G$ and $H$ is
$H$, in particular, $G$ and $H$ have equal rank. Now define the
element $e_1\in \hg_{\CC}^*$ by
\[
e_1\left(\begin{array}{ccc}0&ix&0\\-ix&0&0\\0&0&0\end{array}\right)=x.
\]
Then $\Delta^+=\Delta_{\pg}^+=\{e_1\}$.
\end{example}

\begin{example}\label{ex:2m+1}
Let $n=2m\ge 4$ and $G=\SO{n+1},H=\SO{n}$. Note that $G$ and $H$ are
compact, connected and semisimple for all $n$. (See \cite{Kna} section II.1 and proposition 1.136).
The Lie algebras of $G$ and $H$ are given by
\begin{eqnarray*}
\g&=&\so{n+1}=\{X\in\M{n+1}(\RR)\mid X+X^t=0\}\\
\hg&=&\so{n}=\{X\in\M{n}(\RR)\mid X+X^t=0\}
\end{eqnarray*}
and their complexifications are
\begin{eqnarray*}
\g_{\CC}&=&\so{n+1,\CC}=\{X\in\M{n+1}(\CC)\mid X+X^t=0\}\\
\hg_{\CC}&=&\so{n,\CC}=\{X\in\M{n}(\CC)\mid X+X^t=0\}.
\end{eqnarray*}
Note that
\[
\dim\g=\frac{(n+1)^2-(n+1)}{2}=\frac{n(n+1)}{2},\quad\dim\hg=\frac{n^2-n}{2}=\frac{n(n-1)}{2}.
\]
We may think of $H$ as a subgroup
of $G$ when we identify the element $h\in \SO{n}$ with the matrix
\[
h=\left(\begin{array}{cc}h&0\\0&1\end{array}\right)\in \SO{n+1},
\]
and under this identification, we identify an element $X\in\hg$ with
\[
X=\left(\begin{array}{cc}X&0\\0&0\end{array}\right).
\]
We now consider the homogeneous space $G/H$ which can be 
 identified with the $n$-sphere $S^n$ (see \cite{War} 3.65 (a)). In the
 following we show that $G/H$ is a symmetric space. 
Any element of $\g$ is of the form
\[
Y=\left(\begin{array}{cccc}y_{1, 1}&\cdots&y_{1, n}&y_{1,n+1}\\
\vdots& &\vdots&\vdots\\
y_{n, 1}&\cdots&y_{n,n}&y_{n,n+1}\\
y_{n+1, 1}&\cdots&y_{n+1,n}&0
\end{array}\right)
\]
where $(y_{i,j})_{1\le i,j\le n}\in\hg$ and $y_{i,n+1}=-y_{n+1,i}$ for
$i=1,\ldots,n$. 
Now let $\theta\colon \g\rightarrow \g$ be given by
\[
\theta(Y)=\left(\begin{array}{cccc}y_{1, 1}&\cdots&y_{1, n}&-y_{1,n+1}\\
\vdots& &\vdots&\vdots\\
y_{n, 1}&\cdots&y_{n,n}&-y_{n,n+1}\\
-y_{n+1, 1}&\cdots&-y_{n+1,n}&0
\end{array}\right).
\]
 It is easy to see that $\theta$ is an involution of $\g$ where $\hg$ is the $+1$
 eigenspace. Hence $G/H$ is a symmetric space. We see that
 the elements of $\pg$ are of the form 
\[
Y=\left(\begin{array}{cccc}0&\cdots&0&y_{1,n+1}\\
\vdots& &\vdots&\vdots\\
0&\cdots&0&y_{n,n+1}\\
-y_{1, n+1}&\cdots&-y_{n,n+1}&0
\end{array}\right).
\]
The dimension of $\pg$ is
\[
\dim\pg=\dim\g-\dim\hg=n=2m.
\]
Because of the assumption that $n=2m$ is even, $G$ and $H$ are of equal
rank since in this case $\g_{\CC}$ and
$\hg_{\CC}$ have the same Cartan subalgebra $\tg_{\CC}$ given by the
elements of $\g_{\CC}$ of the form
\begin{equation}\label{cartan}
X=\left(\begin{array}{cccccc}0&ih_1&&&&\\
-ih_1&0&&&&\\
&&\ddots&&&\\
&&&0&ih_m&\\
&&&-ih_m&0&\\
&&&&&0\end{array}\right)
\end{equation}
where $h_1,\ldots h_m\in\CC$. Let $\tg$ consist of the elements of the
form (\ref{cartan}) for which $h_1,\ldots,h_m\in i\RR$. Then $\tg$ is
the Lie algebra of the maximal torus $T$ of $G$ and $H$. $T$ consists
of the elements of $\SO{n+1}$ of the form
\[
t=\left(\begin{array}{cccccc}\cos\theta_1&\sin\theta_1&&&&\\
-\sin\theta_1&\cos\theta_1&&&&\\
&&\ddots&&&\\
&&&\cos\theta_m&\sin\theta_m&\\
&&&-\sin\theta_m&\cos\theta_m\\
&&&&&1\end{array}\right)
\] 
where $\theta_1,\ldots,\theta_m\in\RR$. 
For each $X\in\tg_{\CC}$ of the form given in (\ref{cartan}), define 
\[
e_j(X)=h_j\quad\textrm{for }j\in\{1,\ldots,m\}.
\]
Then $e_j\in\tg_{\CC}^*$ and the positive roots of $G$ and $H$ are given by
\begin{eqnarray*}
\Delta^+&=&\{e_i\pm e_j\mid i<j\}\cup\{e_k\}\\
\Delta_{\hg}^+&=&\{e_i\pm e_j\mid i<j\}.
\end{eqnarray*}
\end{example}

\subsection{The spinor representation}
Assume that $\Ad\colon H\rightarrow \SO{\pg}$ lifts to
$\SP{\pg}$. We then 
have the spinor representation
$\chi=\sigma\circ\overline{\rho}\colon H\rightarrow \GL(S)$ of
$H$ (as in (\ref{chi})). The differential $\chi_*$ of $\chi$ at the identity is given by the
composition
\begin{eqnarray}\label{dchi}
\xymatrix{
\hg\ar[r]^{\ad\hspace{2mm}}&\so{\pg}\ar[r]^{(\psi_*)^{-1}}_{\cong}&\spin{\pg}\ar[r]^{\sigma}&\End(S)
}.
\end{eqnarray}
Using proposition \ref{prop:sigmaweights} we get the following result.
\begin{proposition}\label{prop:chiweights}
The weights of the spinor representation $\chi\colon H\rightarrow
\GL(S)$ for $\dim\pg=2m$ are given by
\[
\tilde{\lambda}_{\varepsilon}=\tfrac{1}{2}\sum_{k=1}^m\varepsilon_k
\alpha_k,\quad\textrm{for }\varepsilon=(\varepsilon_1,\ldots,\varepsilon_m)\in\{\pm 1\}^m
\]
where $\{\alpha_k\mid 1\le k\le m\}$ is an enumeration of the roots of
$\Delta_{\pg}^+$.

The weights of the representation $\chi^+$ are
$\{\tilde{\lambda}_{\varepsilon}\mid \varepsilon\in E^+\}$ and the weights of $\chi^-$ are
$\{\tilde{\lambda}_{\varepsilon}\mid \varepsilon\in E^-\}$ where
$E^{\pm}$ are as in (\ref{e+-}).

The multiplicity of each $\tilde{\lambda}_{\varepsilon}$ is the number
of ways in which $\tilde{\lambda}_{\varepsilon}$ can be written in the
above form.
\end{proposition}
\proof We think of $\Ad\colon H\rightarrow \SO{\pg}$ as the composition
\[
\xymatrix{H\ar[r]^{\Ad\hspace{2mm}}&\SO{\pg}\ar[r]^{\operatorname{Id}}&\SO{\pg}
}
\]
where $\operatorname{Id}$ denotes the standard representation of
$\SO{\pg}$ on $\pg$. Since $\Ad(T)$ is a compact abelian Lie subgroup of
$\SO{\pg}$, it must lie in a maximal torus $T'$ of $\SO{\pg}$. So we
get the following commutative diagram.
\[
\xymatrix{\tg\ar[d]_{\exp}\ar[r]^{\ad}&\tg'\ar[d]^{\exp}\ar[r]^{\operatorname{Id}_*\hspace{2mm}}&\so{\pg}\ar[d]^{\exp}\\
T\ar[r]^{\Ad}&T'\ar[r]^{\operatorname{Id}\hspace{2mm}}&\SO{\pg}
}
\]
where
$\operatorname{Id}_*=\ad_{\operatorname{Cl}}\circ\psi_*^{-1}$. Let $\bigoplus_{j=1}^m\pg_{\pm j}$ be the weight space
decomposition of $\pg_{\CC}$ with respect to $\operatorname{Id}$,
i.e.,  $\pg_{\pm j}$ are all one-dimensional and
\[
\operatorname{Id}_*(X)v=\pm\eta_j(X)v\quad\textrm{for
  }v\in\pg_{\pm j},X\in\tg'_{\CC}.
\]
We then have that
\[
\operatorname{Id}_*\circ\ad(\xi)v=\pm\eta_j(\ad(\xi))v\quad\textrm{for
  }v\in\pg_{\pm j},\xi\in\tg_{\CC}.
\]
If $\alpha_1,\ldots\alpha_m$ is an enumeration of the
roots of $\Delta_{\pg}^+$, then $\pg_{\CC}=\bigoplus_{j=1}^m\g_{\pm\alpha_j}$
where $\g_{\pm\alpha_j}$ are all one-dimensional and
\[
\operatorname{Id}_*\circ\ad(\xi)v=\pm\alpha_j(\xi)v\quad\textrm{for
  }v\in\g_{\pm \alpha_j},\xi\in\tg_{\CC}.
\]
Hence by re-enumeration, we may assume that
\[
(\eta_j\circ\ad)(\xi)=\alpha_j(\xi)\quad\textrm{for }\xi\in \tg_{\CC}, j\in\{1,\ldots,m\}.
\]
Since $\chi_*=\sigma\circ\psi_*^{-1}\circ\ad$, proposition
\ref{prop:sigmaweights} shows that
\begin{eqnarray*}
\chi_*(\xi)s&=&\sigma(\psi_*^{-1}\circ\ad(\xi))s=\lambda_{\varepsilon}(\psi_*^{-1}\circ\ad(\xi))s\\&=&\tfrac{1}{2}\sum_{j=1}^m\varepsilon_j\eta_j(\ad(\xi))s=\tfrac{1}{2}\sum_{j=1}^m\varepsilon_j\alpha_j(\xi)s
\end{eqnarray*}
for $s\in S_{\varepsilon},\xi\in\tg_{\CC}$. This proves the
proposition.  \qed

Now we compute an explicit expression for the spinor representation
$\chi_*$ which becomes useful in section \ref{sec:square} where we
find the square of the Dirac operator. Let
$\{X_1,\ldots,X_{2m}\}$ be an orthonormal basis of $\pg$. Recall that
$\{X_iX_j\in \CL{\pg}\mid i<j\}$ is a basis of $\spin{\pg}$. We have the
following result (corresponding to lemma 2.1 of \cite{Pa}).
\begin{lemma}\label{lemma:explicit}
For all $Y\in \hg$
\[
\chi_*(Y)=\tfrac{1}{4}\sum_{k,l=1}^{2m}\langle[Y,X_k],X_l\rangle \sigma(X_k)\sigma( X_l).
\]
\end{lemma}
\proof We have the following commutative diagram.
\begin{equation}\label{chistar}
\xymatrix{
&\so{\pg}\\
\hg\ar[ru]^{\ad}\ar[r]^{\overline{\rho}_*}\ar[rd]^{\chi_*}&\spin{\pg}\ar[u]_{\ad_{\operatorname{Cl}}}\ar[d]^{\sigma}\\
&\End(S)
}
\end{equation}
Let $Y\in\hg$. We have that 
\[
\overline{\rho}_*(Y)=\sum_{k<l}C_{kl}X_k X_l
\]
where $C_{kl}\in\RR$. For each $i\in\{1,\ldots 2m\}$ we get that
\begin{eqnarray*}
[Y,X_i]&=&\ad(Y)X_i=\ad_{\operatorname{Cl}}((\overline{\rho}_*)(Y))(X_i)\\
&=&\sum_{k<l}C_{kl}(X_k
X_l X_i-X_iX_kX_l)\\
&=&\sum_{l>i}2C_{il}X_l-\sum_{k<i}2C_{ki}X_k.
\end{eqnarray*}
Here we have used the fact that
\[
X_kX_l X_i-X_iX_kX_l=\left\{\begin{array}{ll}0&\textrm{if }i\ne k,i\ne
    j\\
-2X_k&\textrm{if }i=l\\
2X_l&\textrm{if }i=k\end{array}\right..
\]
Hence 
\[
\langle[Y,X_i],X_j\rangle=\left\{\begin{array}{ll}2C_{ij}&\textrm{for
      }i<j\\
-2C_{ij}&\textrm{for
      }i>j\\
0&\textrm{for }i=j
\end{array}\right. .
\]
This implies that
\[
\overline{\rho}_*(Y)=\tfrac{1}{2}\sum_{k<l}\langle[Y,X_k],X_l\rangle X_k X_l
\]
and that for $k>l$ and $i\in\{1,\ldots 2m\}$
\[
\langle[Y,X_k],X_l\rangle X_kX_l=\langle[Y,X_l],X_k\rangle X_lX_k,\quad\langle[Y,X_i],X_i\rangle=0.
\]
Hence
\[
\overline{\rho}_*(Y)=\tfrac{1}{4}\sum_{k,l=1}^{2m}\langle[Y,X_k],X_l\rangle X_k X_l
\] 
and the result follows from (\ref{chistar}). \qed

\subsection{The irreducible parts of $\chi$}
We now study the irreducible parts of $\chi$. In section
\ref{sec:casimir} we use this to find a simple expression for
$\chi_*(\Omega_H)$ where $\Omega_H$ is the Casimir element of $H$ and
this is then used to find $D^2$ in section \ref{sec:square}.

Since $\chi_*$ is a finite-dimensional
representation of a complex reductive Lie algebra, it is completely
reducible in sense that $S$ splits into a direct sum of invariant
subspaces and the restriction of $\chi_*$ to each of these is
irreducible (see \cite{Kna} theorem 5.29). Similarly, $\chi^{\pm}_*$ are
completely reducible. We now study this splitting in more detail. According to the theorem of the highest weight (theorem 5.110 in \cite{Kna}), each
irreducible component of $\chi$ is in one-to-one 
correspondance with the highest weight of that component. A weight of
an irreducible component of $\chi$ is of the form
$\tilde{\lambda_{\varepsilon}}$ as given in proposition
\ref{prop:chiweights}. Now let $W$ denote the Weyl group of $G$ and
let
\[
W_1=\{\sigma\in W\mid \Delta_{\hg}^+\subset \sigma\Delta^+\}.
\]
For each $\sigma\in W_1$ we have that
\[
\sigma\Delta^+=\Delta_{\hg}^+\cup\Delta_{\pg}^{\sigma}
\]
where
$\Delta_{\pg}^{\sigma}=\{\varepsilon_1^{\sigma}\alpha_1,\ldots,\varepsilon_m^{\sigma}\alpha_m\}$
for some $(\varepsilon_1^{\sigma},\ldots,\varepsilon_m^{\sigma})\in\{\pm 1\}^m$ and
where $\alpha_1,\ldots,\alpha_m$ is an enumeration of the elements of
$\Delta_{\pg}^+$. 
\begin{lemma}\label{lemma:bijection}
Let $W_H$ denote the Weyl group of $H$. The map
\[
\begin{array}{ccc} W_H\times W_1&\rightarrow&W\\
(s,\sigma)&\mapsto& s\sigma
\end{array}
\]
is a bijection.
\end{lemma}
\proof Let $\Delta_{\hg}=\Delta_{\hg}^+\cup(-\Delta_{\hg}^+)$ and
$\Delta_{\pg}=\Delta_{\pg}^+\cup(-\Delta_{\pg}^+)$. Suppose that $w\in
W$. We have that
$w\Delta^+=\Delta_{\hg}^{w}\cup\Delta_{\pg}^{w}$ where
$\Delta_{\hg}^{w}\subset\Delta_{\hg}$ and
$\Delta_{\pg}^{w}\subset\Delta_{\pg}$. Since $\Delta_{\hg}^{w}$ defines a
system of positive roots of $\hg_{\CC}$, there is a unique element
$s\in W_{H}$ such that $s\Delta_{\hg}^+=\Delta_{\hg}^{w}$ (see
theorem 1.8 of Humphreys~\cite{Hum}). Let
$\sigma=s^{-1}w$. Then
\[
\sigma\Delta^+=s^{-1}(\Delta_{\hg}^{w}\cup\Delta_{\pg}^{w})=\Delta_{\hg}^+\cup
s^{-1}\Delta_{\pg}^{w}.
\]
So $\sigma\in W_1$ and hence $w=s\sigma$ for $s\in W_H$ and $\sigma\in W_1$. This
shows surjectivity.

Now suppose that $s\sigma=s'\sigma'$ for $s,s'\in W_H,\sigma,\sigma'\in W_1$. We then have
that $\sigma'\sigma^{-1}=(s')^{-1}s\in W_H$. Hence
\begin{equation}\label{deltah}
\sigma'\sigma^{-1}\Delta_{\hg}=\Delta_{\hg}.
\end{equation} 
If $\alpha\in\Delta_{\hg}^+$,
then $\sigma^{-1}\alpha\in\Delta^+$ and therefore
$\sigma'\sigma^{-1}\alpha\in\Delta_{\hg}^+\cup\Delta_{\pg}^{\sigma'}$. But then
(\ref{deltah}) shows that $\sigma'\sigma^{-1}\alpha\in\Delta_{\hg}^+$ and therefore
\[
(s')^{-1}s\Delta_{\hg}^+=\sigma'\sigma^{-1}\Delta_{\hg}^+=\Delta_{\hg}^+
\]
Hence $(s')^{-1}s=1$ (by theorem 1.8 of \cite{Hum}) which implies that
$s'=s,\sigma'=\sigma$.   \qed

Now let
\[
\delta=\tfrac{1}{2}\sum_{\alpha\in\Delta^+}\alpha
\]
and let
\[
\delta_{\hg}=\tfrac{1}{2}\sum_{\alpha\in\Delta_{\hg}^+}\alpha,\quad\delta_{\pg}=\tfrac{1}{2}\sum_{\alpha\in\Delta_{\pg}^+}\alpha.
\]
For each $\sigma\in W_1$, let 
\[
\delta_{\pg}^{\sigma}=\tfrac{1}{2}\sum_{\alpha\in\Delta_{\pg}^{\sigma}}\alpha.
\]
Let $l$ denote
the length function on the Weyl group $W$, i.e.,  if $w\in W$, $l(w)$ is
the smallest number of reflections in simple roots the product of
which is $w$. Corollary 1.7 of \cite{Hum} shows that $l(w)$ is actually the
number of positive roots sent to negative roots by $w$. The sign function $\sign\colon
W\rightarrow \{\pm 1\}$ on $W$ is given by $\sign(w)=\det
w=(-1)^{l(w)}$. Hence if $\sigma\in W_1$, then since
$l(\sigma)$ is the number of negative elements of $\Delta_{\pg}^{\sigma}$, we have
that
\[
\delta_{\pg}^{\sigma}\textrm{ is a weight of
  }\left\{\begin{array}{ll}\chi^+&\textrm{if }\sign(\sigma)=1\\
\chi^-&\textrm{if }\sign(\sigma)=-1\end{array}\right..
\]
We have the following result.
\begin{lemma}\label{lemma:high}
For each $\sigma\in W_1$, the element $\delta_{\pg}^{\sigma}$ is the highest weight of an
irreducible component $\tau_{\sigma}$ of $\chi$. When $\sigma\ne \sigma'$, $\tau_{\sigma}$
and $\tau_{\sigma'}$ are not equivalent. 
\end{lemma}
\proof Any weight $\tilde{\lambda_{\varepsilon}}$ can be written as 
\[
\tilde{\lambda_{\varepsilon}}=\delta_{\pg}^{\sigma}-\sum_{\alpha\in\Phi_{\varepsilon}}\alpha
\] 
where
$\Phi_{\varepsilon}=\{\varepsilon_i^{\sigma}\alpha_i\in\Delta_{\pg}^{\sigma}\mid\varepsilon_i^{\sigma}\ne\varepsilon_i\}$.
Suppose that $\delta_{\pg}^{\sigma}$ is not a highest weight. Then there
would exist some $\varepsilon\in\{\pm 1\}^m$ and an element
$\beta=\sum_i n_i\beta_i\ne 0$  where $\{\beta_i\}\subset\Delta_{\hg}^+$ are simple roots of
$\hg_{\CC}$ and $n_i\in\{0,2,\ldots\}$ such that
\begin{equation}\label{nothigh}
\delta_{\pg}^{\sigma}+\beta=\delta_{\pg}^{\sigma}-\sum_{\alpha\in\Phi_{\varepsilon}}\alpha.
\end{equation}
Now we may redefine the notion of positivity on $\tg_{\CC}^*$ such
that $\sigma\Delta^+$ are the positive roots of $\g_{\CC}$. With this
notion of positivity we have that
\[
\beta'=\sum_{\alpha\in\Phi_{\varepsilon}}\alpha\ge 0.
\]
Since
$\Delta_{\hg}^+\subset \sigma\Delta^+$, we still have that $\beta>0$ and therefore
$\beta+\beta'>0$. However, (\ref{nothigh}) shows that
$\beta+\beta'=0$. This is a contradiction and therefore we conclude
that $\delta_{\pg}^{\sigma}$ is the highest weight of an irreducible component
of $\chi$. 

Now suppose that $\delta_{\pg}^{\sigma}=\delta_{\pg}^{\sigma'}$. Then
$\sigma\delta=\sigma'\delta$, i.e.,  
\[
\langle
\sigma^{-1}\sigma'\delta,\alpha\rangle=\langle\delta,\alpha\rangle>0\textrm{ for
  all }\alpha\in\Delta^+
\]
(see proposition 2.69 of \cite{Kna}). But then theorem 3.10.9 of Wallach~\cite{Wal} shows
that $\sigma^{-1}\sigma'=1$, i.e.,  $\sigma=\sigma'$. Hence $\tau_{\sigma}$ and $\tau_{\sigma'}$ are
inequivalent for $\sigma\ne \sigma'$.
\qed

We also have the following result.
\begin{lemma}\label{lemma:noncompact}
The representations $\chi^+$ and $\chi^-$ have no weights in common.
\end{lemma}
\proof  Let $\{\beta_1,\ldots,\beta_l\}$ be an enumeration of the
simple roots of $\Delta^+$ such that
$\{\beta_1,\ldots,\beta_k\}\subset\Delta_{\hg}^+$ and
$\{\beta_{k+1},\ldots,\beta_l\}\subset\Delta_{\pg}^+$. Note that
$\{\beta_{k+1},\ldots,\beta_l\}\ne\emptyset$ since we have assumed that
$\hg\ne\g$. Any positive 
element $\alpha$ of $\tg_{\CC}^*$ can be written uniquely as
\[
\alpha=\sum_{i=1}^l n_i\beta_i
\]
where $n_i\in\{0,1,\ldots\}$ and we call
$n(\alpha)=\sum_{i=1}^l n_i$ the level of $\alpha$. Let
\[
n_{\pg}(\alpha)=\sum_{i=k+1}^l n_i\beta_i.
\]
We claim that for $\alpha\in\Delta_{\hg}^+$, $n_{\pg}(\alpha)$ is
even and for $\alpha\in\Delta_{\pg}^+$, $n_{\pg}(\alpha)$ is
odd.
We prove this by induction on the level $n(\alpha)$ of $\alpha$. Suppose
$n(\alpha)=1$. If $\alpha\in\Delta_{\hg}^+$, then $\alpha=\beta_i$  for
some $i\in\{1,\ldots,k\}$ and therefore $n_{\pg}(\alpha)=0$. If $\alpha\in\Delta_{\pg}^+$, then $\alpha=\beta_i$  for
some $i\in\{k+1,\ldots,l\}$ and hence $n_{\pg}(\alpha)=1$. This proves the
case $n(\alpha)=1$. Now let $\alpha\in\Delta^+$ with $n(\alpha)>1$ and suppose 
that the result holds for all elements of $\Delta^+$ with level
less than $n(\alpha)$. 
Since $n(\alpha)=\sum_{i=1}^ln_i>1$ and each $n_i\in\{0,1,\ldots\}$,
there must be a $j\in\{1,\ldots,l\}$ such that
\[
\alpha-\beta_j=(n_j-1)\beta_j+\sum_{i\ne j}n_i\beta_i>0.
\] 
Proposition 2.48 (e) of \cite{Kna}
shows that $\alpha-\beta_j\in\Delta^+$. Hence by letting
$\beta=\alpha-\beta_j$ and $\beta'=\beta_j$ we have written $\alpha$
as a sum
$\alpha=\beta+\beta'$ where
$\beta,\beta'\in\Delta^+$ and
$n(\beta)<n(\alpha),n(\beta')<n(\alpha)$. Since
$n_{\pg}(\alpha)=n_{\pg}(\beta)+n_{\pg}(\beta')$ the claim is proved
if we can show the following:
\begin{eqnarray}\label{evenodd}
\beta,\beta'\in\Delta_{\hg}^+&\textrm{implies}&\alpha\in\Delta_{\hg}^+\nonumber\\
\beta,\beta'\in\Delta_{\pg}^+&\textrm{implies}&\alpha\in\Delta_{\hg}^+\\
\beta\in\Delta_{\hg}^+,\beta'\in\Delta_{\pg}^+&\textrm{implies}&\alpha\in\Delta_{\pg}^+.\nonumber
\end{eqnarray}
Corollary 2.35 of \cite{Kna} shows that
$[\g_{\beta},\g_{\beta'}]=\g_{\alpha}$ and hence using the relations
(\ref{symmetric}) we see that
\begin{eqnarray*}
\g_{\beta},\g_{\beta'}\subset\hg&\textrm{implies}&\g_{\alpha}\subset\hg\\
\g_{\beta},\g_{\beta'}\subset\pg&\textrm{implies}&\g_{\alpha}\subset\hg\\
\g_{\beta}\subset\hg,\g_{\beta'}\subset\pg&\textrm{implies}&\g_{\alpha}\subset\pg.
\end{eqnarray*}
This shows (\ref{evenodd}) and we have therefore proved the claim.

Suppose that $\chi^{\pm}$ have some weight in common. Then there are non-empty subsets
$\Phi,\Phi'\subset\Delta_{\pg}^+$ with $\vert\Phi\vert$ even and
$\vert\Phi'\vert$ odd such that
\[
\sum_{\alpha\in\Phi}\alpha=\sum_{\alpha\in\Phi'}\alpha.
\]
But then we have that
\begin{equation}\label{noncompact}
\sum_{\alpha\in\Phi}n_{\pg}(\alpha)=\sum_{\alpha\in\Phi'}n_{\pg}(\alpha).
\end{equation}
This is a contradiction since the above shows that the
left hand side of (\ref{noncompact}) is even and the right hand side is odd.
\qed

Using lemmas \ref{lemma:high} and \ref{lemma:noncompact} we now
find the splitting of $\chi$ into irreducible parts.   
We have the following result (as in \cite{Pa} lemma 2.2).
\begin{proposition}\label{prop:irred}
For each $\sigma\in W_1$, let $\tau_{\sigma}$ denote the irreducible
representation of $H$ with highest weight
$\delta_{\pg}^{\sigma}$. Then the multiplicity of each $\tau_{\sigma}$ in $\chi$ is
one and
\[
\chi^+=\bigoplus_{\sigma\in W_1^+}\tau_{\sigma},\quad \chi^-=\bigoplus_{\sigma\in W_1^-}\tau_{\sigma}
\]
where $W_1^+=\{\sigma\in W_1\mid \sign \sigma= 1\}$, $W_1^-=\{\sigma\in W_1\mid \sign
\sigma=- 1\}$. 
\end{proposition}
\proof We have that $\tau^+=\bigoplus_{\sigma\in W_1^+}\tau_{\sigma}$ is a
subrepresentation of $\chi^+$ and that $\tau^-=\bigoplus_{\sigma\in W_1^-}\tau_{\sigma}$ is a
subrepresentation of $\chi^-$. 
We show that
\begin{equation}\label{trace}
\trace\chi^+-\trace\tau^+=\trace\chi^+-\trace\tau^-=0.
\end{equation}
Recall that
$\delta_{\pg}=\tfrac{1}{2}\sum_{\alpha\in\Delta_{\pg}^+}\alpha$. We
now have that
\begin{eqnarray*}
\trace\chi^+-\trace\chi^-&=&\left(\sum_{\Phi\subset\Delta_{\pg}^+,\vert\Phi\vert
  \textrm{ even}}e^{\delta_{\pg}-\sum_{\alpha\in\Phi}\alpha}\right)-\left(\sum_{\Phi\subset\Delta_{\pg}^+,\vert\Phi\vert
  \textrm{ odd}}e^{\delta_{\pg}-\sum_{\alpha\in\Phi}\alpha}\right)\\
&=&e^{\delta_{\pg}}\sum_{\Phi\subset\Delta_{\pg}^+}(-1)^{\vert\Phi\vert}e^{\sum_{\alpha\in\Phi}\alpha}\\
&=&\left(\prod_{\alpha\in\Delta_{\pg}^+}e^{\frac{\alpha}{2}}\right)\left(\prod_{\alpha\in\Delta_{\pg}^+}(1-e^{-\alpha})\right)\\
&=&\prod_{\alpha\in\Delta_{\pg}^+}(e^{\frac{\alpha}{2}}-e^{-\frac{\alpha}{2}}).
\end{eqnarray*}
Recall that $\delta_{\hg}=\tfrac{1}{2}\sum_{\alpha\in\Delta_{\hg}^+}\alpha$ and $\delta=\tfrac{1}{2}\sum_{\alpha\in\Delta^+}\alpha$. Note that
for $\sigma\in W_1$ we have that
$\delta_{\pg}^{\sigma}+\delta_{\hg}=\sigma\delta$. Using Weyl's character formula
(theorem 5.75 of \cite{Kna})
and the bijection $W_H\times W_1\rightarrow W$ of lemma \ref{lemma:bijection} we get that
\begin{eqnarray*}
\trace\tau^+-\trace\tau^-&=&\left(\sum_{\sigma\in W^+_1}\sum_{s\in
  W_H}\sign(s)e^{s(\delta_{\pg}^{\sigma}+\delta_{\hg})}\right)\left(\sum_{s\in W_H}\sign(s)e^{s\delta_{\hg}}\right)^{-1}\\
& &-\left(\sum_{\sigma\in
  W^-_1}\sum_{s\in
  W_H}\sign(s)e^{s(\delta_{\pg}^{\sigma}+\delta_{\hg})}\right)\left(\sum_{s\in W_H}\sign(s)e^{s\delta_{\hg}}\right)^{-1}\\
&=&\left(\sum_{\sigma\in W_1}\sign(\sigma)\sum_{s\in
    W_H}\sign(s)e^{s\sigma\delta}\right)\left(\sum_{s\in
    W_H}\sign(s)e^{s\delta_{\hg}}\right)^{-1}\\
&=&\left(\sum_{w\in
    W}\sign(w)e^{w\delta}\right)\left(\sum_{s\in
    W_H}\sign(s)e^{s\delta_{\hg}}\right)^{-1}\\
&=&\prod_{\alpha\in\Delta^+}(e^{\frac{\alpha}{2}}-e^{-\frac{\alpha}{2}})\prod_{\alpha\in\Delta^+_{\hg}}(e^{\frac{\alpha}{2}}-e^{-\frac{\alpha}{2}})^{-1}\\
&=&\prod_{\alpha\in\Delta^+_{\pg}}(e^{\frac{\alpha}{2}}-e^{-\frac{\alpha}{2}})=\trace\chi^+-\trace\chi^-.
\end{eqnarray*}
We therefore conclude that
\[
\trace\chi^+-\trace\tau^+=\trace\chi^--\trace\tau^-.
\]
If this does not vanish, $\chi^+$ and $\chi^-$ must have some
irreducible component in common and therefore some weight of $\chi^+$
must also be a weight of $\chi^-$. This contradicts lemma
\ref{lemma:high} and we have
therefore proved (\ref{trace}) and hence that $\chi^+=\tau^+$ and
$\chi^-=\tau^-$. Since the $\tau_{\sigma}$'s are inequivalent, the
multiplicity of each of them is one. \qed

\section{The Dirac operator on symmetric spaces}\label{chap:dsymm}
\subsection{The action of Casimir elements}\label{sec:casimir}
In this section we study the action of the Casimir element of
$G$ on
the sections of an induced bundle $\underline{V}$ and the action of
the Casimir element of $H$ on the space of spinors $S$. Both of these
actions are important in determining the expression for the square of
the Dirac operator.

Let $\{Y_1,\ldots,Y_r,X_1,\ldots,X_{2m}\}$ be an orthonormal basis of $\g$
such that $\{Y_1,\ldots,Y_r\}$ is a basis of $\hg$ and
$\{X_1,\ldots,X_{2m}\}$ is a basis of $\pg$. We then have the
Casimir
element $\Omega_H\in U(\hg_{\CC})$ of $H$ and the Casimir element
$\Omega\in U(\g_{\CC})$ of $G$ given by
\[
\Omega_H=\sum_{i=1}^{r}-Y_i^2,\quad\Omega=\sum_{i=1}^{r}-Y_i^2+\sum_{i=1}^{2m}-X_i^2
\]
where $U(\hg_{\CC})$ and $U(\g_{\CC})$ denote the universal enveloping
algebras of $\hg_{\CC}$ and $\g_{\CC}$ respectively.

Recall that we have an action $\lr_*$ of $\g$ on $\Gamma(\underline{V})$ given by
\begin{equation}\label{lstarxi}
\lr_{*}(\xi)(\tilde{\varphi}\otimes
v)(g')=\frac{d}{dt}\arrowvert_{t=0}\tilde{\varphi}(\exp(-t\xi)g')\otimes
v=\tilde{\xi}\tilde{\varphi}\otimes v(g')
\end{equation}
where $\tilde{\xi}$ denotes the right invariant vector field on $G$ induced
by $\xi$. We may identify the tangent bundle of $G$ with the induced
vector bundle $G\times \g$ under the trivial representation $\{e\}$ of
$G$ and a vector $X=(l_g)_*\eta\in T_g(G)$ is identified with
$(g,\eta)$. In order to describe the vector field $\tilde{\xi}$ in
these terms, we would like to find $\eta(g)\in\g$ for each $g\in G$ such that
\[
\tilde{\xi}\tilde{\varphi}(g)=(l_g)_*\eta(g)\quad.
\]
We have that
\begin{eqnarray*}
\eta(g)\tilde{\varphi}(e)&=&(l_{g^{-1}})_*\tilde{\xi}\tilde{\varphi}(g)=\frac{d}{dt}\arrowvert_{t=0}\tilde{\varphi}(g^{-1}\exp(-t\xi)g)\\
&=&\frac{d}{dt}\arrowvert_{t=0}\tilde{\varphi}(\exp (-t\Ad(g^{-1})\xi))\\
&=&\Ad(g^{-1})\xi\tilde{\varphi}(e).
\end{eqnarray*}
Hence we identify $\tilde{\xi}(g)$ with $\Ad(g^{-1})\xi$.
Proposition 5.24 of \cite{Kna} shows that $\Omega$ is in the center of
$U(\g_{\CC})$. So if we denote
the element $\sum_{i=1}^{t}-\tilde{Y_i}^2+\sum_{i=1}^{2m}-\tilde{X_i}^2$ by
$\tilde{\Omega}$, then we may identify $\tilde{\Omega}(g)$ with
$\Ad(g^{-1})\Omega=\Omega$. Therefore
\begin{eqnarray*}
\tilde{\Omega}\tilde{\varphi}(g)&=&\Omega\tilde{\varphi}(e)\\
&=&\frac{d}{dt}\arrowvert_{t=0}\sum_{i=1}^{r}-\tilde{\varphi}(\exp(-2tY_i))+\sum_{i=1}^{2m}-\tilde{\varphi}(\exp(-2tX_i))\\
&=&\left(\sum_{i=1}^r-(\tilde{r}_*(Y_i))^2+\sum_{i=1}^{2m}-(\tilde{r}_*(X_i))^2\right)\tilde{\varphi}(e)\\
&=&\rs(\Omega)\tilde{\varphi}(g)
\end{eqnarray*}
where $\rs$ is the differential of the right regular
action of $\g$ on $C^{\infty}(G)$ which we extend to an action
of $U(\g_{\CC})$. Combining this with (\ref{lstarxi}) we get the
following result.
\begin{lemma}\label{lemma:lromega}
Let $\Omega$ be the Casimir element of $\g$. For each section $\tilde{\varphi}\otimes
v$ of an induced bundle $\underline{V}$ we have that
\[
\lr_*(\Omega)(\tilde{\varphi}\otimes
v)
=(\rs(\Omega)\otimes 1)(\tilde{\varphi}\otimes v)
\]
where $\lr_*$ denotes the left regular action of $U(\g_{\CC})$ on
$\Gamma(\underline{V})$ and $\rs$ denotes the
right regular action of $U(\g_{\CC})$ on $C^{\infty}(G)$.
\end{lemma}

Now we see how $\chi_*(\Omega_H)$ acts on $S$. We have the
following result.
\begin{lemma}\label{lemma:chiomega}
Let $\Omega_H$ be the Casimir element of $\hg$. $\chi_*(\Omega_H)$ acts on $S$ as scalar multiplication by $\langle
\delta,\delta\rangle-\langle\delta_{\hg},\delta_{\hg}\rangle$.
\end{lemma}
\proof By proposition \ref{prop:irred} we know that $\chi_*$ is a sum of
irreducible parts $(\tau_{\sigma})_*$ where $\sigma\in W_1$. Suppose that $H$ is
semisimple. By proposition
5.28 of \cite{Kna} we have that $(\tau_{\sigma})_*(\Omega_H)$ acts
as scalar multiplication by
$\langle\delta_{\pg}^{\sigma},\delta_{\pg}^{\sigma}+2\delta_{\hg}\rangle$. Since
\begin{eqnarray*}
\langle\delta_{\pg}^{\sigma},\delta_{\pg}^{\sigma}+2\delta_{\hg}\rangle&=&\langle\delta_{\pg}^{\sigma}+\delta_{\hg},\delta_{\pg}^{\sigma}+\delta_{\hg}\rangle-\langle\delta_{\hg},\delta_{\hg}\rangle=\langle
\sigma\delta,\sigma\delta\rangle-\langle\delta_{\hg},\delta_{\hg}\rangle\\
&=&\langle
\delta,\delta\rangle-\langle\delta_{\hg},\delta_{\hg}\rangle
\end{eqnarray*}
for all $\sigma\in W_1$, $\chi_*(\Omega_H)$ acts
as scalar multiplication by $\langle
\delta,\delta\rangle-\langle\delta_{\hg},\delta_{\hg}\rangle$ on
$S$. 
Suppose $H$ is not
semisimple, i.e.,  $\hg=[\hg,\hg]\oplus Z_{\hg}$ where $Z_{\hg}\ne
0$. Since (\ref{orthogonal}) shows that $Z_{\hg}=[\hg,\hg]^{\perp}\cap
\hg$, we may assume that
\[
[\hg,\hg]=\operatorname{span}\{Y_1,\ldots Y_t\},\quad Z_{\hg}=\operatorname{span}\{Y_{t+1},\ldots Y_r\}
\]
for some $t$. Let
\[
\Omega_{[\hg,\hg]}=\sum_{i=1}^t-Y_i^2,\quad\Omega_{Z_{\hg}}=\sum_{i=t+1}^r-Y_i^2.
\]
The irreducible
representations of $\hg$ are the irreducible representations of
$[\hg,\hg]$ extended to $\hg$ by being $0$ on $Z_{\hg}$ and the
irreducible representations of $Z_{\hg}$ extended to $\hg$ by being $0$ on
$[\hg,\hg]$. If $(\tau_{\sigma})_*$ is an irreducible representation
of $[\hg,\hg]$ with highest weight $\delta_{\pg}^{\sigma}$, then as in
the semisimple case we get that
\[
(\tau_{\sigma})_*(\Omega_H)=(\tau_{\sigma})_*(\Omega_{[\hg,\hg]})=\langle
\delta,\delta\rangle-\langle\delta_{\hg},\delta_{\hg}\rangle.
\]
The irreducible representations of $Z_{\hg}$ are just
linear functionals on $Z_{\hg}$ times $i$ and each representation has
itself as a weight. So if $(\tau_{\sigma})_*$ is an irreducible
representation of $Z_{\hg}$ with highest weight $\delta_{\pg}^{\sigma}$, then
\[
(\tau_{\sigma})_*(\Omega_H)=(\tau_{\sigma})_*(\Omega_{Z_{\hg}})=\vert\delta_{\pg}^{\sigma}\vert^2.
\]
Since $(Z_{\hg})_{\CC}\subset\tg_{\CC}$, we have that $(Z_{\hg})_{\CC}\perp
\g_{\alpha}$ for all $\alpha\in\Delta_{\hg}^+$. Hence $\langle\delta_{\pg}^{\sigma},\alpha\rangle=0$ for all
$\alpha\in\Delta_{\hg}^+$ and
therefore
$\vert\delta_{\pg}^{\sigma}\vert^2=\langle\delta_{\pg}^{\sigma},\delta_{\pg}^{\sigma}+2\delta_{\hg}\rangle$.
\qed


\subsection{The square of the Dirac operator}\label{sec:square}
We now show that in the symmetric case, the square of the Dirac
operator has a simple expression, namely, it consists of a costant
term plus the left regular action of the Casimir element of $G$. Note that the restriction of $D^2$ to
$\Gamma(\underline{S^+\otimes V})$ is $D^-\circ D^+$ and that the
restriction of $D^2$ to
$\Gamma(\underline{S^-\otimes V})$ is $D^+\circ D^-$. The expression
for $D^2$ shows that elements in the kernels of $D^{\pm}$ are
eigenvectors of the left regular action of the Casimir element of
$G$. We therefore obtain a criterion that must be met by the highest
weights of the irreducible subrepresentations which appear in the
kernels of $D^{\pm}$ (see corollary \ref{cor:criterion}). This is an
important step in determining the representations $\tilde{\pi}^{\pm}$
on the kernels of $D^{\pm}$ in section \ref{sec:ker}.

An element of $\tg_{\CC}^*$ is said to be analytically integral
if it
induces a character on the maximal torus $T$ of $H$ and $G$. Let
\[
\mathcal{F}=\left\{\mu\in\tg_{\CC}^*\mid
\mu\textrm{ is analytically integral }\right\}.
\] 
According to the theorem of the highest weight (theorem 5.110 of \cite{Kna}), the irreducible representations of $H$ stand in
one-to-one correspondance with the dominant analytically integral
forms of $H$, i.e.,  the elements
\[
\mathcal{F}_H=\left\{\mu\in\mathcal{F}\mid
\langle\mu,\alpha\rangle\ge
0 \textrm{ for all }\alpha\in\Delta_{\hg}^+\right\},
\]
the correspondance being that $\mu\in\mathcal{F}_H$ is the highest
weight of the representation. Similarly, the irreducible
representations of $G$ stand in one-to-one correspondance with 
\[
\mathcal{F}_G=\left\{\mu\in\mathcal{F}\mid
\langle\mu,\alpha\rangle\ge
0 \textrm{ for all }\alpha\in\Delta^+\right\}.
\]
Since $\delta_{\pg}$ is the highest
weight of an irreducible component of
$\chi^+$, $\delta_{\pg}$ is an element of $\mathcal{F}_H$ so if 
$\mu$ is any element of $\mathcal{F}_H$, then
$\lambda=\mu-\delta_{\pg}\in\mathcal{F}$. We now
have the following result (corresponding to proposition 3.2 of \cite{Pa}).
\begin{proposition}\label{prop:square}
Let $G$ and $H$ be compact connected
Lie groups of equal rank where $G$ is semisimple and $G/H$ is a
symmetric space and suppose that $\Ad\colon H\rightarrow \SO{\pg}$ lifts to
$\SP{\pg}$. For $\mu\in\mathcal{F}_H$, let $(V_{\mu},\tau_{\mu})$ be the irreducible complex representation of $H$ with highest
weight $\mu$. Let $\lambda=\mu-\delta_{\pg}$. Then the square of the Dirac operator on
$\Gamma(\underline{S\otimes V_{\mu}})$ is given by
\[
D_{\mu}^2=(\lr_{\mu})_*(\Omega)-\langle\lambda+2\delta,\lambda\rangle
1
\]
where $(\lr_{\mu})_*(\Omega)$ denotes the left regular
action of the Casimir element of $\g$ on $\Gamma(\underline{S\otimes V_{\mu}})$.
\end{proposition}
\proof For notational convenience we omit the subscript
$\mu$ in this proof. Proposition \ref{prop:dirac} showed that
in terms of the isomorphism 
(\ref{sections}), $D$ is  given by
\[
D=\sum_{i=1}^{2m}\rs(X_i)\otimes \sigma(X_i)\otimes 1.
\]
Since
\[
\sigma(X_i)^2=-1,\quad\sigma(X_i)\sigma(X_j)=-\sigma(X_j)\sigma(X_i)\textrm{
  for }i\ne j
\]
we get that
\begin{eqnarray}
D^2&=&\sum_{i=1}^{2m}\rs(X_i)^2\otimes\sigma(X_i)^2\otimes 1
+\sum_{i\ne j}\rs(X_i)\rs(X_j)\otimes\sigma(X_i)\sigma(X_j)\otimes 1\nonumber\\
&=&\sum_{i=1}^{2m}\rs(X_i)^2\otimes (-1)\otimes 1+\tfrac{1}{2}\sum_{i,j=1}^{2m}\rs[X_i,X_j]\otimes\sigma(X_i)\sigma(X_j)\otimes 1.\label{dsquare}
\end{eqnarray}
Since $[X_i,X_j]\in\hg$ for all $i,j$ and $\ad$ is skew-symmetric on $\g$ with respect to
$\langle\cdot,\cdot\rangle$, we get that
\begin{eqnarray*}
[X_i,X_j]&=&\sum_{q=1}^r\langle[X_i,X_j],Y_q\rangle
Y_q=\sum_{q=1}^r\langle\ad(X_i)X_j,Y_q\rangle Y_q\\
&=&\sum_{q=1}^r\langle-\ad(X_i)Y_q,X_j\rangle Y_q=\sum_{q=1}^r\langle[Y_q,X_i],X_j\rangle Y_q.
\end{eqnarray*}
Substituting this into (\ref{dsquare}) and using lemma
\ref{lemma:explicit} gives us that
\begin{eqnarray*}
D^2&=&\sum_{i=1}^{2m}\rs(X_i)^2\otimes
(-1)\otimes 1+\tfrac{1}{2}\sum_{i,j=1}^{2m}\sum_{q=1}^r\langle[Y_q,X_i],X_j\rangle\rs(Y_q)\otimes\sigma(X_i)\sigma(X_j)\otimes
1\\
&=&\sum_{i=1}^{2m}\rs(X_i)^2\otimes (-1)\otimes
1+2\sum_{q=1}^r\rs(Y_q)\otimes\chi_*(Y_q)\otimes 1.
\end{eqnarray*}
Now observe that
\begin{eqnarray*}
(\rs\otimes\chi_*)(Y_q)^2\otimes 1&=&\left(\rs(Y_q)\otimes 1\otimes
1+1\otimes\chi_*(Y_q)\otimes 1\right)^2\\
&=&\rs(Y_q)^2\otimes 1\otimes
1+1\otimes\chi_*(Y_q)^2\otimes 1+2\rs(Y_q)\otimes\chi_*(Y_q)\otimes 1.
\end{eqnarray*}
Differentiating the $H$-invariance condition on 
$C^{\infty}(G)\otimes S\otimes V$ of (\ref{sections}), we get that
\[
\rs(Y_q)\otimes 1\otimes 1+1\otimes\chi_*(Y_q)\otimes 1+1\otimes
1\otimes \tau_*(Y_q)=0
\]
and therefore
\[
(\rs\otimes\chi_*)(Y_q)\otimes 1=\rs(Y_q)\otimes 1\otimes 1+1\otimes\chi_*(Y_q)\otimes 1=-1\otimes
1\otimes \tau_*(Y_q).
\]
Hence
\begin{eqnarray*}
2r_*(Y_q)\otimes\chi_*(Y_q)\otimes 1&=&-\rs(Y_q)^2\otimes 1\otimes
1-1\otimes\chi_*(Y_q)^2\otimes 1+(\rs\otimes\chi_*)(Y_q)^2\otimes 1\\
&=&-\rs(Y_q)^2\otimes 1\otimes 1- 1\otimes\chi_*(Y_q)^2\otimes 1+1\otimes1\otimes\tau_*(Y_q)^2
\end{eqnarray*}
and therefore
\begin{eqnarray*}
D^2&=&\sum_{i=1}^{2m}-\rs(X_i)^2\otimes 1\otimes
1\\
& &+\sum_{q=1}^r\left(-\rs(Y_q)^2\otimes 1\otimes 1- 1\otimes\chi_*(Y_q)^2\otimes 1+1\otimes1\otimes\tau_*(Y_q)^2\right)\\
&=&\rs(\Omega)\otimes 1\otimes 1-\left(1\otimes
1\otimes \tau_*(\Omega_H)-1\otimes\chi_*(\Omega_H)\otimes 1\right).
\end{eqnarray*}
We know from lemma \ref{lemma:chiomega} that $\chi_*(\Omega_H)=(\langle
\delta,\delta\rangle-\langle\delta_{\hg},\delta_{\hg}\rangle) 1$ and
since $\tau$ is an irreducible representation of $H$ with highest
weight $\mu=\lambda+\delta_{\pg}$,
$\tau_*(\Omega_H)=\langle\lambda+\delta_{\pg}+2\delta_{\hg},\lambda+\delta_{\pg}\rangle 1$. So
\begin{eqnarray*}
1\otimes
1\otimes \tau_*(\Omega_H)-1\otimes\chi_*(\Omega_H)\otimes 1&=&\left(\langle\lambda+\delta_{\pg}+2\delta_{\hg},\lambda+\delta_{\pg}\rangle-\langle
\delta,\delta\rangle+\langle\delta_{\hg},\delta_{\hg}\rangle\right)
1\\
&=&\left(\langle\lambda+2\delta-\delta_{\pg},\lambda+\delta_{\pg}\rangle-\langle\delta,\delta\rangle+\langle\delta_{\hg},\delta_{\hg}\rangle\right)1.
\end{eqnarray*}
Now observe that
\begin{eqnarray*}
\langle\lambda+2\delta-\delta_{\pg},\lambda+\delta_{\pg}\rangle-\langle\delta,\delta\rangle+\langle\delta_{\hg},\delta_{\hg}\rangle
&=&\langle\lambda+2\delta,\lambda\rangle+\langle 2\delta-\delta_{\pg},\delta_{\pg}\rangle-\langle\delta,\delta\rangle+\langle\delta_{\hg},\delta_{\hg}\rangle\\
&=&\langle\lambda+2\delta,\lambda\rangle+2(\langle\delta_{\pg},\delta_{\pg}\rangle+\langle\delta_{\hg},\delta_{\pg}\rangle)-\langle\delta_{\pg},\delta_{\pg}\rangle\\
&&-(\langle\delta_{\pg},\delta_{\pg}\rangle+\langle\delta_{\hg},\delta_{\hg}\rangle+2\langle\delta_{\hg},\delta_{\pg}\rangle)+\langle\delta_{\hg},\delta_{\hg}\rangle\\
&=&\langle\lambda+2\delta,\lambda\rangle.
\end{eqnarray*}
We conclude that the square of the Dirac operator on
$\Gamma(\underline{S\otimes V_{\mu}})$ in terms of
(\ref{sections}) is given by
\[
D^2=\rs(\Omega)\otimes 1\otimes
1-\langle\lambda+2\delta,\lambda\rangle
1.
\]
Now using
lemma \ref{lemma:lromega} we get the the desired result. \qed

We now show that the assumption that $\Ad\colon H\rightarrow
\SO{\pg}$ lifts to $\SP{\pg}$ is not always necessary. Assume that $G$ and $H$ are compact
connected 
Lie groups of equal rank such that $G$ is semisimple and $G/H$ is a
symmetric space. Let $\psi_{H}\colon\tilde{H}\rightarrow H$ denote the
universal covering homomorphism of $H$. The kernel
$\operatorname{ker}\psi_H\cong \pi_1(H)$ is a discrete subgroup of the center
$Z(\tilde{H})$ of $\tilde{H}$. If $\psi_G\colon\tilde{G}\rightarrow G$ denotes the
universal covering homomorphism of $G$ and $\iota\colon H\rightarrow
G$ is the inclusion of $H$ in $G$, then $\iota$ lifts to
$\tilde{\iota}\colon\tilde{H}\rightarrow\tilde{G}$, i.e.,  the following
diagram commutes
\[
\xymatrix{
\tilde{H}\ar[d]^{\psi_H}\ar[r]^{\tilde{\iota}}&\tilde{G}\ar[d]^{\psi_G}\\
H\ar[r]^{\iota}&G}.
\]
Suppose $\dim\pg\ge 3$. Then $\psi\colon\SP{\pg}\rightarrow\SO{\pg}$
is the universal covering homomorphism of $\SO{\pg}$ and therefore
$\Ad\colon H\rightarrow\SO{\pg}$ lifts to a homomorphism
$\widetilde{\Ad}$ as follows
\[
\xymatrix{
\tilde{H}\ar[d]^{\psi_H}\ar[r]^{\widetilde{\Ad}}&\SP{\pg}\ar[d]^{\psi}\\
H\ar[r]^{\Ad}&\SO{\pg}}.
\]
Consider
$\widetilde{\Ad}\arrowvert_{\operatorname{ker}\psi_H}\colon\operatorname{ker}\psi_H\rightarrow\{\pm 1\}$
and let $K=\operatorname{ker}\widetilde{\Ad}\arrowvert_{\operatorname{ker}\psi_H}$. Then $K$ is
a subgroup of $\operatorname{ker}\psi_H$ of index either $1$ or
$2$ contained in $Z(\tilde{H})$. So if $H_1=\tilde{H}/K$, then $H_1$ is either a two-fold cover of
$H$ or it is $H$ itself. We get the following diagram
\[
\xymatrix{
\tilde{H}\ar@/_1pc/[dd]_{\psi_H}\ar[d]\ar[r]^{\tilde{\iota}}&\tilde{G}\ar[d]\ar@/^1pc/[dd]^{\psi_G}\\
H_1\ar[d]\ar[r]^{\iota_1}&G_1\ar[d]\\
H\ar[r]^{\iota}&G\\}
\]
where $G_1=\tilde{G}/\tilde{\iota}(K)$ and $\iota_1$ lifts
$\iota$. When $\iota_1$ is injective, then $H_1$ lies naturally
inside $G_1$ and $G_1/H_1=G/H$. Lifting $\Ad\colon H\rightarrow\SO{\pg}$ to the map $\Ad_1\colon
H_1\rightarrow \SP{\pg}$ we get the spinor representations 
$\chi=\sigma\circ\Ad_1\colon H_1\rightarrow\End(S)$ and
$\chi^{\pm}=\sigma^{\pm}\circ\Ad_1\colon H_1\rightarrow\End(S^{\pm})$ of $H_1$. 

Similarly, if $\dim\pg=2$ then $\RR$ is the universal cover of
$\SO{\pg}$ and if $\psi\circ\phi\colon\RR\rightarrow\SO{\pg}$
denotes the covering homomorphism, $\Ad\colon H\rightarrow\SO{\pg}$
lifts as follows
\[
\xymatrix{
\tilde{H}\ar[dd]^{\psi_H}\ar[r]^{\widetilde{\Ad}}&\RR\ar[d]^{\phi}\\
&\SP{\pg}\ar[d]^{\psi}\\
H\ar[r]^{\Ad}&\SO{\pg}}.
\]
We have that $\operatorname{ker}\phi=2\ZZ$. Consider
$\Ad\arrowvert_{\operatorname{ker}\psi_H}\colon
\operatorname{ker}\psi_H \rightarrow \ZZ$ and let
$K=\operatorname{ker}\Ad\arrowvert_{\operatorname{ker}\psi_H}$. Then
as before $H_1=\tilde{H}/K$ is either the trivial or two-fold cover of
$H$ and if $G_1=\tilde{G}/K$ and $\iota_1$, then $G_1/H_1=G/H$
and the lifting $\Ad_1\colon
H_1\rightarrow \SP{\pg}$ of $\Ad\colon H\rightarrow\SO{\pg}$ gives the
spinor representations 
$\chi=\sigma\circ\Ad_1\colon H_1\rightarrow\End(S)$ and $\chi^{\pm}=\sigma^{\pm}\circ\Ad_1\colon H_1\rightarrow\End(S^{\pm})$ of $H_1$.

We conclude that by going to the covering group $H_1$ of $H$ which is
either a two-fold cover or $H$ itself, we in some cases (namely when
$\iota_1$ is injective) get a spin structure on  
$G/H=G_1/H_1$, even if $\Ad\colon H\rightarrow
\SO{\pg}$ does not lift to $\SP{\pg}$. In the rest of this paper,
whenever $G_1$ and $H_1$ are mentioned, we assume that $\iota_1$ is injective.

Since $G$ and $H$ are compact, $G_1$ and $H_1$ are compact. Note that the Lie
algebras of $G_1$ and $H_1$ are $\g$ and $\hg$ respectively and that
the complexification of the Lie algebra of the maximal torus $T_1$ of
$H_1$ and $G_1$ is $\tg_{\CC}$. As before, we define the following
\begin{eqnarray*}
\mathcal{F}_1=\{\mu\in\tg_{\CC}^*\mid\mu\textrm{ induces a character
  on }T_1\}\\
\mathcal{F}_{H_1}=\{\mu\in\mathcal{F}_1\mid\langle\mu,\alpha\rangle\ge
0\textrm{ for all }\alpha\in\Delta_{\hg}^+\}\\
\mathcal{F}_{G_1}=\{\mu\in\mathcal{F}_1\mid\langle\mu,\alpha\rangle\ge
0\textrm{ for all }\alpha\in\Delta^+\}.
\end{eqnarray*}
Now let
\[
\mathcal{F}_{H_1}'=\{\mu\in\mathcal{F}_{H_1}\mid\mu-\delta_{\pg}\in\mathcal{F}\}.
\]
Note that if $\mu$ induces a character on $T$, then composing with the
covering homomorphism $H_1\rightarrow H$, we get a character on $T_1$
and hence $\mathcal{F}\subset\mathcal{F}_1$ and
$\mathcal{F}_H\subset\mathcal{F}_{H_1}$. If 
$\mu\in\mathcal{F}_H$, then $\mu-\delta_{\pg}\in\mathcal{F}$ and
therefore we have that
\begin{equation}\label{fh}
\mathcal{F}_H\subset\mathcal{F}_{H_1}'.
\end{equation}
Suppose
that $\mu\in\mathcal{F}_{H_1}'$ and let
$\lambda=\mu-\delta_{\pg}\in\mathcal{F}$. Let $(V_{\mu},\tau_{\mu})$ be
the irreducible representation of $H_1$ with highest weight $\mu$. This gives the representation
$(S\otimes V_{\mu},\chi\otimes\tau_{\mu})$ and we use this to form the induced bundle $\underline{S\otimes V_{\mu}}$ on $G_1/H_1$. As in remark 3.2 of \cite{Pa} we now
show that 
$\chi\otimes\tau_{\mu}$ actually is a
representation of $H$. Any weight $\eta$ of $\chi\otimes\tau_{\mu}$ is of the form $\eta=\gamma+\nu$ where $\gamma$ and $\nu$
are weights of $\chi$ and $\tau_{\mu}$ respectively.  Hence 
\[
\eta=\delta_{\pg}-\sum_{\alpha\in\Phi}\alpha+\mu-\sum_{\alpha\in\Delta_{\hg}^+}m_{\alpha}\alpha
\]
where $m_{\alpha}\in\{0,1,\ldots\}$ and
$\Phi\subset\Delta_{\pg}^+$. Therefore
\[
\eta=\lambda+2\delta_{\pg}-\sum_{\alpha\in\Delta^+}n_{\alpha}\alpha
\] 
where $n_{\alpha}\in\{0,1,\ldots\}$. Since each of the terms above
gives rise to a character on $T$, $\chi\otimes\tau_{\mu}$ is a
representation of $H$. So the following diagram commutes
\[
\xymatrix{G_1\times_{H_1}(S\otimes
V_{\mu})\ar[r]^{\cong}\ar[d]^{\pi'}&G\times_{H}(S\otimes
V_{\mu})\ar[d]^{\pi}\\
G_1/H_1\ar[r]^{\cong}&G/H
}
\]
where $\pi$ and $\pi'$ denote the usual projections. Hence we may identify the sections of $G_1/H_1$,
(i.e.,  the elements of $C^{\infty}(G_1)\otimes S\otimes
V_{\mu}$ that are invariant under the representation
$\tilde{r}\otimes\chi\otimes\tau_{\mu}$ of $H_1$)
with the sections of $G/H$, (i.e.,  the elements of $C^{\infty}(G)\otimes S\otimes
V_{\mu}$ that are invariant under the representation
$\tilde{r}\otimes\chi\otimes\tau_{\mu}$ of
$H$). Under this identification, the left regular action 
$\tilde{l}_{\mu}$ of $G_1$ on $\Gamma(\underline{S\otimes
V_{\mu}})$ corresponds to the left regular action of
$G$ on $\Gamma(\underline{S\otimes
V_{\mu}})$ which commutes with the Dirac operator $D_{\mu}\colon\Gamma(\underline{S\otimes
V_{\mu}})\rightarrow\Gamma(\underline{S\otimes
V_{\mu}})$. We also denote this action by
$\tilde{l}_{\mu}$. Note that the differential of both
actions $\tilde{l}_{\mu}$ gives the left
regular action of $\g$ on $\Gamma(\underline{S\otimes
V_{\mu}})$. In a similar way we use the representations
$\chi^{\pm}\otimes\tau_{\mu}$ to get the Dirac
operators $D_{\mu}^{\pm}\colon\Gamma(\underline{S^{\pm}\otimes
V_{\mu}})\rightarrow\Gamma(\underline{S^{\mp}\otimes
V_{\mu}})$ and the left regular action of $G_1$ on $\Gamma(\underline{S^{\pm}\otimes
V_{\mu}})$ give actions 
$\tilde{l}_{\mu}^{\pm}$ of $G$ that commute with
$D^{\pm}_{\mu}$.

We can now reformulate proposition \ref{prop:square} as follows.
\begin{proposition}\label{prop:square2}
Let $G$ and $H$ be compact connected 
Lie groups of equal rank where $G$ is semisimple and $G/H$ is a symmetric space.
For $\mu\in\mathcal{F}_{H_1}'$ let $(V_{\mu},\tau_{\mu})$ be the
irreducible representation of $H_1$ with highest weight
$\mu$. Let $\lambda=\mu-\delta_{\pg}$. Then the operator $D_{\mu}^2\colon
\Gamma(\underline{S\otimes V_{\mu}})\rightarrow
\Gamma(\underline{S\otimes V_{\mu}})$ is given by
\[
D_{\mu}^2=(\lr_{\mu})_*(\Omega)-\langle\lambda+2\delta,\lambda\rangle
1
\]
where $(\lr_{\mu})_*(\Omega)$ denotes the left regular
action of the Casimir element of $\g$ on $\Gamma(\underline{S\otimes V_{\mu}})$.
\end{proposition}
Note that in the case $G_1=G$, $H_1=H$, we have that
$\mathcal{F}_{H_1}'=\mathcal{F}_H$ and so in this case, proposition
\ref{prop:square2} is exactly proposition \ref{prop:square}.

The following result, which is a direct consequence of proposition
\ref{prop:square2}, is an important step in determining the
representations $\tilde{\pi}_{\mu}^{\pm}$ on the kernel of
$D_{\mu}^{\pm}$ in section \ref{sec:ker}.
\begin{corollary}\label{cor:criterion}
Let $\mu\in\mathcal{F}_{H_1}'$ and $\lambda=\mu-\delta_{\pg}$. If $\pi_{\nu}$ is an irreducible subrepresentation of $\tilde{\pi}_{\mu}^+$ or $\tilde{\pi}_{\mu}^-$ with highest
weight $\nu$, then
\[
\langle\nu+2\delta,\nu\rangle=\langle\lambda+2\delta,\lambda\rangle.
\]
\end{corollary}
\proof Suppose that $\pi_{\nu}$
is an irreducible subrepresentation of $\tilde{\pi}_{\mu}^+$ or $\tilde{\pi}_{\mu}^-$ with highest
weight $\nu$. Then since
$\operatorname{ker}D_{\mu}^{\pm}\subset\operatorname{ker}D_{\mu}^2$, we must have
that
\[
\pi_{\nu}(\Omega)-\langle\lambda+2\delta,\lambda\rangle=0.
\] 
On the other hand, proposition
5.28 of \cite{Kna} shows that
\[
\pi_{\nu}(\Omega)=\langle\nu+2\delta,\nu\rangle.
\]
We conclude that $\nu$ must satisfy the equation
\[
\langle\nu+2\delta,\nu\rangle=\langle\lambda+2\delta,\lambda\rangle.
\]  \qed
\subsection{The representations $\tilde{\pi}_{\mu}^{\pm}$}\label{sec:ker}
Recall that when we have assumptions as in proposition
\ref{prop:square2}, the left regular action $\tilde{l}_{\mu}^{\pm}$ of $G$ on
$\Gamma(\underline{S^{\pm}\otimes V_{\mu}})$ restricts to representations
$\tilde{\pi}_{\mu}^{\pm}$ on $\operatorname{ker}D_{\mu}^{\pm}$.  
In the rest of this section $G,H,G_1,H_1,\mu$ and $\lambda$ will be as in
proposition \ref{prop:square2}.  The main result of this
paper is as follows. 
\begin{theorem}\label{thm:repr}
Let $\mu\in\mathcal{F}_{H_1}'$, i.e., $\mu\in\mathcal{F}_{H_1}$ and
$\lambda=\mu-\delta_{\pg}\in\mathcal{F}$, and let $\dim\pg=2m$. 

For $m$ even we have that
\[
\begin{array}{ll}
\tilde{\pi}^+_{\mu}=\pi_{\sigma^{-1}(\lambda+\delta)-\delta}\textrm{
  and }\tilde{\pi}^-_{\mu}=0&\textrm{if
  }\sigma^{-1}(\lambda+\delta)-\delta\in\mathcal{F}_G, \sigma\in W_1^+\\
\tilde{\pi}^-_{\mu}=\pi_{\sigma^{-1}(\lambda+\delta)-\delta}\textrm{
  and }\tilde{\pi}^+_{\mu}=0&\textrm{if
  }\sigma^{-1}(\lambda+\delta)-\delta\in\mathcal{F}_G, \sigma\in
W_1^-  \\
\tilde{\pi}^+_{\mu}=\tilde{\pi}^-_{\mu}=0&\textrm{otherwise}.
\end{array}
\]

For $m$ odd we have
that
\[
\begin{array}{ll}
\tilde{\pi}^-_{\mu}=\pi_{\sigma^{-1}(\lambda+\delta)-\delta}\textrm{
  and }\tilde{\pi}^+_{\mu}=0&\textrm{if
  }\sigma^{-1}(\lambda+\delta)-\delta\in\mathcal{F}_G, \sigma\in W_1^+\\
\tilde{\pi}^+_{\mu}=\pi_{\sigma^{-1}(\lambda+\delta)-\delta}\textrm{
  and }\tilde{\pi}^-_{\mu}=0&\textrm{if
  }\sigma^{-1}(\lambda+\delta)-\delta\in\mathcal{F}_G, \sigma\in
W_1^-\\
\tilde{\pi}^+_{\mu}=\tilde{\pi}^-_{\mu}=0&\textrm{otherwise}.
\end{array}
\]

Here $\pi_{\sigma^{-1}(\lambda+\delta)+\delta}$ denotes the
irreducible representation of $G$ with highest weight
$\sigma^{-1}(\lambda+\delta)+\delta$ and $W_1^{\pm}$ are as in
proposition \ref{prop:irred}.
\end{theorem}

Now since for any $\nu\in\mathcal{F}_G$ we have that
$\nu\in\mathcal{F}_H\subset\mathcal{F}_{H_1}'$ (see (\ref{fh})) and
\[
(\nu+\delta)-\delta=\nu\in\mathcal{F}_G,
\]
theorem \ref{thm:repr} characterizes all irreducible representations
of $G$ as representations in the kernel of certain $D^+$ or $D^-$. We
prove theorem \ref{thm:repr} as follows. 
\begin{description}
\item[Step 1:] We study how $\tilde{l}^{\pm}_{\mu}$ breaks up into
  irreducible parts. 
\item[Step 2:] We prove proposition \ref{prop:subrep} which enables us
  to determine the multiplicities of
some of the irreducible parts of $\tilde{l}^{\pm}_{\mu}$. 
\item[Step 3:] Corollary \ref{cor:criterion}, which is a direct
consequence of the expression for $D^2_{\mu}$,
gives a criterion which must be satisfied by the subrepresentations of
$\tilde{\pi}_{\mu}^{\pm}$. This shows that the
irreducible parts of $\tilde{\pi}_{\mu}^{\pm}$ can
have multiplicity at most one and that
$\tilde{\pi}_{\mu}^+$ and
$\tilde{\pi}_{\mu}^-$ can have no irreducible parts
in common. 
\item[\textsc{Step 4:}] We complete the proof of theorem \ref{thm:repr} by
  combining step 3 with proposition \ref{prop:trace} which shows that
  except for possibly one irreducible part, the representations 
$\tilde{\pi}_{\mu}^+$ and
$\tilde{\pi}_{\mu}^-$ have all irreducible parts in
common.
\end{description} 

\textsc{Step 1. }
Consider the actions 
$\tilde{l}_{\mu}^{\pm}$. Extending
$\tilde{l}_{\mu}^{\pm}$ to the square integrable
sections $L_2(\underline{S\otimes
  V})$ of $\underline{S\otimes V}$, the Peter-Weyl theorem (section 3
of Bott~\cite{Bot}) shows that
\[
\tilde{l}_{\mu}^{\pm}=\sum_{\nu\in\mathcal{F}_G}\dim\Hom_G(V_{\nu},\Gamma(\underline{S^{\pm}\otimes
  V_{\mu}}))\pi_{\nu}
\]
where $(\pi_{\nu},V_{\nu})$ denotes the irreducible representation of
$G$ with highest weight $\nu$, $\Hom_G(V,V')$ denotes the continuous linear
maps respecting the action of $G$ on $V,V'$ and the sum is a unitary direct sum of
representations. Frobenius reciprocity (proposition 2.1 of \cite{Bot}) shows
that
\begin{eqnarray}\label{frob}
\Hom_G(V_{\nu},\Gamma(\underline{S^{\pm}\otimes
  V_{\mu}}))&=&\Hom_H(V_{\nu},S^{\pm}\otimes
  V_{\mu})\nonumber\\
&=&\Hom_{H_1}(V_{\nu}\otimes (S^{\pm})^*,
  V_{\mu})
\end{eqnarray}
where $(S^{\pm})^*$ denotes the duals of $S^{\pm}$. Note that $\Hom_G(V_{\nu},\Gamma(\underline{S^{\pm}\otimes
  V_{\mu}}))<\infty$. 
We combine this with the following lemma.
\begin{lemma}\label{lemma:dual}
Let $\dim\pg=2m$ and let $(S^{\pm})^*$ be the duals of $S^{\pm}$. As
representations of $H_1$ under the contragradient representations
$(\chi^{\pm})^c$ we have that
\[
(S^{\pm})^*=\left\{\begin{array}{ll}S^{\pm}&\textrm{if }m\textrm{ is even}\\
S^{\mp}&\textrm{if }m\textrm{ is odd}\end{array}\right..
\]
\end{lemma}
\proof Let $\langle\cdot,\cdot\rangle_S$ be the inner product on $S$
with respect to which $\chi$ is unitary. We identify $S^*$ with $S$
through the linear anti-isomorphism $*\colon S\rightarrow S^*$ given by
\[
s^*(s')=\langle s',s\rangle_S\quad\textrm{for all }s,s'\in S.
\]
Hence $s\mapsto s^*$ is an isomorphism
$\overline{S}\rightarrow S^*$ where $\overline{S}$ denotes the vector
space $S$ but where scalar multiplication is given by multiplication
with the conjugate scalars. The contragradient representation $\chi^c$ is given by
\begin{eqnarray*}
\chi^c(h)(s^*)(s')&=& s^*(\chi(h^{-1})s')=\langle
\chi(h^{-1})s',s\rangle_S\\
&=&\langle s',\chi(h)s\rangle_S=(\chi(h)s)^*(s')
\end{eqnarray*}
for $s,s'\in S,h\in H_1$.
 
Recall that $S$ is the unique representation on $\SP{n}$ which is the
restriction of an irreducible representation
of the Clifford algebra $\CCL{n}$. Since $\overline{S}$ as a
representation of $\CCL{n}$ is irreducible, we must have that
$\overline{S}=S$ as representations of $\SP{n}$. So it suffices to show that
\begin{equation}\label{s+-}
\overline{S^{\pm}}=\left\{\begin{array}{ll}S^{\pm}&\textrm{if }m\textrm{ is even}\\
S^{\mp}&\textrm{if }m\textrm{ is odd}\end{array}\right..
\end{equation}
 
Recall that $S^{\pm}$ are the $\pm 1$ eigenspaces of
$\sigma(\omega_{\CC}')$ where
\[
\omega_{\CC}'=i^{m}X_1\cdots X_{2m}.
\]
Now since for $s\in S$
\[
\overline{\sigma(\omega_{\CC}')s}=\left\{\begin{array}{ll}\sigma(\omega_{\CC}')s&\textrm{if }m\textrm{ is even}\\
-\sigma(\omega_{\CC}')s&\textrm{if }m\textrm{ is odd}\end{array}\right.,
\]
we see that (\ref{s+-}) holds.

\qed

In the following $\dim\pg=2m$. We conclude that
\begin{equation}\label{ltilde}
\tilde{l}^{\pm}_{\mu}=\left\{\begin{array}{ll}\sum_{\nu\in\mathcal{F}_G}\dim\Hom_{H_1}(S^{\pm}\otimes
V_{\nu},V_{\mu})\pi_{\nu}&\textrm{if }m\textrm{ is
  even}\\
\sum_{\nu\in\mathcal{F}_G}\dim\Hom_{H_1}(S^{\mp}\otimes
V_{\nu},V_{\mu})\pi_{\nu}&\textrm{if }m\textrm{ is odd}\end{array}\right..
\end{equation}

\textsc{Step 2.}
Proposition \ref{prop:subrep} below (as in \cite{Pa} lemma 8.1)
enables us to
determine the multiplicities of some of the irreducible parts
$\pi_{\nu}$ of
$\tilde{l}^{\pm}_{\mu}$. In the proof of proposition
\ref{prop:subrep} we use the following lemmas which can be found in 
Kostant~\cite{Kos}.
\begin{lemma}\label{lemma:kos1}
Let $\pi_{\delta}$ be the irreducible representation of $G_1$ with
highest weight $\delta$. The weights of $\pi_{\delta}$ are given by
\[
\rho_{\Phi}=\delta-\sum_{\alpha\in\Phi}\alpha
\]
where $\Phi\subset\Delta^+$. The multiplicity of each $\rho_{\Phi}$ is
the number of ways in which $\rho_{\Phi}$ can be expressed in this form.
\end{lemma}
\proof See lemma 5.9 of \cite{Kos}.
\begin{lemma}\label{lemma:kos2}
Let $\nu_1,\nu_2\in\mathcal{F}_{G_1}$ and let
$\pi_{\nu_1},\pi_{\nu_2}$ denote the irreducible representations of
$G_1$ with highest weights $\nu_1$ and $\nu_2$ respectively. Suppose
$\xi_1$ is a weight of $\pi_{\nu_1}$ and $\xi_2$ is a weight of
$\pi_{\nu_2}$. Then
\begin{equation}\label{inequality}
\vert\nu_1+\nu_2\vert\ge\vert\xi_1+\xi_2\vert
\end{equation}
and equality holds exactly when
\[
\nu_1=w\xi_1,\quad\nu_2=w\xi_2
\]
for some $w\in W$.
\end{lemma}
\proof (See lemma 5.8 of \cite{Kos}). If $\nu_1=w\xi_1,\nu_2=w\xi_2$
for some $w\in W$,
then clearly equality holds in (\ref{inequality}). We now show the
inequality (\ref{inequality}) and that equality only holds if
$\nu_1=w\xi_1,\nu_2=w\xi_2$ for some $w\in W$. Let $w\in W$ be such
that
\begin{equation}\label{dom}
\langle w(\xi_1+\xi_2),\alpha\rangle\ge 0\quad\textrm{for all }\alpha\in\Delta^+.
\end{equation}
Such a $w$ exists according to corollary 2.68 of \cite{Kna}. Now let
\[
\varphi_1=\nu_1-w\xi_1,\quad\varphi_2=\nu_2-w\xi_2.
\]
Since $w\xi_1$ is a weight of $\pi_{\nu_1}$ and $w\xi_2$ is a weight
of $\pi_{\nu_2}$ (see \cite{Kna} theorem 5.5(e)), we have that
\[
w\xi_1=\nu_1-\sum_{\alpha\in\Delta^+}n_{\alpha}\alpha,\quad w\xi_2=\nu_2-\sum_{\alpha\in\Delta^+}m_{\alpha}\alpha
\]
where $n_{\alpha},m_{\alpha}\in\{0,1,\ldots\}$. Hence
\begin{equation}\label{varphi12}
\varphi_1\ge0,\quad\varphi_2\ge0.
\end{equation}
Let
\[
\varphi=\varphi_1+\varphi_2=\nu_1+\nu_2-w(\xi_1+\xi_2).
\]
Then $\varphi\ge 0$ and we see that
\[
\vert\nu_1+\nu_2\vert^2=\vert\xi_1+\xi_2\vert^2+\vert\varphi\vert^2+2\langle
w(\xi_1+\xi_2)\varphi\rangle.
\]
Because of (\ref{dom}), $\langle w(\xi_1+\xi_2),\varphi\rangle\ge 0$,
and therefore we have proved the inequality
(\ref{inequality}). Furthermore, we see that if equality holds in
(\ref{inequality}), then $\vert\varphi\vert=0$, i.e., $\varphi=0$. But by
(\ref{varphi12}), this implies that $\varphi_1=\varphi_2=0$ , i.e.,
$\nu_1= w\xi_1$ and $\nu_2= w\xi_2$. This completes the proof. \qed
\begin{proposition}\label{prop:subrep}
Let $(\pi_{\nu},V_{\nu})$ be an irreducible representation of $G$ with
highest weight $\nu\in\mathcal{F}_G$. If $\xi\in\mathcal{F}_{H_1}$ is the
highest weight of an irreducible subrepresentation $\tau_{\xi}$ of
$\chi\otimes\pi_{\nu}$, then
\[
\vert\nu+\delta\vert\ge\vert\xi+\delta_{\hg}\vert
\]
and equality holds exactly when $\xi$ is of the form
\[
\xi_{\sigma}=\sigma(\nu+\delta)-\delta_{\hg}\quad\textrm{for some
  }\sigma\in W_1.
\]
The multiplicity of $\tau_{\xi_{\sigma}}$ in $\chi\otimes\pi_{\nu}$ is one
and $\tau_{\xi_{\sigma}}$ is a subrepresentation of
$\chi^+\otimes\pi_{\nu}$ when $\sigma\in W_1^+$ and a subrepresentation
of $\chi^-\otimes\pi_{\nu}$ when $\sigma\in W_1^-$ where $W_1^{\pm}$
are as in proposition \ref{prop:irred}.
\end{proposition}
\proof According to \cite{Kna} proposition 9.72, $\xi=\nu+\eta$ where $\eta$
is a weight of $\chi$. Hence
\[
\xi+\delta_{\hg}=\nu+\delta-\sum_{\alpha\in\Phi}\alpha
\]
where $\Phi\subset\Delta_{\pg}^+$. Lemma \ref{lemma:kos1} shows that $\delta-\sum_{\alpha\in\Phi}\alpha$ is a weight of
$\pi_{\delta}$ where $\pi_{\delta}$ denotes
the irreducible representation of $G_1$ with highest weight
$\delta$. Since $\nu$ is a weight of
$\pi_{\nu}$, lemma \ref{lemma:kos2} shows that $\vert\nu+\delta\vert\ge\vert\xi+\delta_{\hg}\vert$ and equality holds
if and only if $\xi+\delta_{\hg}=w(\nu+\delta)$ for some $w\in
W $.
Suppose that such a $w$ exists. Lemma \ref{lemma:bijection}
shows that there are unique $s\in W_H,\sigma\in W_1$ such that
$w=s\sigma$. We have that
\[
\langle s\sigma(\nu+\delta)-\delta_{\hg},\alpha\rangle=\langle\xi,\alpha\rangle\ge
0\quad\textrm{for all }\alpha\in\Delta_{\hg}^+.
\]
Proposition 2.69 of \cite{Kna} then shows that
\begin{equation}\label{ssigmanu}
\langle s\sigma(\nu+\delta),\alpha\rangle>0\quad\textrm{for all }\alpha\in\Delta_{\hg}^+.
\end{equation}
Similarly,
\[
\langle\nu+\delta,\alpha\rangle>0\quad\textrm{for all }\alpha\in\Delta^+
\]
and therefore
\[
\langle\sigma(\nu+\delta),\sigma\alpha\rangle>0\quad\textrm{for all }\alpha\in\Delta^+.
\]
Since $\Delta_{\hg}^+\subset\sigma\Delta^+$ we have that
\begin{equation}\label{sigmanu}
\langle\sigma(\nu+\delta),\alpha\rangle>0\quad\textrm{for all }\alpha\in\Delta^+_{\hg}.
\end{equation}
Theorem 3.10.9 of \cite{Wal}, (\ref{ssigmanu}) and (\ref{sigmanu}) now shows that $s=1$. We conclude that
$\xi=\xi_{\sigma}=\sigma(\nu+\delta)-\delta_{\hg}$.

Finally, we show that the multiplicity of $\tau_{\xi_{\sigma}}$ in
$\chi\otimes\pi_{\nu}$ is one and that $\tau_{\xi_{\sigma}}$ is a subrepresentation of
$\chi^{\pm}\otimes\pi_{\nu}$ when $\sigma\in W_1^{\pm}$. Using proposition
9.72 of \cite{Kna} again, we get that the weights of any irreducible
subrepresentation of $\chi^{\pm}\otimes\pi_{\nu}$ are of the form
$\lambda^{\pm}+\gamma$ where $\lambda^{\pm}$ is a weight of
$\chi^{\pm}$ and $\gamma$ is a weight of
$\pi_{\nu}$. In proposition \ref{prop:chiweights} we saw that
the weights of $\chi$ are of the form
$\tilde{\lambda}_{\varepsilon}=\delta_{\pg}-\sum_{\alpha\in\Phi_{\varepsilon}}\alpha$ where
$\Phi_{\varepsilon}\subset\Delta_{\pg}^+$ and each of the weight
spaces $S_{\varepsilon}$ are
one-dimensional. Hence we can find a basis of $S\otimes V_{\nu}$
consisting of elements $s_{\varepsilon}\otimes v_{\gamma}$ where
$s_{\varepsilon}$ is a weight vector of $\chi$ with weight $\tilde{\lambda}_{\varepsilon}$
and $v_{\gamma}$ is a weight vector of $\pi_{\nu}$ with weight
$\gamma$.  So if we can show that
$\xi_{\sigma}=\sigma\delta-\delta_{\hg}+\sigma\nu=\tilde{\lambda}_{\varepsilon}+\gamma$ implies that
$\tilde{\lambda}_{\varepsilon}=\sigma\delta-\delta_{\hg}$ and
$\gamma=\sigma\nu$, then we have proved the desired result. Suppose that
  $\xi_{\sigma}=\tilde{\lambda}_{\varepsilon}+\gamma$. Let
\[
\rho=\tilde{\lambda}_{\varepsilon}+\delta_{\hg}=\delta-\sum_{\alpha\in\Phi_{\varepsilon}}\alpha.
\]
This is a
  weight of $\pi_{\delta}$ according to lemma \ref{lemma:kos1} and we
  have that
\[
\sigma(\nu+\delta)=\rho+\gamma.
\]
Let 
\[
\psi_1=\delta-\sigma^{-1}\rho,\quad\psi_2=\nu-\sigma^{-1}\gamma.
\]
Since by theorem 5.5 (e) of \cite{Kna}, $\sigma^{-1}\rho$ is a weight of $\pi_{\delta}$ and
$\sigma^{-1}\gamma$ is a weight of $\pi_{\nu}$, we must have that
\[
\sigma^{-1}\rho=\delta-\sum_{\alpha\in\Delta^+}n_{\alpha}\alpha,\quad\sigma^{-1}\gamma=\nu-\sum_{\alpha\in\Delta^+}m_{\alpha}\alpha
\]
where $n_{\alpha},m_{\alpha}\in\{0,1,\ldots\}$ and therefore $\psi_1\ge 0,\psi_2\ge 0$. Since also
\[
\psi_1+\psi_2=\nu+\delta-\sigma^{-1}(\rho+\gamma)=0,
\]
we conclude that $\psi_1=\psi_2=0$. Hence $\rho=\sigma\delta$
and $\gamma=\sigma\nu$, i.e.,
$\tilde{\lambda}_{\varepsilon}=\sigma\delta-\delta_{\hg}$ and $\gamma=\sigma\nu$.  \qed

\textsc{Step 3. }
Now we turn to the subrepresentations $\tilde{\pi}_{\mu}^{\pm}$ of
$\tilde{l}^{\pm}_{\mu}$. Note that due to
(\ref{ltilde}) we have that
\[
\tilde{\pi}_{\mu}^{\pm}=\sum_{\nu\in\mathcal{F}_G}[\tilde{\pi}_{\mu}^{\pm}:\pi_{\nu}]\pi_{\nu}
\]
where the multiplicity $[\tilde{\pi}_{\mu}^{\pm}:\pi_{\nu}]$ of
$\pi_{\nu}$ in $\tilde{\pi}_{\mu}^{\pm}$ for each
$\nu\in\mathcal{F}_G$ satisfies 
\begin{equation}\label{multpinu}
[\tilde{\pi}_{\mu}^{\pm}:\pi_{\nu}]\le \left\{\begin{array}{ll}\dim\Hom_{H_1}(S^{\pm}\otimes
V_{\nu},V_{\mu})&\textrm{if }m\textrm{ is even}\\
\dim\Hom_{H_1}(S^{\mp}\otimes
V_{\nu},V_{\mu})&\textrm{if }m\textrm{ is odd}\end{array}\right..
\end{equation}

We are interested in the multiplicity of the representation
$\tau_{\xi}$ in $\chi\otimes\pi_{\nu}$ when
$\xi=\mu=\lambda+\delta_{\pg}$ and $\pi_{\nu}$ is a subrepresentation of
$\tilde{\pi}_{\mu}^+$ or
$\tilde{\pi}_{\mu}^-$.
Corollary \ref{cor:criterion} shows that in this case
$\vert\nu+\delta\vert=\vert\xi+\delta_{\hg}\vert$ since
\begin{eqnarray*}
\vert\nu+\delta\vert^2-\vert(\lambda+\delta_{\pg})+\delta_{\hg}\vert^2&=&\vert\nu+\delta\vert^2-\vert\lambda+\delta\vert^2\\
&=&\langle\nu,\nu\rangle+\langle\delta,\delta\rangle+2\langle\delta,\nu\rangle\\&
&-(\langle\lambda,\lambda\rangle+\langle\delta,\delta\rangle+2\langle\delta,\lambda\rangle)\\
&=&\langle\nu+2\delta,\nu\rangle-\langle\lambda+2\delta,\lambda\rangle=0.
\end{eqnarray*}
Hence if $\pi_{\nu}$ is a subrepresentation of
$\tilde{\pi}^+_{\mu}$,  proposition \ref{prop:subrep} and
(\ref{multpinu}) show that
\begin{equation}\label{multpi+}
[\tilde{\pi}^+_{\mu}:\pi_{\nu}]=1\textrm{ and }[\tilde{\pi}^-_{\mu}:\pi_{\nu}]=0
\end{equation}
and if $\pi_{\nu}$ is a subrepresentation of
$\tilde{\pi}^-_{\mu}$
\begin{equation}\label{multpi-}
[\tilde{\pi}^-_{\mu}:\pi_{\nu}]=1\textrm{ and }[\tilde{\pi}^+_{\mu}:\pi_{\nu}]=0.
\end{equation}

\textsc{Step 4. }
We now show (proposition \ref{prop:trace} below) that for each $\nu\in\mathcal{F}_G$
\[
[\tilde{\pi}^-_{\mu}:\pi_{\nu}]=[\tilde{\pi}^+_{\mu}:\pi_{\nu}]
\]
except possibly for one particular
$\nu_0\in\mathcal{F}_G$.
Combining this with (\ref{multpi+}) and (\ref{multpi-}) shows
that either
\begin{eqnarray*}
\tilde{\pi}^+_{\mu}=\tilde{\pi}^-_{\mu}=0&\textrm{or}&\\
\tilde{\pi}^+_{\mu}=\pi_{\nu_0}\textrm{ and }\tilde{\pi}^-_{\mu}=0&\textrm{or}&\\
\tilde{\pi}^-_{\mu}=\pi_{\nu_0}\textrm{ and }\tilde{\pi}^+_{\mu}=0.
\end{eqnarray*}
This completes the proof of theorem \ref{thm:repr}.
\begin{proposition}\label{prop:trace}
Let $\mu\in\mathcal{F}_{H_1}'$, $\lambda=\mu-\delta_{\pg}$ and $\dim\pg=2m$. Suppose that
$\sigma^{-1}(\lambda+\delta)-\delta\in\mathcal{F}_G$ for some
$\sigma\in W_1$. Then $\sigma$ is unique and 
\[
\trace \tilde{\pi}^+_{\mu}-\trace\tilde{\pi}^-_{\mu}=(-1)^m\sign(\sigma)\trace\pi_{\sigma^{-1}(\lambda+\delta)-\delta}.
\]
If no such $\sigma$ exists, then
\[
\trace \tilde{\pi}^+_{\mu}-\trace\tilde{\pi}^-_{\mu}=0.
\]
\end{proposition}
\proof First we prove the uniqueness of $\sigma$. Let
$\sigma^{-1}(\lambda+\delta)-\delta\in\mathcal{F}_G$. Then because of
proposition 2.69 of \cite{Kna},
\[
\langle\sigma^{-1}(\lambda+\delta),\alpha\rangle>0\quad\textrm{for }\alpha\in\Delta^+.
\]
Theorem 3.10.9 of \cite{Wal} now shows that $\sigma$ is unique.

For each
$\nu\in\mathcal{F}_G$ let
\[
\Gamma_{\nu}^{\pm}=\{\varphi(v)\mid v\in V_{\nu},\varphi\in\Hom_G(V_{\nu},\Gamma(\underline{S^{\pm}\otimes
  V_{\mu}}))\},
\]
where $(V_{\nu},\pi_{\nu})$ is the irreducible representation of $G$
with highest weight $\nu$, i.e.,  we have that
\[
L_2(\underline{S^{\pm}\otimes V_{\mu}})=\sum_{\nu\in\mathcal{F}_G}\Gamma_{\nu}^{\pm}
\]
where the sum is a unitary sum of representations under the left
regular representation. Note that
\[
\dim\Gamma_{\nu}^{\pm}=\dim V_{\nu}\cdot\dim\Hom_G(V_{\nu},\Gamma(\underline{S^{\pm}\otimes
  V_{\mu}}))<\infty.
\] 
Now let
\[
D_{\nu}^{\pm}=D_{\mu}^{\pm}\arrowvert_{\Gamma_{\nu}^{\pm}}.
\]
We note that
\[
\trace \tilde{\pi}^+_{\mu}-\trace\tilde{\pi}^-_{\mu}=\sum_{\nu\in\mathcal{F}_G}\ch(\operatorname{ker}D^+_{\nu})-\ch(\operatorname{ker}D^-_{\nu})
\]
where $\ch(V)$ denotes the trace of the representation of $G$ on $V$
under the left regular action. 
Since for all $v\in V_{\mu}$ and $\varphi\in\Hom_G(V_{\nu},\Gamma(\underline{S^{\pm}\otimes
  V_{\mu}}))$,
\[
D_{\nu}^{\pm}\varphi(v)=\varphi'(v)
\]
where $\varphi'=D^{\pm}_{\mu}\circ\varphi\in\Hom_G(V_{\nu},\Gamma(\underline{S^{\mp}\otimes
  V_{\mu}}))$, we have that
\[
D_{\nu}^{\pm}(\Gamma_{\nu}^{\pm})\subset\Gamma_{\nu}^{\mp}
\]
and $D_{\nu}^{\pm}$ are the formal adjoints of each other. So we have that
$D_{\nu}^+\colon\Gamma_{\nu}^+\rightarrow\Gamma_{\nu}^-$ gives an isomorphism
\[
(\operatorname{ker}D_{\nu}^+)^{\perp}\cong\operatorname{im}D_{\nu}^+=(\operatorname{ker}D_{\nu}^-)^{\perp}
\]
where $\perp$ is taken within $\Gamma_{\nu}^{\pm}$ with respect to the
inner product on $L_2(\underline{S^{\pm}\otimes
  V_{\mu}})$. Therefore
\[
\ch(\Gamma_{\nu}^{\pm})=\ch((\operatorname{ker}D_{\nu}^{\pm})^{\perp})+\ch(\operatorname{ker}D^{\pm}_{\nu})=\ch((\operatorname{ker}D_{\nu}^+)^{\perp})+\ch(\operatorname{ker}D^{\pm}_{\nu}).
\]
Hence using (\ref{frob}) we see that
\begin{eqnarray*}
\ch(\operatorname{ker}D^+_{\nu})-\ch(\operatorname{ker}D^-_{\nu})&=&\ch(\Gamma_{\nu}^+)-\ch(\Gamma_{\nu}^-)\\
&=&(-1)^m(\dim\Hom_{H_1}(S^+\otimes
V_{\nu},V_{\mu})\\
& &-\dim\Hom_{H_1}(S^-\otimes
V_{\nu},V_{\mu}))\trace\pi_{\nu}.
\end{eqnarray*}
In order to compute $\dim\Hom_{H_1}(S^+\otimes
V_{\nu},V_{\mu})-\dim\Hom_{H_1}(S^-\otimes
V_{\nu},V_{\mu})$ we calculate
\begin{equation}\label{tracechio}
\trace(\chi^+\otimes\pi_{\nu})-\trace(\chi^-\otimes\pi_{\nu})=(\trace\chi^+-\trace\chi^-)\trace\pi_{\nu}
\end{equation}
where we think of $\pi_{\nu}$ as a representation on $H_1$. Using the Weyl character formula (theorem 5.75 of \cite{Kna}) and lemma \ref{lemma:bijection} we get that
\begin{eqnarray*}
\trace{\pi_{\nu}}&=&\prod_{\alpha\in\Delta^+}(e^{\frac{\alpha}{2}}-e^{-\frac{\alpha}{2}})^{-1}\sum_{w\in
  W }\sign(w)e^{w(\nu+\delta)}\\
&=&(\trace\chi^+-\trace\chi^-)^{-1}\prod_{\alpha\in\Delta^+_{\hg}}(e^{\frac{\alpha}{2}}-e^{-\frac{\alpha}{2}})^{-1}\sum_{\sigma\in
  W_1}\sign(\sigma)\sum_{s\in
  W_H}\sign(s)e^{s\sigma(\nu+\delta)}.
\end{eqnarray*}
Here we have used the expression for $\trace\chi^+-\trace\chi^-$
found in the proof of proposition \ref{prop:irred}. Now since
$\nu\in\mathcal{F}_G$, we have that for each $\sigma\in W_1$
\[
\langle\nu+\delta,\alpha\rangle\ge\langle\delta,\alpha\rangle\quad\textrm{for
  all }\alpha\in\Delta^+
\]
and therefore
\[
\langle\sigma(\nu+\delta),\alpha\rangle\ge\langle\sigma\delta,\alpha\rangle\quad\textrm{for
  all }\alpha\in\sigma\Delta^+\supset\Delta_{\hg}^+.
\]
Hence
\[
\langle\sigma(\nu+\delta)-\delta_{\hg},\alpha\rangle\ge\langle\sigma\delta-\delta_{\hg},\alpha\rangle=\langle\delta_{\pg}^{\sigma},\alpha\rangle\ge
0\quad\textrm{for
  all }\alpha\in\Delta_{\hg}^+,
\]
i.e.,  $\sigma(\nu+\delta)-\delta_{\hg}\in\mathcal{F}_{H_1}$. The Weyl
character formula therefore shows that
\[
\trace\pi_{\nu}=(\trace\chi^+-\trace\chi^-)^{-1}\sum_{\sigma\in
  W_1}\sign(\sigma)\trace\tau_{\sigma(\nu+\delta)-\delta_{\hg}}
\]
where $\tau_{\sigma(\nu+\delta)-\delta_{\hg}}$ denotes the irreducible
representation of $H_1$ with highest weight
$\sigma(\nu+\delta)-\delta_{\hg}$. Inserting this into
(\ref{tracechio}) we get
\[
\trace(\chi^+\otimes\pi_{\nu})-\trace(\chi^-\otimes\pi_{\nu})=\sum_{\sigma\in
  W_1}\sign(\sigma)\trace\tau_{\sigma(\nu+\delta)-\delta_{\hg}}.
\]
Hence
\begin{eqnarray*}
\lefteqn{\dim\Hom_{H_1}(S^+\otimes
V_{\nu},V_{\mu})-\dim\Hom_{H_1}(S^-\otimes
V_{\nu},V_{\mu}){}}\\
& & \quad=\sum_{\sigma\in
  W_1}\sign(\sigma)\dim\Hom_{H_1}(V_{\sigma(\nu+\delta)-\delta_{\hg}},V_{\mu}).
\end{eqnarray*}
Therefore
\[
\trace \tilde{\pi}^+_{\mu}-\trace\tilde{\pi}^-_{\mu}=(-1)^m\sum_{\nu\in\mathcal{F}_G}\sum_{\sigma\in
  W_1}\sign(\sigma)\dim\Hom_{H_1}(V_{\sigma(\nu+\delta)-\delta_{\hg}},V_{\mu}).
\]
Note that
\[
\dim\Hom_{H_1}(V_{\sigma(\nu+\delta)-\delta_{\hg}},V_{\mu})=\left\{\begin{array}{ll}1&\textrm{if
      }\sigma(\nu+\delta)-\delta_{\hg}=\mu\\
0&\textrm{otherwise}\end{array}\right.
\]
and $\sigma(\nu+\delta)-\delta_{\hg}=\mu=\lambda+\delta_{\pg}$ exactly
when $\sigma^{-1}(\lambda+\delta)-\delta=\nu$. Since this can happen
for at most one $\sigma\in W_1$, we conclude that
\[
\trace
\tilde{\pi}^+_{\mu}-\trace\tilde{\pi}^-_{\mu}=\left\{\begin{array}{ll}(-1)^m\sign(\sigma)\trace\pi_{\sigma^{-1}(\lambda+\delta)-\delta}&\textrm{if
      }\sigma^{-1}(\lambda+\delta)-\delta\in\mathcal{F}_G\\
0&\textrm{otherwise}\end{array}\right..
\]
\qed

\begin{remark}
Note that since $\mathcal{F}_G$ is a lattice,
only finitely many $\nu\in\mathcal{F}_G$ can satisfy the equation
\[
\vert\nu+\delta\vert=\vert\lambda+\delta\vert.
\]
Hence $\tilde{\pi}_{\mu}^{\pm}$ can only have finitely many
subrepresentations $\pi_{\nu}$. This 
shows that in the symmetric case, the fact that the kernels
$\operatorname{ker}D_{\mu}^{\pm}$ of
$D_{\mu}^{\pm}$ are
finite-dimensional follows directly from corollary \ref{cor:criterion},
(i.e.,  without the ellipticity of $D_{\mu}$).
\end{remark}

\section{Examples}\label{chap:ex}
In the following we study some specific examples. The notation is as in section \ref{chap:dsymm}.
\subsection{$\SO{3} / \SO{2}$}
Let $G=\SO{3}$ and $H=\SO{2}$ as in example \ref{ex:3}. 
We now show that we can identify the representations
$\tilde{\pi}_{\mu}^+$ and
$\tilde{\pi}_{\mu}^-$ for all $\mu\in\mathcal{F}_{H_1}'$
as in theorem \ref{thm:repr}. Although theorem \ref{thm:repr}
gives us the result we need, we will in this case do some more direct
computations based on theorem 9.16 of \cite{Kna}. Using the
same notation as in example \ref{ex:3} we get that
\[
\delta=\delta_{\pg}=\tfrac{1}{2}e_1
\]
and
\[
W =W_1=\{1,s_{e_1}\}
\]
where $s_{e_1}$ denotes reflection in the root $e_1$. Since $\SO{\pg}=\SO{2}=H$,
$\Ad\colon H\rightarrow \SO{\pg}$ is the standard representation
and this lifts to the identity map
$H_1=\SP{2}\rightarrow\SP{\pg}=\SP{2}$. We get the $1$-dimensional spinor representations
$\chi^{\pm}$ of $H_1$ with highest weights
$\pm\tfrac{1}{2}e_1$. Since $H=S^1$, we must have that the
analytically integral forms of $H$ are exactly of the form $\lambda_1
e_1$ where $\lambda_1\in\ZZ$. Since $H$ has no roots, we therefore have that
\[
\mathcal{F}=\mathcal{F}_H=\{\lambda_1 e_1\mid \lambda_1\in\ZZ\}.
\]
Note that
\[
\mathcal{F}_G=\{\nu_1 e_1\mid \nu_1\in\ZZ,\quad\nu_1\ge 0\}.
\]
Since $H_1=\SP{2}$ is the two-fold cover of $H=S^1$, we have that
\[
\mathcal{F}_1=\mathcal{F}_{H_1}=\{\lambda_1 e_1\mid \lambda_1\in\ZZ
\textrm{ or }
\lambda_1\in \ZZ+\tfrac{1}{2}\}.
\]
Hence
\[
\mathcal{F}_{H_1}'=\{\mu_1e_1\mid\mu_1\in\ZZ+\tfrac{1}{2}\}.
\]
\begin{proposition}
For $G=\SO{3}$ and $H=\SO{2}$ let $\mu\in\mathcal{F}_{H_1}'$ and $\lambda=\mu-\delta_{\pg}=\lambda_1e_1$. We have that
\[
\begin{array}{ll}
\tilde{\pi}_{\mu}^-=\pi_{\lambda}\textrm{ and }\tilde{\pi}_{\mu}^+=0&\textrm{for
  }\lambda_1\ge 0\\
\tilde{\pi}_{\mu}^+=\pi_{-(\lambda+e_1)}\textrm{ and }\tilde{\pi}_{\mu}^-=0&\textrm{for
  }\lambda_1\le -1\end{array}.
\] 
\end{proposition}
\proof Let 
$(V_{\mu},\tau_{\mu})$ be an irreducible representation of $H_1$
with highest weight $\mu=\lambda+\delta_{\pg}=(\lambda_1+\tfrac{1}{2})e_1$. Then
$\dim V_{\mu}=1$ and therefore $\chi^{\pm}\otimes \tau_{\mu}$
are 1-dimensional and hence irreducible representations of $H$ with
highest weights $(\pm\tfrac{1}{2}+\lambda_1+\tfrac{1}{2})e_1$. In the following
we show that any irreducible component of
$\tilde{\pi}_{\mu}^{\pm}$ has multiplicity at most $1$. Let $\pi_{\nu}$ be an irreducible representation of $G$
with highest weight $\nu=\nu_1 e_1\in\mathcal{F}_G$. Frobenius reciprocity (theorem 9.9 of \cite{Kna}) gives us that
\[
[\tilde{\pi}_{\mu}^{\pm}:\pi_{\nu}]\le [\tilde{l}_{\mu}^{\pm}:\pi_{\nu}]=[\pi_{\nu}\arrowvert_H:\chi^{\pm}\otimes\tau_{\mu}].
\] 
Theorem 9.16 of \cite{Kna} shows that a restriction of an irreducible
representation of $G$ to $H$ decomposes with multiplicity one into
irreducible representations of $H$. Hence
\[
[\tilde{\pi}_{\mu}^{\pm}:\pi_{\nu}]\le 1.
\]
Now we show that $\tilde{\pi}_{\mu}^{\pm}$ has at most one irreducible
component. If $\pi_{\nu}$ is an irreducible subrepresentation of
$\tilde{\pi}_{\mu}^+$ or $\tilde{\pi}_{\mu}^-$, then
corollary
\ref{cor:criterion} shows that
\[
(\nu_1^2+\nu_1)-(\lambda_1^2+\lambda_1)=(\nu_1-\lambda_1)(\nu_1+\lambda_1+1)=0,
\]
i.e., 
\begin{equation}\label{nu1}
\nu_1=\lambda_1\quad\textrm{or}\quad\nu_1=-(\lambda_1+1).
\end{equation}
Using Frobenius
reciprocity again we see that in order for $[\tilde{\pi}_{\mu}^{\pm}:\pi_{\nu}]$ to be
non-zero, $\chi^{\pm}\otimes\tau_{\mu}$ must be an irreducible
subrepresentation of $\pi_{\nu}\arrowvert_H$. Recall that
$\chi^+\otimes\tau_{\mu}$ has highest weight
$(\tfrac{1}{2}+\lambda_1+\tfrac{1}{2})e_1=(\lambda_1+1)e_1$ and $\chi^-\otimes\tau_{\mu}$ has highest
weight $(-\tfrac{1}{2}+\lambda_1+\tfrac{1}{2})e_1=\lambda_1 e_1$. Theorem 9.16 of \cite{Kna} shows
that the highest weights of the irreducible
subrepresentations of $\pi_{\nu}\arrowvert_H$ are exactly of the form
$\eta=\eta_1 e_1$ where $\vert\eta_1\vert\le\nu_1$. So in order that $[\tilde{\pi}_{\mu}^+:\pi_{\nu}]\ne 0$ we
must have that $\vert\lambda_1+1\vert\le\nu_1$ and in order that
$[\tilde{\pi}_{\mu}^-:\pi_{\nu}]\ne 0$ we
must have that $\vert\lambda_1\vert\le\nu_1$. Using (\ref{nu1}) we therefore see that
\begin{eqnarray*}
\lbrack\tilde{\pi}_{\mu}^+:\pi_{\nu}\rbrack&=&0\quad\textrm{unless
  possibly when }\lambda_1\le -1\textrm{ and }\nu_1=-(\lambda_1+1)\\
\lbrack\tilde{\pi}_{\mu}^-:\pi_{\nu}\rbrack&=&0\quad\textrm{unless
  possibly when }\lambda_1\ge0\textrm{ and }\nu_1=\lambda_1.
\end{eqnarray*}
We conclude that if $\lambda_1\ge 0$, then
\begin{eqnarray*}
\lbrack\tilde{\pi}_{\mu}^+:\pi_{\nu}\rbrack&=&0\quad\textrm{for all }\nu\\
\lbrack\tilde{\pi}_{\mu}^-:\pi_{\nu}\rbrack&=&0\quad\textrm{for }\nu\ne\lambda,
\end{eqnarray*}
and if $\lambda_1\le -1$, then
\begin{eqnarray*}
\lbrack\tilde{\pi}_{\mu}^+:\pi_{\nu}\rbrack&=&0\quad\textrm{for }\nu\ne -(\lambda+e_1)\\
\lbrack\tilde{\pi}_{\mu}^-:\pi_{\nu}\rbrack&=&0\quad\textrm{for all }\nu.
\end{eqnarray*}

Observe that
\[
\begin{array}{ll}1(\lambda+\delta)-\delta=\lambda\in\mathcal{F}_G&\textrm{for
    }\lambda_1\ge 0\\
s_{e_1}(\lambda+\delta)-\delta=-(\lambda+e_1)\in\mathcal{F}_G&\textrm{for
    }\lambda_1\le -1\end{array}.
\]
So proposition \ref{prop:trace} shows that
\[
\trace\tilde{\pi}_{\mu}^+-\trace\tilde{\pi}_{\mu}^-=\left\{\begin{array}{ll}-\pi_{\lambda}&\textrm{for
      }\lambda_1\ge 0\\
\pi_{-(\lambda+e_1)}&\textrm{for }\lambda_1\le -1\end{array}\right..
\]
We conclude that
\[
\begin{array}{ll}
\tilde{\pi}_{\mu}^-=\pi_{\lambda}\textrm{ and }\tilde{\pi}_{\mu}^+=0&\textrm{for
  }\lambda_1\ge 0\\
\tilde{\pi}_{\mu}^+=\pi_{-(\lambda+e_1)}\textrm{ and }\tilde{\pi}_{\mu}^-=0&\textrm{for
  }\lambda_1\le -1\end{array}.
\]
\qed

\subsection{$\SO{2m+1} / \SO{2m}$}
As in example \ref{ex:2m+1} let
$G=\SO{2m+1}$ and $H=\SO{2m}$ where $m\ge 2$. We have that $G_1=\SP{2m+1}$
and $H_1=\SP{2m}$ are the universal covers of $G$ and $H$
respectively. As in the case of
$\SO{3}/\SO{2}$, we use theorem 9.16 of \cite{Kna} to study the
representations $\tilde{\pi}_{\mu}^{\pm}$ for $\mu\in\mathcal{F}_{H_1}'$. Using
the notation of example \ref{ex:2m+1}, we see that 
\[
\delta=\tfrac{1}{2}\sum_{k=1}^m(m-k+1)e_k,\quad\delta_{\pg}=\tfrac{1}{2}\sum_{k=1}^m
e_k.
\]
The simple roots of $G$ are
\[
e_1-e_2,e_2-e_3,\ldots,e_{m-1}-e_m,e_m
\]
and the simple roots of $H$ are
\[
e_1-e_2,e_2-e_3,\ldots,e_{m-1}-e_m,e_{m-1}+e_m.
\]
(See \cite{Kna} section 2.5). We now find the analytically integral
forms $\mathcal{F}$ and $\mathcal{F}_1$.  
Since $G_1$ is compact, semisimple and has center $\{\pm 1\}$,
proposition 4.68 of \cite{Kna} shows that the analytically integral
forms $\mathcal{F}_1$ are given by
\[
\mathcal{F}_1=\{\sum_{k=1}^m\lambda_ke_k\mid\lambda_k\in\ZZ\textrm{ for
all }k\textrm{ or }\lambda_k\in\ZZ+\tfrac{1}{2}\textrm{ for
all }k\}.
\]
Hence proposition 4.67 of \cite{Kna} shows that the analytically
integral forms $\mathcal{F}$ are given by
\[
\mathcal{F}=\{\sum_{k=1}^m\lambda_ke_k\mid\lambda_k\in\ZZ\textrm{ for
all }k\}.
\]
When checking
dominance of an analytically integral form, it is enough to check the
dominance with respect to simple roots.
An element of $\sum_{k=1}^m\nu_ke_k\in\mathcal{F}_1$ is dominant with respect to $G$ exactly when
\[
\nu_1\ge\nu_2\ge\cdots\ge\nu_m\ge 0.
\]
An element $\mu=\sum_{k=1}^m \mu_k e_k\in\mathcal{F}_1$  is dominant with respect to $H$ exactly when
\[
\mu_1\ge\mu_2\ge\cdots\ge\vert\mu_m\vert.
\]
Hence
\[
\mathcal{F}_{H_1}'=\{\sum_{k=1}^m \mu_k e_k\mid
\mu_k\in\ZZ+\tfrac{1}{2}\textrm{ for all }k\textrm{ and }\mu_1-\tfrac{1}{2}\ge\mu_2-\tfrac{1}{2}\ge\cdots\ge\vert\mu_m-\tfrac{1}{2}\vert\}
\]
We now prove the following.
\begin{proposition}
Let $G=\SO{2m+1}$ and $H=\SO{2m}$ where $m\ge 2$. Let
$\mu\in\mathcal{F}_{H_1}'$ and let 
$\lambda=\mu-\delta_{\pg}=\sum_{k=1}^m\lambda_k e_k$ and $\lambda'=\sum_{k=1}^{m-1}\lambda_k e_k-(\lambda_m+1)e_m$.

If $m$ is even
\[
\begin{array}{ll}
\tilde{\pi}_{\mu}^+=\pi_{\lambda}\textrm{ and }\tilde{\pi}_{\mu}^-=0&\textrm{for
  }\lambda_m\ge 0\\
\tilde{\pi}_{\mu}^-=\pi_{\lambda'}\textrm{ and }\tilde{\pi}_{\mu}^+=0&\textrm{for
  }\lambda_m\le -1\end{array}
\]

and if  $m$ is odd
\[
\begin{array}{ll}
\tilde{\pi}_{\mu}^-=\pi_{\lambda}\textrm{ and }\tilde{\pi}_{\mu}^+=0&\textrm{for
  }\lambda_m\ge 0\\
\tilde{\pi}_{\mu}^+=\pi_{\lambda'}\textrm{ and }\tilde{\pi}_{\mu}^-=0&\textrm{for
  }\lambda_m\le -1\end{array}.
\]
\end{proposition}
\proof 
Let $(V_{\mu},\tau_{\mu})$ be the irreducible
representation of $H_1=\SP{n}$ with highest weight 
\[
\mu=\lambda+\delta_{\pg}=\sum_{k=1}^m (\lambda_k+\tfrac{1}{2})e_k.
\]
The half spinor representations $\chi^{\pm}$ of $H_1$ have weights
\[
\tfrac{1}{2}\sum_{k=1}^m\varepsilon_{k}e_k\textrm{ where
  }\varepsilon=(\varepsilon_1,\ldots,\varepsilon_m)\in\{\pm 1\}^m
\]
and such a weight is a weight of $\chi^+$ exactly when
$\varepsilon_k=-1$ for an even number of $\varepsilon_k$. Proposition
9.72 of \cite{Kna} shows that the highest weights of the irreducible parts of
$\chi^{\pm}\otimes\tau_{\mu}$ have the form
\[
\eta=\sum_{k=1}^m(\lambda_k+\tfrac{1}{2}(1+\varepsilon_{k}))e_k\textrm{ where
  }\varepsilon=(\varepsilon_1,\ldots,\varepsilon_m)\in\{\pm 1\}^m
\]
with an even number of $\varepsilon_k=-1$ if $\eta$ is a weight
of $\chi^+\otimes\tau_{\mu}$ and an odd number of $\varepsilon_k=-1$ if $\eta$ is a weight
of $\chi^-\otimes\tau_{\mu}$.
Note that
\[
\eta_k=\left\{\begin{array}{ll}\lambda_k&\textrm{if }\varepsilon_k=-1\\
\lambda_k+1&\textrm{if }\varepsilon_k=1\end{array}\right..
\]
Let $\pi_{\nu}$ be an irreducible
representation of $G$ with highest weight $\nu$, i.e. $\nu=\sum_{k=1}^m\nu_k
e_k\in\mathcal{F}_G$ where $\nu_k\in\ZZ$ for all $k$ and
$\nu_1\ge\nu_2\ge\cdots\ge\nu_m\ge 0$. Suppose that $\pi_{\nu}$ is a subrepresentation of either
$\tilde{\pi}^+_{\mu}$ or $\tilde{\pi}^-_{\mu}$. Then
corollary
\ref{cor:criterion} shows that 
\begin{equation}\label{equal}
\sum_{k=1}^m (\nu_k^2+(m-k+1)\nu_k)=\sum_{k=1}^m (\lambda_k^2+(m-k+1)\lambda_k).
\end{equation}
The irreducible
part $\gamma_{\eta^{\pm}}$ of $\chi^{\pm}\otimes\tau_{\mu}$ with highest
weight $\eta^{\pm}$ induces a subrepresentation of
$\tilde{l}^{\pm}_{\mu}$. Now let $\tilde{\gamma}_{\eta^{\pm}}$
be the restriction of this subrepresentation to the kernel of
$D^{\pm}_{\mu}$. Note that
\[
[\tilde{\pi}^{+}_{\mu}:\pi_{\nu}]=\sum_{\eta^{+}}[\tilde{\gamma}_{\eta^{+}}:\pi_{\nu}],\quad[\tilde{\pi}^{-}_{\mu}:\pi_{\nu}]=\sum_{\eta^{-}}[\tilde{\gamma}_{\eta^{-}}:\pi_{\nu}].
\]
Frobenius reciprocity now shows that
\[
[\tilde{\gamma}_{\eta^{\pm}}:\pi_{\nu}]\le[\pi_{\nu}\arrowvert_H:\gamma_{\eta^{\pm}}].
\]
Theorem 9.16 of \cite{Kna} shows that the irreducible parts of
$\pi_{\nu}\arrowvert_H$ have multiplicity one and are exactly the
representations of $H$ with highest weights $\sum_{k=1}^ma_ke_k$
satisfying
\[
\nu_1\ge a_1\ge\nu_2\ge a_2\ge\cdots\ge\nu_m\ge\vert a_m\vert.
\] 
So in order for $\pi_{\nu}$ to be a subrepresentation of $\tilde{\pi}^+_{\mu}$, we must
have that
\begin{equation}\label{vurd+}
\nu_1\ge\lambda_1+\tfrac{1}{2}(1+\varepsilon_{1})\ge\cdots\ge\nu_m\ge\vert\lambda_m+\tfrac{1}{2}(1+\varepsilon_{m})\vert
\end{equation}
for some $\varepsilon$ with an even number of $\varepsilon_k=-1$ and in order for $\pi_{\nu}$ to be a subrepresentation of $\tilde{\pi}^-_{\mu}$, we must
have that
\begin{equation}\label{vurd-}
\nu_1\ge\lambda_1+\tfrac{1}{2}(1+\varepsilon_{1})\ge\cdots\ge\nu_m\ge\vert\lambda_m+\tfrac{1}{2}(1+\varepsilon_{m})\vert
\end{equation}
for some $\varepsilon$ with an odd number of $\varepsilon_k=-1$.

Suppose that $\lambda_m\ge 0$, (i.e.,  $\lambda\in\mathcal{F}_G$). For
any $\nu$ satisfying 
(\ref{vurd+}) or (\ref{vurd-}), we have in particular that
\[
\nu_k\ge\lambda_k\ge 0\quad\textrm{for }k\in\{1,\ldots,m\}.
\]
When $\nu$ also satisfies (\ref{equal}), then (\ref{vurd+}) or (\ref{vurd-}) can
only hold if
\[
\varepsilon_1=\cdots=\varepsilon_m=-1\textrm{ and
  }\nu_k=\lambda_k\textrm{ for }k\in\{1,\ldots,m\}.
\]
So we have that
\begin{displaymath}
  \begin{array}{ll}
    [\tilde{\gamma}_{\eta^{\pm}}:\pi_{\nu}]=0&\textrm{for }\nu\ne\lambda\\{}
    [\tilde{\gamma}_{\eta^{\pm}}:\pi_{\lambda}]=0&\textrm{unless possibly for one
    }\eta^+\textrm{ if }m\textrm{ is even}\\{}
    [\tilde{\gamma}_{\eta^{\pm}}:\pi_{\lambda}]=0&\textrm{unless possibly for one
    }\eta^-\textrm{ if }m\textrm{ is odd}.
  \end{array}
\end{displaymath}
Hence
\[
\begin{array}{ll}
[\tilde{\pi}_{\mu}^{\pm}:\pi_{\nu}]=0&\textrm{for }\nu\ne\lambda\\{}
[\tilde{\pi}_{\mu}^+:\pi_{\lambda}]\le 1,[\tilde{\pi}_{\mu}^-:\pi_{\lambda}]=0&\textrm{if }m\textrm{ is even}\\{}
[\tilde{\pi}_{\mu}^-:\pi_{\lambda}]\le 1,[\tilde{\pi}_{\mu}^+:\pi_{\lambda}]=0&\textrm{if }m\textrm{ is odd}.
\end{array}
\] 
Now using proposition \ref{prop:trace}, we see that since
$\lambda=1(\lambda+\delta)-\delta\in\mathcal{F}_G$, the representations
$\tilde{\pi}_{\mu}^{\pm}$ cannot both be zero. We conclude that
\[
\begin{array}{ll}
\tilde{\pi}_{\mu}^+=\pi_{\lambda}\textrm{ and }\tilde{\pi}_{\mu}^-=0&\textrm{if
  }m\textrm{ is even}\\
\tilde{\pi}_{\mu}^-=\pi_{\lambda}\textrm{ and }\tilde{\pi}_{\mu}^+=0&\textrm{if
  }m\textrm{ is odd}
\end{array}.
\]

Suppose that $\lambda_m\le -1$. If $\nu$ satisfies 
(\ref{vurd+}) or (\ref{vurd-}), we have in particular that
\begin{equation}\label{vurd}
\nu_k\ge\lambda_k\ge 0\textrm{ for }k\in\{1,\ldots,m-1\},\quad\nu_m\ge\vert\lambda_m\vert-1=-(\lambda_m+1).
\end{equation}
The condition (\ref{equal}) is given by
\begin{equation}\label{sumnulambda}
\sum_{k=1}^{m-1}\left((\nu_k^2-\lambda_k^2)+(m-k+1)(\nu_k-\lambda_k)\right)+
\left((\nu_m^2-\lambda_m^2)+(\nu_m-\lambda_m)\right)=0.
\end{equation}
Now (\ref{vurd}) implies that the first term of (\ref{sumnulambda}) is
non-negative and therefore that the second term is non-positive. The
second term is a second order polynomial in $\nu_m$ with zeroes at
$\lambda_m$ and $-(\lambda_m+1)$ so it is non-positive if $\nu_m$ lies
in the interval $[\lambda_m,-(\lambda_m+1)]$. Since also
$\nu_m\in\{-(\lambda_m+1),-\lambda_m,\ldots\}$ we conclude that for (\ref{vurd+}) or (\ref{vurd-}) and (\ref{equal}) to hold we must have
that
\[
\varepsilon_k=-1, \nu_k=\lambda_k\textrm{ for
  }k\in\{1,\ldots,m-1\}\quad \textrm{and}\quad\varepsilon_m=1,\nu_m=-(\lambda_m+1).
\]
Hence if $\lambda'=\sum_{k=1}^{m-1}\lambda_ke_k-(\lambda_m+1)e_m$, we
have that
\[
  \begin{array}{ll}
    [\tilde{\gamma}_{\eta^{\pm}}:\pi_{\nu}]=0&\textrm{for }\nu\ne\lambda'\\{}
    [\tilde{\gamma}_{\eta^{\pm}}:\pi_{\lambda'}]=0&\textrm{unless possibly for one
    }\eta^-\textrm{ if }m\textrm{ is even}\\{}
    [\tilde{\gamma}_{\eta^{\pm}}:\pi_{\lambda'}]=0&\textrm{unless possibly for one
    }\eta^+\textrm{ if }m\textrm{ is odd}
  \end{array}
\]
and
\[
\begin{array}{ll}
[\tilde{\pi}_{\mu}^{\pm}:\pi_{\nu}]=0&\textrm{for }\nu\ne\lambda'\\{}
[\tilde{\pi}_{\mu}^-:\pi_{\lambda'}]\le 1,[\tilde{\pi}_{\mu}^-:\pi_{\lambda}]=0&\textrm{if }m\textrm{ is even}\\{}
[\tilde{\pi}_{\mu}^+:\pi_{\lambda'}]\le 1,[\tilde{\pi}_{\mu}^+:\pi_{\lambda}]=0&\textrm{if }m\textrm{ is odd}.
\end{array}
\] 
Now using proposition \ref{prop:trace}, we see that since
$\lambda'=s_{e_m}(\lambda+\delta)-\delta\in\mathcal{F}_G$, the
representations 
$\tilde{\pi}_{\mu}^{\pm}$ cannot both be zero. We conclude that
\[
\begin{array}{ll}
\tilde{\pi}_{\mu}^-=\pi_{\lambda'}\textrm{ and }\tilde{\pi}_{\mu}^+=0&\textrm{if
  }m\textrm{ is even}\\
\tilde{\pi}_{\mu}^+=\pi_{\lambda'}\textrm{ and }\tilde{\pi}_{\mu}^-=0&\textrm{if
  }m\textrm{ is odd}
\end{array}.
\]
\qed
\pagebreak

\end{document}